%% file: 550.tex


\ifx\shlhetal\undefinedcontrolsequence\let\shlhetal\relax\fi

\input amstex
\input mathdefs \input citeadd

\sectno=-1   
\localtags
\NoBlackBoxes
\define\mr{\medskip\roster}
\define\sn{\smallskip\noindent}
\define\mn{\medskip\noindent}
\define\bn{\bigskip\noindent}
\define\ub{\underbar}
\define\wilog{\text{without loss of generality}}
\define\ermn{\endroster\medskip\noindent}

\define\dbcu{\dsize\bigcup}
\define\nl{\newline}
\documentstyle {amsppt}
\newbox\noforkbox \newdimen\forklinewidth
\forklinewidth=0.3pt   
\setbox0\hbox{$\textstyle\bigcup$}
\setbox1\hbox to \wd0{\hfil\vrule width \forklinewidth depth \dp0
                        height \ht0 \hfil}
\wd1=0 cm
\setbox\noforkbox\hbox{\box1\box0\relax}
\def\unionstick{\mathop{\copy\noforkbox}\limits}
\def\nonfork#1#2_#3{#1\unionstick_{\textstyle #3}#2}
\def\nonforkin#1#2_#3^#4{#1\unionstick_{\textstyle #3}^{\textstyle #4}#2}     
%
\setbox0\hbox{$\textstyle\bigcup$}
\setbox1\hbox to \wd0{\hfil{\sl /\/}\hfil}
\setbox2\hbox to \wd0{\hfil\vrule height \ht0 depth \dp0 width
                                \forklinewidth\hfil}
\wd1=0cm
\wd2=0cm
\newbox\doesforkbox
\setbox\doesforkbox\hbox{\box1\box0\relax}
\def\nunionstick{\mathop{\copy\doesforkbox}\limits}

\def\fork#1#2_#3{#1\nunionstick_{\textstyle #3}#2}
\def\forkin#1#2_#3^#4{#1\nunionstick_{\textstyle #3}^{\textstyle #4}#2}     
\topmatter
\title {0.1 Laws \\
Sh550} \endtitle
\author {Saharon Shelah \thanks{\null\newline
We thank Alice Leonhardt for the beautiful typing.
First version - Fall '91 \newline
Written 94/April/7 \newline
Partially supported by the United States-Israel Binational Science
Foundation and by the NSF Grant [Keisler?] NSF Grant 144-EF67 \null\newline
Publ.550 \newline
Latest Revision 98/Apr/10} \endthanks} \endauthor
\affil {Institute of Mathematics \\
The Hebrew University \\
Jerusalem, Israel
\medskip
Rutgers University \\
Department of Mathematics \\
New Brunswick, NJ USA \\
\medskip
University of Wisconsin \\
Department of Mathematics \\
Madison, Wisconsin USA} \endaffil
\bigskip

\abstract   We give a framework for dealing with 0-1 laws (for first order
logic) such that expanding by further random structure tend to give us another
case of the framework.  From another perspective we deal with 0-1 laws when 
the number of solutions of first order formulas with parameters behave
dichotomically. \endabstract
\endtopmatter
\document  
\input alice2jlem

\newpage

\head {Anotated Content} \endhead \bigskip

\noindent
\S0 $\quad$ Introduction
\bigskip

\noindent
\S1 $\quad$ The context: probability and model theory
\roster
\item "{${}$}"  [We define the context: random models ${\Cal M}_n$ for which
we may ask if the 0-1 law holds.  Having in mind cases the random
structure there are relations which were drawn with constant probability
as well as one decaying with $n$. We give sufficient condition for the
satisfaction of the 0-1 laws, semi-nice].
\endroster
\bigskip

\noindent
\S2 $\quad$ More accurate measure and drawing again
\roster
\item "{${}$}"  [We restrict ourselves further by fixing, when 
$A <_{\text{pr}} B$, the number of copies of $B$ over a copy of $A$ in
${\Cal M}_n$, up to some factor. The clearest case is when this approximate
value is $\sim \|{\Cal M}_n\|^{\alpha(A,B)}$ (the polynomial case). 
We also give the framework for redrawing].
\endroster
\bigskip

\noindent
\S3 $\quad$ Regaining the context for ${\frak K}^+$: reduction
\roster
\item "{${}$}"  [We deal with ${\frak K}^+$; i.e. having a 0-1 context
${\frak K}$ we expand ${\Cal M}_n$ randomly to ${\Cal M}^+_n$.  Assume in
random enough ${\Cal M}^+_n$, all relevant cases behave close enough to the
expected value, we can see what should be $\le^+_a,\le^+_i,\le^+_s,
\le^+_{\text{pr}}$, so we define such relations called for the time being 
$\le^+_b,\le^+_j,\le^+_t,\le^+_{\text{qr}}$ respectively and 
investigate them, proving they behave
similarly enough to what we hope to prove for ${\frak K}^+$.  Then we
restrict ourselves to the polynomial case and phrase a sufficient
condition for succeeding: $\le^+_b = \le^+_a$, etc. and for preservation
of semi-niceness.  We point out two extreme cases: when the probabilities are
essentially constant and when the probabilities are of the form
$\|{\Cal M}_n\|^\alpha,\alpha \in (0,1)_{\Bbb R}$.  The restriction to the
polynomial case is for simplicity.]
\endroster
\bigskip

\noindent
\S4 $\quad$ Clarifying the probability problem
\mr
\item "{{}}"  [We deal more explicitly with what is needed for showing that
the 0-1 content we get by redrawing is again nice enough.  We restrict
ourselves to what is needed.]
\endroster
\bigskip

\noindent
\S5 $\quad$ The probability argument
\mr
\item "{{}}"  [We replace the content of \S4 by a slightly more general one
and then prove the required inequality.]
\endroster
\bigskip

\noindent
\S6 $\quad$ Free amalgamation
\mr
\item "{{}}"  [We axiomatize the ``edgeless disjoint amalgamation" used in
earlier works to ``free amalgamation".]
\endroster
\bigskip

\noindent
\S7 $\quad$ Variant of niceness
\mr
\item "{{}}"  [We consider some variants of semi-nice (and semi-good) and
their relations.  In earlier verion \S7, \S8 we done inside \S1, \S2.]
\endroster
\newpage

\head {\S0 Introduction} \endhead  \resetall
\bigskip

Here we continue Shelah Spencer \cite{ShSp:304}, 
Baldwin Shelah \cite{BlSh:528}, and (in the model theory) Shelah 
\cite{Sh:467} (see earlier \cite{GKLT}, \cite{Fa76}, Lynch \cite{Ly85}, see
the survey \cite{Sp}); we are
trying to get results of reasonable generality.  In particular, we want:
if our random model ${\Cal M}^+$ expands a given ${\Cal M}$, where the 
class of ${\Cal M}$'s is in our context, then also the class of 
${\Cal M}^+$'s is in our context.  We shall be most detailed on the ``nice 
polynomial" case (see definitions in the text).

Let us turn to trying to explain the results. \newline
${\frak K}$ is a 0-1-context (see Definition \scite{2.1}) give for each $n$, 
a probability distribution for the $n$-th random model ${\Cal M}_n$.  
The paper is self contained and describes a
reasonable result of drawing (finite) models, where a finite sequence
$\bar a$ has a closure, $c \ell^k(\bar a,{\Cal M}_n$) in the model 
${\Cal M}_n$, presumably random enough relative to $k + \ell g(\bar a)$, 
and the number of elements of $c \ell^m(\bar a,{\Cal M})$ has an a priori 
bound; i.e. depending on $\ell g (\bar a)$ and $m$ only.
And as far as first order formulas are concerned, this is all there is (see
\S1).  However, here we allow relations with constant probability; this 
implies that the limit theory are not necessarily stable (as proved in
\cite{BlSh:528} on the theories from \cite{ShSp:304} and more), but they are
simple (\cite{Sh:93}).)  Model theoretically our approach allows 
non-symmetric relations $R(\bar x)$; i.e.
$\neg(R(\dotsc,x_i,\ldots)_{i=1,n} \equiv R(\dotsc,x_{\sigma(i)},
\ldots)_{i=1,n})$ for permutation $\sigma$ of $\{1,\dotsc,n\}$.

An extreme case is when for some $m = m_{\ell g(\bar a)},k > m \Rightarrow
c \ell^k(\bar a,M) = c \ell^m(\bar a,M)$ (when $M$ is random enough relative 
to $\ell g(\bar a) + k$).  Still a very nice case is when for every $m$ and 
$k$ for some $\ell$, for $M$ random enough and $\bar a \in {}^k M$ we have: if
$A_0 = \text{ Rang}(\bar a)$ and $A_{i + 1} = c \ell^m(A_i,M)$ then
$A_\ell = A_{\ell + 1}$ can be proved for $\ell = \ell(k,\ell g(\bar a))$.  
But we may have a successor function to begin with and then this is 
impossible.
\medskip

Note that if there are $\sim \| {\Cal M}_n\|$ uniformly definable sets with
$\sim \sqrt{\text{ log} \| {\Cal M}_n \|}$ elements and we draw a two-place 
relation with the probability of each pair being e.g. $1/2$
we get example of all two-place relations on such definable small sets (up to 
isomorphism), not what we want.   So it is natural to ask that definable sets
$\varphi(x,\bar a)$ with no apriori bound, are quite large, say of size
$\in [h^d_{\varphi,\bar a}({\Cal M}_n),h^u_{\varphi,\bar a}
({\Cal M}_n)]$ (with $d$ for down, $u$ for up.  The size of
$h^d_{\varphi,\bar a}$ is discussed below.  Actually instead of dealing with 
such $\varphi$'s, we deal with $B$'s with $A <_s B$; for each copy of $A$ in 
${\Cal M}_n$ we look at the number of copies of $B$ above it.  In addition, 
concerning the drawing we can weaken the natural demand of independency to:
if $A \subseteq {\Cal M}_n$ ($|A|$ small compared to $n$) the probability of
${\Cal M}^+_n \restriction A = A^+$ is
$\in \left[ p^d_{A / \cong}[{\Cal M}_n],p^u_{A / \cong}[{\Cal M}_n] \right]$ 
even knowing
${\Cal M}^+_n \restriction B$ for every $B \subseteq {\Cal M}_n,B \nsubseteq
A$ (i.e. the conditional probability).  Again the probabilities should be
such that the phenomena described in the first sentence of this paragraph 
do occur.
\medskip

We can ask more: the weak independence described in the last paragraph is 
even knowing all other instances of relations.  We also can ask for true 
independence.

We may look at more strict cases:
\demo{Case 1 - Polynomial Case}\footnote{note: here $n^\alpha,\alpha \in 
\Bbb R$, is considered polynomial, $2^{\sqrt{\log n}}$ is not, so not only
integer powers are considered polynomial}
$h^d_{A/\cong}[{\Cal M}_n]$ and
$h^u_{A/\cong}[{\Cal M}_n]$ are $\cong \| {\Cal M}_n \|^{\alpha(A/\cong)}$ 
(for \newline
$A \in {\Cal K}_\infty$) and $p^d_{A/\cong}[{\Cal M}_n]$ and
$p^u_{A/\cong}[{\Cal M}_n]$ are $\cong \|{\Cal M}_n\|^{\beta(A/\cong)}$ where
$A \in {\Cal K}^\oplus_\infty$. \newline
\medskip

\noindent
In this case the main danger in drawing is that some $\alpha(A)$'s and
$\beta(A)$'s are linearly dependent (looking at $\Bbb R$ as a vector 
space over $\Bbb Q$) similar to [ShSp304]; i.e. for having a $0-1$ law (or
convergence) we need the so-called irrationality condition.

Note that here the resulting framework also falls under Case 1.
\newline
The reader can concentrate on Case 1, and we usually do it here.  
\enddemo
\bigskip

\demo{Case 2}  The $h$'s are near $1$.
\smallskip

\noindent
Here $h^d_{A/\cong}[{\Cal M}_n],h^u_{A/\cong}[{\Cal M}_n]$ are constants
$\in (0,1)_{\Bbb R}$ or at least are in \newline
$({\frac 1{h({\Cal M}_n)}},
1-{\frac 1{h({\Cal M}_n)}})$ where $h$ goes to infinity more slowly than
any $n^\varepsilon$ (where \newline
$\varepsilon \in \Bbb R^{> 0}$).
\medskip

Here we may allow the $p^d_{A/\cong}[{\Cal M}_n],p^u_{A/\cong}[{\Cal M}_n]$
more freedom - just to be large enough (and/or small enough than $1$)
compared to log$(\| {\Cal M}_n \|$).  Here the irrationality condition 
includes requiring that
the $p^d_{A/\cong}[{\Cal M}_n],p^u_{A/\cong}[{\Cal M}_n]$ are close enough
so that we get definite answers. 
\enddemo
\bigskip

\demo{Case 3}  The $p$'s are near $1$.

Here $p^d_{A/\cong}[{\Cal M}_n],p^u_{A/\cong}[{\Cal M}_n]$ are constant
$\in (0,1)_{\Bbb R}$ or at least are in \newline
$({\frac 1{h({\Cal M}_n)}},
1-{\frac 1{h({\Cal M}_n)}})$ where $h$ goes to infinity more slowly than any
$n^\varepsilon$ (where \newline
$\varepsilon \in \Bbb R^{>0}$).

Here we may allow the $h^d_{A/\cong}[{\Cal M}_n],h^u_{A/\cong}[{\Cal M}_n]$
more freedom - just the fraction of $\|{\Cal M}_n\|^{|A|}$ they are far 
enough from $0$ and $1$ compared to log$(\|M_n\|$.  The irrationality 
condition will include requirements that the 
$h^d_{A/\cong}[{\Cal M}_n],h^u_{A/\cong}[{\Cal M}_n]$ are close enough to 
get definite answers.
\medskip

To carry the probability argument (and get that the resulting class
${\frak K}^+$ is also in our framework), we have to deal with the following 
problem.  \newline
For ${\Cal M}_n$ large enough it is natural to define a tree $T$ (i.e. 
a set $T$ and a partial order $< = <_T$ such that $x \in T \Rightarrow 
\{ y:y <_T x\}$ is linearly ordered by
$<_T$ and have lev$(x) \in \Bbb N$ elements).  Let $T_\ell = \{ x \in T:
\text{lev}(x) = \ell \}$, assume $T$ has $m$ levels, and for $\ell < m$, and
$\eta \in T_\ell$, the number $k_\eta = |\text{Suc}_T(\eta)|$ is 
$\in [k^d_\ell,k^u_\ell]$.  Suppose further we
have $A_0 <_{pr} A_1 <_{pr} \cdots <_{pr} A_m$ (see Definition 
\scite{2.3}(2)(d),
$A_\ell$ is a ``small model" in the original vocabulary, increasing with
$\ell$) and we have $\langle f_\eta:\eta \in T \rangle$ such that
$f_\eta$ embeds $A_{\ell g(\eta)}$ into ${\Cal M}_n$, and
$\eta <_T \nu \in T \Rightarrow f_\eta \subseteq f_\nu$.  
Assume further $A^+_i$ is an expansion of $A_i$ to the larger vocabulary,
increasing with $i$, so formally
$A^+_i \le A^+_{i+1},A^+_i \restriction \tau = A_i$. \newline
Lastly assume that for each $c \in {\Cal M}_n$ the set $|\{ \nu \in 
\text{Suc}_T(\eta):c \in \text{ Rang}(f_\nu) \backslash \text{ Rang}
(f_\eta)\}|$ has an apriori bound (i.e. not depending
on $n$).  Under the condition
that $f_{\text{Min}(T)}$ embeds $A^+_0$ into ${\Cal M}^+_n$, what is the
number of $\eta \in T_m$ for which $f_\eta$ embeds $A^+_m$ into
${\Cal M}^{+2}_n$? (at least approximately). \newline

We can easily give bounds to the expected value, and, working harder, for 
variance.  But we want to show that the probability of ``large" deviation 
from the
expected value are $< \| {\Cal M}_n\|^{-\alpha}$ for every $\alpha \in 
\Bbb R^+$.
This is enough to really ignore those cases (not just say they occur very
rarely but show that, for ${\Cal M}_n$ random enough they do not occur at
all).
This is easy if $m = 1$, and so we have, essentially, many independent
events.

Note the following obstruction: the independence is violated as possibly
$\eta_1 \ne \eta_2,\nu_1 \in \text{ Suc}_T(\eta_1),\nu_2 \in \text{ Suc}_T
(\eta_2)$ and Rang$(f_{\nu_1}) \backslash \text{Rang}(f_{\eta_1}),
\text{Rang}(f_{\nu_2}) \backslash \text{Rang}(f_{\eta_2})$ are not disjoint.
Particularly disturbing case is when for some $x \ne y$ from \newline
$A_{\ell g(\nu_1)},A_{\ell g(\nu_2)}$ respectively we have $f_{\nu_1}(x) = 
f_{\nu_2}(y)$. 
However, this chaotic obstruction can be overcome by shrinking somewhat 
the tree, so we can get:
\medskip
\roster
\item "{$\bigoplus$}"  $f_{\eta_1}(a_1) = f_{\eta_2}(a_2) \Rightarrow
a_1 = a_2$.
\endroster
\medskip

\noindent
Still we have two extreme cases which are quite different
\medskip
\roster
\item "{$\bigoplus_1$}"  for each $a_1$ there is $\ell$ such that
$a_1 \in A_{\ell + 1}$ and for every $a_2$: \newline
$f_{\eta_1}(a_1) = f_{\eta_2}(a_2) \Leftrightarrow (a_1 = a_2 \in 
A_{\ell + 1} \and \eta_1 \restriction (\ell + 1) = \eta_2 \restriction
(\ell + 1))$
\medskip
\noindent
\item "{$\bigoplus_2$}"  $T = \dsize \bigcup_{\ell \le m} \dsize
\prod_{i < \ell} m_i$ and: if $a \in A_{\ell + 1} \backslash A_\ell \and
\eta_1,\eta_2 \in T_{\ell + 1}$ then: \newline
$\eta_1(\ell) = \eta_2(\ell) 
\Leftrightarrow f_{\eta_1}(a) = f_{\eta_2}(a)$.
\endroster
\medskip

What we need are good upper bounds on deviation uniformly under the
circumstances (including, in particular, those two cases).

We deal with probability only as required: anything with probability
\newline
$< \| {\Cal M}_n \|^c$ (i.e. for every $c \in \Bbb R^+$, for every random
enough ${\Cal M}_n$) can be discarded; also we need to eliminate the
extreme cases ($=$ large deviation) but do not care about the exact
distribution in the middle.

Having explained the probability side problems, let us turn to the model
theoretic ones.  First, we want to include cases like the successor
function so that, possibly there are no set $A \subseteq {\Cal M}_n$ 
satisfying $A \notin \{ \emptyset,|{\Cal M}_n|\}$ with
$c \ell^k(A,{\Cal M}_n) = A$.  Hence we have to suffice ourselves with 
semi-niceness (similarly to \cite{Sh:467}) rather than niceness (as is 
done in \cite{ShSp:304} and more generally in \cite{BlSh:528}).  \newline
Secondly, the ``free amalgamation" is no longer the disjoint amalgamation
with no new edges (or new cases of the relations) as in \cite{ShSp:304},
\cite{BlSh:528}, \cite{Sh:467}.  This change enables us to deal, e.g. with
the random graph structure with two kinds of edges ($=$ two 2-place relations)
red with probability $n^{-\alpha}(0 < \alpha < 1)$ and green with edge 
probability 1/2.  But, as in \cite{ShSp:304},\cite{BlSh:528}
(and unlike \newline
\cite{Sh:467}),
$c \ell^k(A,{\Cal M}_n)$ has, for random enough ${\Cal M}_n$, a bound
depending on $|A|,k$ only.  We could generalize our framework so that
\cite{Sh:467} is included but decide not to do it here.
But there is a price: the definition of ``nice"
(see \scite{2.6C}) becomes more complicated; having both is dealt
with in \cite{Sh:637}.
\medskip

An interesting phenomena is the dichotomy: either the limit theory is very
simple, analyzable or is complicated (and related to \cite{Sh:c}).  Lately
Tyszkiewicz proves related theorem for monadic theories
of for the classes of groups; I think the right starting point should be
the parallel infinite problem with no probability, what can be the 
(monadic) infinitary
theory of models of first order $T$ monadically expanded, see Baldwin
Shelah \cite{BlSh:156}, survey [Baldwin, handbook of model theoretic logics].
I do not know enough to conjecture the dividing line, but can note
that all the ``complicated" limit theories $T_\infty$ are ``complex" which we
define as: there is a formula $\varphi_t(\bar x,\bar y,\bar z,t),
\varphi_x(\bar x,\bar y,\bar z,t)$ such that
$T_\infty \Vdash$ ``for some $t,\varphi_+(\bar x,\bar y,\bar z,t),\psi_\times
(\bar x,\bar y,\bar z,t)$ define a model which satisfies the axioms of PA
including induction but adding ``there is a last element, there are $\ge n$ 
elements" (can add induction scheme for formulas of quantifier depth $\le n$.
See more in Baldwin \cite{B}[nearly model complete and 0-1].  See Baldwin
Shelah \cite{Sh:639}.

We can deal with the example of directed graphs, edge probability
$n^{-\alpha}$ both directions has same probability.
\bigskip

\centerline {$* \qquad * \qquad *$}
\bigskip

\noindent
Baldwin has asked me to explicate here how the theory of $([n],S,R)$ where
$S$ is the successor relation, $R$ a random graph with edge probability
being $n^{-\alpha}$, i.e. how the fact that every element has an immediate
successor, is reflected in the present treatment (that is, got into the 
axiomatization
we get; this is unlike \cite{BlSh:528}).  We can consider three cases: with
successor being modulo $n$, with usual successor and (e.g.) with 
successor being
$(i,i+1)$ for $i$ not divisible by the square root of $n$ (rounded).  The
limit models are: in the first case it has copies of $\Bbb Z$, in the second
case one copy of $\omega$, one copy of $\omega^*$ and copies of $\Bbb Z$, and
in the third case many copies of $\omega$, many copies of $\omega^*$ and
many copies of $\Bbb Z$.  Now in the semi-nice case (see 
Definition \scite{2.6G}) we should look at the set of semi-$(k,r)$-good 
quadruples, now the pairs of $(A,A^+)$ which appear in such 
quadruples (so $A \subseteq A^+$) gives us the distinctions.

For $x \in A$ let 

$$
\align
\ell^+(x,A) = \text{ Min}\{\ell:&\text{ there are no } x_0,\dotsc,x_\ell \in
A \text{ such that } x = x_0 \text{ and} \\
  &A \models S(x_i,x_{i+1}) \text{ for } i < \ell\}
\endalign
$$

$$
\align
\ell^-(x,A) = \text{ Min}\{\ell:&\text{ there are no } x_0,\dotsc,x_\ell \in
A^+ \text{ such that } x = x_\ell \text{ and} \\
  &A \models S(x_i,x_{i+1}) \text{ for } i < \ell\}.
\endalign
$$
\medskip

\noindent
Let $m(k,r)$ be large enough.

Now in the first case, in any such pair for every $x \in A$ we have
$\ell^+(x,A^+) \ge m(k,r)$ and $\ell^-(x,A^+) \ge m(k,r)$, in the third case 
for no $x \in A,\ell^+(x,A^+) < m(k,r) \and \ell^-(x,A^+) < m(k,r)$ and in 
the second case, there may be at most one $x \in A$ with $\ell^+(x,A) <
m(k,r)$ and at most one $x \in A$ with $\ell^-(x,A^+) < m(k,r)$ 
(but they are not the same).

In the strict polynomial (or even less) case we can also deal with properties
suggested by Lynch \cite{Ly90}.  He asks for the results in Shelah
Spencer \cite{ShSp:304} for more accurate numerical (asymptotic) results, 
particularly in the case the probability is zero he proved
\mr
\item "{$(*)_1$}"  for every first order sentence $\psi$ such that
Prob$({\Cal M}_n \models \psi)$ converge to zero one of the following occurs:
{\roster
\itemitem{ $(\alpha)$ }  Prob$({\Cal M}_n \models \psi) = 
c\|{\Cal M}_n\|^{- \beta} + 0 (n^{-\beta-\varepsilon})$ for some 
$c,\beta,\varepsilon \in \Bbb R^+$,
\sn
\itemitem{ $(\beta)$ }  Prob$({\Cal M}_n \models \psi) = 0(\|{\Cal M}_n\|
^\varepsilon)$ for every $\varepsilon \in \Bbb R^+$. \nl
Confirming his conjecture in \cite{Sh:551} we prove
\itemitem{ $(\beta)^+$ }  Prob$({\Cal M}_n \models \psi) = 
0(e^{-\|{\Cal M}_n\|^\varepsilon})$ for every $\varepsilon \in \Bbb R^+$.
\endroster}
\ermn
We shall explicate this elsewhere.

The starting point of this research was a question of Lynch communicated to
me by Spencer in Fall '91 on whether we can do \cite{ShSp:304} starting with
a successor function; but I thought the real problem was to have 
a general framework and I lectured on it in Rutgers Fall '95; see
\cite{BlSh:528}.

See more \cite{Sh:467}, \cite{Sh:581} and [Baldwin, near model complete and
0-1 laws]. \newline
I thank Shmuel Lefschas and John Baldwin for comments and corrections.
Earlier we have used $\nonfork{}{}_{}$ and version of niceness from the
beginning of the paper.
\enddemo
\bigskip

\demo{\stag{1.1} Notation}  $\Bbb R$ set of reals, 
$\Bbb R^{> 0} = \{ \alpha \in 
\Bbb R: \alpha > 0 \},\Bbb R^{\ge 0} = \Bbb R^{> 0} \cup \{ 0 \},\mu_n$ 
the nth probability (= distribution).  
Here $n$ is the index for ${\Cal M}_n$, which is always used 
for the model chosen $\mu_n$ - randomly (we do not assume ${\Cal M}_n$
necessarily has exactly $n$ elements).  $\Bbb N$ is the set of natural
numbers.  We use $k,\ell,m,n,i,j,r,s$ for natural numbers and we use
$\alpha,\beta,\gamma$ for reals, $\varepsilon,\zeta$ for positive reals.

Let $A,B,C,D$ denote finite models $M,N$ models and $f,g$ denote embeddings.
Let $h$ denote a function with range $\subseteq \Bbb R$.
\enddemo
\bigskip

\demo{\stag{1.2} Notation}  1) We use $\tau$ for vocabularies, $\tau$
consisting of predicates only (for simplicity), $n(R)$ the number 
of places of $R$ ($=$ the arity of $R$). \newline
2)  In general treatment we can demand that each $R \in \tau$ 
will be interpreted as irreflexive relation; i.e. $\bar a \in R^M 
\Rightarrow \bar a \text{ without repetition}$; (by this demand we 
do not lose any generality as we can add suitable predicates).
We call such $\tau$ irreflexive, but we do not require symmetry (so directed
graphs are allowed).  \newline
We use $A,B,C,D$ for models which are finite, if not explicitly said
otherwise, and $M,N$ for models; we notationally do not strictly distinguish 
between a model and its universe.  Those are
$\tau$-models and $A^+,\dotsc,N^+$ are $\tau^+$-models if not explicitly said
otherwise.  We use $a,b,c,d,e$ for elements of models, bar signifies a
finite sequence.  \newline
3)  We call $\tau$ locally finite if for every $n$ the set $\{ R \in \tau:
n(R) = n\}$ is finite.  Note:  the number of $\tau$-model with
the finite universe $A$, is finite when $\tau$ is finite or locally finite
irreflexive. \newline
4)  Let $f:A \rightarrow B$ mean both are models with the same vocabulary, and
$f$ is an embedding, i.e. $f$ is one to one and for any predicate $R$ 
(in $\tau$, the vocabulary of $A$ and $B$) which is $k$-place we have:
$a_1,\dotsc,a_k \in A \Rightarrow [\langle a_1,\dotsc,a_k \rangle \in R^A
\leftrightarrow \langle f(a_1),\dotsc,f(a_k) \rangle \in R^B]$.  Let id$_A$ 
be the identity map on $A$.  Let $A \subseteq B$ mean id$_A:A \rightarrow B$
and we say: $A$ is a submodel of $B$. \newline
5)  If $A,B$ are submodels of $C$ then $A \cup B$ means $C \restriction
(A \cup B)$. \newline
6)  We say $A,C$ are freely amalgamated over $B$ in $M$ if $B \subseteq A
\subseteq M,B \subseteq C \subseteq M,A \cap C = B$ and: 
if $R$ is a predicate of $M$, for
no $\bar a \in R$, do we have Rang$(\bar a) \nsubseteq B \cup A$, Rang$(\bar
a) \nsubseteq B \cup C$ and Rang$(\bar a) \subseteq B \cup A \cup C$; we
also say $A,C$ are free over $B$ inside $M$.  (But $\nonfork{}{}_{}$-free 
means according to the definition of $\nonfork{}{}_{}$, but this 
generalization is done only in \S6, \S7).
\bigskip

\centerline {$* \qquad * \qquad *$}
\bigskip

\demo{\stag{1.3} Notation}  1)  If $f$ is an embedding of $A$ into $M$ and $A 
\subseteq B$, we say \newline
$\bar g = \langle g_i:i = 1,k \rangle$ are disjoint extensions of
$f$ (for ($A,B$)) \underbar{if}:
\medskip
\roster
\item "{(a)}"  $g_i$ is an extension of $f$ to an embedding of $B$ into $M$
\item "{(b)}"  $1 \le i < j \le k \Rightarrow \text{ Rang}(g_i) \cap
\text{ Rang}(g_j) = \text{ Rang}(f)$.
\endroster
\medskip

We say $\bar g$ is a disjoint $k$-sequence of extensions of $f$ if the above
holds; we also say: $\bar g$ is a sequence disjoint over $A,\bar g$ of length 
$k$. \newline
2)  ex$(f,A,B,M)$ where $A \subseteq B,f$ an embedding of $A$ into $M$ is the
set of extensions $g$ of $f$ to embedding of $B$ into $M$. \newline
nu$(f,A,B,M)$ is the number of elements in ex$(f,A,B,M)$. \newline
3)  Let $\Bbb N$ and also $\omega$ denote the set of natural numbers.
\enddemo
\newpage
\enddemo
\head {\S1 The Context: probability and model theory} \endhead  \resetall
\bigskip

We start by defining a 0-1 context (in Definition \scite{2.1}), defining
the derived $A \le_i B$ ($B$ algebraic over $A$), $A \le_s B$ (dual),
$c \ell^k(A,M)$, (in Definition \scite{2.3}, \scite{2.4}) and point out the
basic properties (in \scite{2.6}).  We define ``${\frak K}$ is weakly nice"
and state its main property, (see \scite{2.6C}, \scite{2.6E}).  Then we
define our main version of nice, (Definition \scite{2.6G}, semi-nice) 
and investigate some variant (\scite{2.6H}, \scite{2.6I}) and define the 0-1
laws and variants (\scite{2.9}, \scite{2.10}, \scite{2.11}).
We prove that ${\frak K}$ is semi-nice implies elimination of quantifiers 
and phrase what this gives for
0-1 laws (in \scite{2.14}, see Definition \scite{2.9}, \scite{2.10}).
\bigskip

\definition{\stag{2.1} Definition}  1) ${\frak K}$ is a $0-1$ context if it 
consists of $\tau,{\Cal K},\le$ and \newline
$\bar{\Cal K} = \langle {\Cal K}_n,\mu_n:n \in \Bbb N \rangle$
satisfying (a)-(c) below where:
\newline
(a)  $\tau$ a vocabulary consisting of predicates only (for simplicity,
$\tau$ irreflexive, see \scite{1.2}), $\tau$ finite or at least locally finite.
\newline
(b)  ${\Cal K}$ a family of finite $\tau$-models, 
closed under isomorphisms and submodels; we denote members by $A,B,C,D$ (sometimes $M,N$) and
for notational simplicity the empty model belongs to ${\Cal K}$. \newline
(c)  ${\Cal K}_n \subseteq {\Cal K},\mu_n$ is a probability measure on 
${\Cal K}_n,{\Cal M}_n$ varies on ${\Cal K}_n$; for notational simplicity 
assume $n_1 \ne n_2 \Rightarrow 
{\Cal K}_{n_1} \cap {\Cal K}_{n_2} = \emptyset$ and $M \in {\Cal K}_n
\Rightarrow \| M \| > 1$; also we sometime ``forget" the possibility
$M_{n_1} \in {\Cal K}_{n_1} \and M_{n_2} \in K_{n_2} \and n_1 \ne n_2 \and
\| M_{n_1}\| = \|M_{n_2}\|$ but no confusion should arise.
\enddefinition
\bigskip

\definition{\stag{2.2} Definition}  Let: 
``every random enough ${\Cal M}_n$ satisfies $\psi$" mean \newline
$1 = \text{ Lim inf}_n \text{ Prob}_{\mu_n}({\Cal M}_n \models \psi$).
Similarly ``almost surely ${\Cal M}_n \models \psi$" and
``$a.s.\,{\Cal M}_n$ satisfies $\psi$".
\enddefinition
\bigskip

\definition{\stag{2.3} Definition}  For ${\frak K}$ as in \scite{2.1}(1) 
we define:
\newline
1)  ${\Cal K}_\infty = \{A \in {\Cal K}:0 < \text{ Lim sup}_n 
\text{ Prob}_{\mu_n}
(A \text{ is embeddable into } {\Cal M}_n)\}$. \newline
2)  We define some two place relations on ${\Cal K}$ (mostly on 
${\Cal K}_\infty$):
\medskip
\roster
\item "{(a)}"  $A \le B$ \underbar{if} $A \subseteq B$ (being submodels) and
$B \in {\Cal K}_\infty$ (hence $A \in {\Cal K}_\infty$, see \scite{2.6}(1))
\smallskip
\noindent
\item "{(b)}"  $A \le_i B$ \underbar{if} $A \le B \in {\Cal K}_\infty$ and 
for some $m \in \Bbb N$ we have 
$1 = \text{ Lim}_n \text{ Prob}_{\mu_n}$(for 
every embedding $f$ of $A$ into ${\Cal M}_n$, there are at most 
$m$ extensions of $f$ to an embedding of $B$ into ${\Cal M}_n$) \newline
[the intention is: $B$ is algebraic over $A$]
\smallskip
\noindent
\item "{(c)}"  $A \le_s B$ \underbar{if} $A \le B \in {\Cal K}_\infty$ and 
for no $B'$ do we have $A <_i B' \le B$ \newline
[the intention is: ``strongly" in some sense $A$ is very closed inside $B$]
\smallskip
\noindent
\item "{(d)}"  $A <_{pr} B$ \underbar{if} $A,B \in {\Cal K}_\infty$ 
and $A <_s B$ and for no $C$ do we have $A <_s C <_s B$
\smallskip
\noindent
\item "{(e)}"  $\bar A$ is a decomposition of $A <_s B$ if
$\bar A = \langle A_\ell:\ell \le k \rangle$ and \newline
$A = A_0 <_{pr} A_1 < \cdots <_{pr} A_k = B$ \newline
[the intention of pr is ``primitive", cannot be decomposed]
\smallskip
\noindent 
\item "{(f)}"  $A \le_a B$ means $A \in {\Cal K}_\infty,B \in {\Cal K}_\infty,
A \subseteq B$ and $A=B$ or for some $m \in \Bbb N$ we have 
$1 = \text{ Lim}_n \text{ Prob}_{\mu_n}$ (for every
embedding $f$ of $A$ into ${\Cal M}_n$, there is no sequence $\langle
g_\ell:\ell < m \rangle$ of embeddings of $B$ into ${\Cal M}_n$, pairwise
disjoint over $A$) \newline
[the intention is: $B$ is algebraic in a weak sense over $A$; more
accurately $A$ is not strongly a submodel of $A$].
\endroster
\enddefinition
\bigskip

\definition{\stag{2.4} Definition}  1) 
$T_\infty = \{ \varphi:\varphi \text{ a first order sentence such
that }$ \newline

$\qquad \qquad \qquad \qquad \qquad \qquad \quad
1 = \text{ Lim}_n \text{ Prob}_{\mu_n}({\Cal M}_n \models \varphi)\}$.\newline
2)  ${\Cal K}^{\text{lim}} = \{M:M \text{ a model of } T_\infty\}$. \newline
3) If $A \subseteq M \in {\Cal K}$ and $k \in \Bbb N$ let \newline
$c \ell^k(A,M) = \cup\{B:B \cap A \le_i B \subseteq 
M \text{ and } |B| \le k\}$. \newline
4) For $A \subseteq M \in {\Cal K},A \in {\Cal K}_\infty$ and
$k,m \in \Bbb N$ we define $c \ell^{k,m}(A,M)$ by induction on $m$ as follows
$c \ell^{k,0}(A,M) = A$ and $c \ell^{k,m+1}(A,M) = c \ell^k(c \ell^{k,m}
(A,M),M)$.  Also let $c \ell^{k,\infty}(A,M) = \dsize \bigcup_m
c \ell^{k,m}(A,M)$ and $c \ell^\infty(A,M) = \dbcu_k c \ell^k(A,M)$. \nl
5) For any $(A,B,D)$ satisfying $A \le D \in {\Cal K}_\infty,B \le D$ and
$k \in \Bbb N$ and embedding $f:A \rightarrow M$ let

$$
\align
ex^k(f,A,B,D,M) = \biggl\{g:&\,g \text{ is an embedding of } D \text{ into } M
\text{ extending } f \\
  &\text{ such that } c \ell^k(g(B),M) = g(c \ell^k(B,D)) \biggr\}
\endalign
$$

$$
nu^k(f,A,B,D,M) = |ex^k(f,A,B,D,M)| 
$$
\medskip

\noindent
(compare with \scite{1.3}(2)).
\enddefinition
\bigskip

\remark{\stag{2.5} Remark}  We have chosen the present 
definition of $c \ell^k$ so as to have more cases where iterating the 
operation lead shortly to a fix point and to be compatible with
\cite{Sh:467} where the possibility $c \ell^k(A,M) \notin {\Cal K}_\infty$
exists.
\endremark
\bigskip

\proclaim{\stag{2.6} Claim}  1) ${\Cal K}_\infty \subseteq {\Cal K}$ is 
closed under isomorphisms and submodels. \newline
[why? reread Definition \scite{2.3}(1)]. \newline
2) $\le_i$ is a partial order on ${\Cal K}_\infty$; if 
$A \le B \le C,A \le_i C$ \underbar{then} $B \le_i C$.  Also $A \le_i B$
\underbar{iff} $B \in {\Cal K}_\infty,A \subseteq B$ and for some 
$m \in \Bbb N$, for no $D \in {\Cal K}_\infty,A \subseteq D$, is there a 
sequence of $m$ (distinct) embeddings of $B$ into $D$ over $A$.  \newline
[why? reread Definition \scite{2.3}(2)(b)]. \newline
3) If $A_1 \le_i A_2 \le C \in {\Cal K}_\infty$ and $B \le C \in 
{\Cal K}_\infty$ and $A_1 \subseteq B$ \underbar{then} $B \le_i B \cup A_2$ 
\newline
[why? see Definition \scite{2.3}(2)(b)]. \newline
4) If $A \le C$ are in ${\Cal K}_\infty$ \underbar{then} for one and only one 
$B \in K_\infty$ we have $A \le_i B \le_s C$ [why?  let $B \le C$ be maximal 
such that $A \le_i B$, it exists as $A$ satisfies the demand and $C$ is 
finite, now $B \le_s C$ by \scite{2.6}(2) + Definition \scite{2.3}(2)(c).  
Hence at least one $B$ exists, so suppose
$A \le_i B_\ell \le_s C$ for $\ell = 1,2$ and $B_1 \ne B_2$ so without loss
of generality $B_2 \backslash B_1 \ne \emptyset$.  Now by \scite{2.6}(3), 
$B_1 \le_i B_1 \cup B_2$ hence $B_1 <_i B_1 \cup B_2 \le C$, but this 
contradicts $B_1 \le_s C$ (see Definition \scite{2.3}(2)(c))].
\newline
5)  If $A <_s B$ \underbar{then} there is a decomposition $\bar A$ of 
$A <_s B$ \newline
[why? see Definition \scite{2.3}(2)(e) and Definition \scite{2.3}(2)(d), 
remembering $B$ is finite]. \newline
6)  If $A \le_s B$ and $C \le B$ \underbar{then} $C \cap A \le_s C$; (note
also \newline
$A \le_s C \and A \le B \le C \Rightarrow A \le_s B$) \newline
[why? otherwise for some $C'$ we have $C \cap A <_i C' \le C$, 
then by \scite{2.6}(3) we have $A <_i A \cup C'$, contradiction to 
$A \le_s B$.  The second phrase holds by Definition \scite{2.3}(2)(c)].
\newline
7)  The relations $\le_i,\le_s,\le_{pr},\le_a$ are preserved by isomorphisms.
\newline
[Why?  read Definition \scite{2.3}(2)].
\newline
8)  If $A <_{pr} B$ \underbar{then} for every $b \in B \backslash A$ we have
$(A \cup \{b\}) \le_i B$; also $A < C \le B \Rightarrow C \le_i B$ \newline
[why? if not, then by \scite{2.6}(4) for some $C,(A \cup \{a\}) 
\le_i C <_s B$, but  \newline
$A <_{pr} B \Rightarrow A <_s B \Rightarrow A <_s C$ (by Definition 
\scite{2.3}(2)(b) and \scite{2.6}(6) respectively) so $A <_s C <_s B$ 
contradicting $A <_{pr} B$.  The second phrase is proved similarly.] \newline
9)  $A \le_i B$ \underbar{iff} $A \le B$ and for every $A'$ we have
$A \le A' < B \Rightarrow (A' <_a B)$ \newline
[why?  trivially $A <_i B \Rightarrow A <_a B$; now
the implication $\Rightarrow$ by \scite{2.6}(2) second phrase + 
Definition \scite{2.3}(2)(b),(f); the implication $\Leftarrow$ by the 
$\triangle$-system lemma and the definitions]. \newline
10)  $\le_s$ is a partial order on ${\Cal K}_\infty$ \newline
[why? by the definition $A \le_s A$, as $A \le_s B \Rightarrow A \subseteq B$
and clearly $A \le_s B \le_s A \Rightarrow A=B$, so the problem is 
transitivity.
So assume $A \le_s B \le_s C$ but $\neg(A \le_s C)$ and we shall derive a
contradiction.  As $\neg(A \le_s C)$ by Definition \scite{2.3}(2)(c) there 
is $B_1$ such that $A <_i B_1 \le C$.  If $B_1 \subseteq B$ we get 
contradiction to $A \le_s B$ by Definition \scite{2.3}(2)(c).  
If $B_1 \nsubseteq B$, then by \scite{2.6}(3)
we get $B \le_i (B_1 \cup B)$, but as $B_1 \nsubseteq B$ we have
$B <_i (B_1 \cup B)$, but clearly $(B_1 \cup B) \le C$ so we get contradiction
to $B \le_s C$ by Definition \scite{2.3}(2)(c), so in any case we have gotten 
the desired contradiction]. \newline
11) $A <_s B$ if for every $m$, $0 < \text{ Lim sup Prob}_{\mu_n}$ (for some
embedding $f$ of $A$ into ${\Cal M}_n$ there are $m$ disjoint extensions
$g:B \rightarrow {\Cal M}_n$ of the embedding $f$).  
If $A <_{\text{pr}} B$ then the inverse statement holds. \newline
[why? read Definition \scite{2.3}(2)(b),(c), see details in the proof of
\scite{2.6E}(1)]. \newline
12)  If $A <_{\text{pr}} B \le D,A \le C \le D$, \underbar{then} 
$C <_{\text{pr}} B \cup C$ or $C \le_i B \cup C$. \newline
[why? by \scite{2.6}(4) for some $C_1$ we have $C \le_i C_1 \le_s B \cup
C$.  If $C_1 \cap B \ne A$, then $A < C_1 \cap B \le B$ hence by
\scite{2.6}(8) we have $C_1 \cap B \le_i B$ so by \scite{2.6}(3) we have
$C_1 \le_i B \cup C_1$ and, of course, $B \cup C_1 = B \cup C$, so
$C \le_i C_1 \le_i C \cup B$ so by \scite{2.6}(2) we have $C \le_i C \cup B$,
one of the possible conclusions.  So assume $B \cap C_1 = A$ hence $C =
C_1$, so $C \le_s B \cup C$, now if $C = B \cup C$ clearly $C \le_i B \cup
C$.  Hence we assume $C \ne B \cup C$, so $C <_s B \cup
C$ and if $C <_{\text{pr}} B \cup C$ we get one of the possible conclusions 
as above.
So assume $\neg(C <_{\text{pr}} B \cup C)$, necessarily for some 
$C_2,C <_s C_2 <_s B \cup C$.  By
\scite{2.6}(6) we have $C \cap B <_s C_2 \cap B <_s B$ so as 
$C \cap B = C_1 \cap B = A$ clearly $A <_s C_2 \cap B <_s B$ hence we get a 
contradiction finishing the proof.]
\newline
13) For every $\ell,k \in \Bbb N$ there is $m(k,\ell) \in \Bbb N$ such that:

$$
\text{if } A \le B \in {\Cal K}_\infty \text{ and } |A| \le \ell
\text{ then } c \ell^k(A,B) \text{ has } \le m(k,\ell) \text{ elements}.
$$
\noindent
[why? read the definitions noting that (even if $\tau$ is only locally
finite) the number of pairs $(C_1,C_2),C_i \subseteq A,C_1 \le C_2,
|C_2| \le k$ up to isomorphism over $C_1$ has a bound depending only on
$C_1$]. \newline
14) $\quad$(a) \,\, $c \ell^k(A,M) \subseteq M$ and \newline

$\qquad \qquad c \ell^k(A,M) \in {\Cal K}_\infty 
\Rightarrow A \le_i c \ell^k(A,M)$
\newline
\smallskip
$\quad$ (b) \,\, $A \subseteq B \subseteq M \Rightarrow c \ell^k(A,M) 
\subseteq c \ell^k(B,M)$ \newline
\smallskip
$\quad$ (c) \,\, if $c \ell^k(A,M) \subseteq N \subseteq M$ \underbar{then}
$c \ell^k(A,N) = c \ell^k(A,M)$ \newline
\smallskip
$\quad$ (d) \,\, if $A \subseteq N \subseteq M$ then $c \ell^k(A,N) \subseteq
c \ell^k(A,M)$ \newline
\smallskip
$\quad$ (e) \,\, for $k \le \ell$ we have $c \ell^k(A,M) \subseteq
c \ell^\ell(A,M)$ \newline
\smallskip
$\quad$ (f) \,\, for every $A \in {\Cal K}_\infty$ and $k$ for every random
enough ${\Cal M}_n$ and  \newline

$\qquad \quad$ embedding $f:A \rightarrow {\Cal M}_n$ we have \newline

$\qquad \quad {\Cal M}_n \restriction c \ell^k(f(A),{\Cal M}_n) \in 
{\Cal K}_\infty$ \newline
\smallskip
$\quad$ (g) \,\, for every $k,m$ for some $\ell$ for every $A \subseteq M \in
{\Cal K}_\infty$ \newline

$\qquad \quad$ we have 
$c \ell^m(c \ell^k(A,M),M) \subseteq c \ell^\ell(A,M)$. 
\mn
[Why?  Just check.] 
\sn
15) $T_\infty$ is a consistent (first order) theory which has infinite 
models if \newline
$0 < \text{ lim sup Prob}_{\mu_n}(\|{\Cal M}_n\| > k)$ for every $k$.
\endproclaim
\bigskip

\remark{Remark}  Note that not necessarily in \scite{2.6}(11), we have
``iff".  Why?  e.g. if \nl
$\tau = \{P_1,P_2\}$ with $P_1,P_2$ unary predicates,
$B = \{b_1,b_2\} \in K_\infty,A = \emptyset,
P^B_\ell = \{b_\ell\}$, and for every $M \in {\Cal K},|P^M_1| = 0 \vee 
|P^M_2| = 0$ and
$M \in K_n \and n$ even $\Rightarrow |P^M_1| \ge n$ and $M \in K_n 
\and n$ odd $\Rightarrow |P^M_2| \ge n$.
\endremark
\bigskip

\definition{\stag{2.6C} Definition}  For ${\frak K}$ as in 
Definition \scite{2.1}, we say ${\frak K}$ is \underbar{weakly nice} 
if we have:
\medskip
\roster
\item "{$(*)_1$}"  for every $A <_{\text{pr}} B$ 
(from ${\Cal K}_\infty$) and $m \in \Bbb N$ we have \newline 
$1 = \text{ Lim}_n \text{ Prob}_{\mu_n} \biggl($for every embedding $f$
of $A$ into ${\Cal M}_n$ there are $m$ \newline

$\qquad \qquad \qquad \qquad$ disjoint extensions of $f$ to embedding of 
$B$ into ${\Cal M}_n \biggr)$.
\endroster
\enddefinition
\bigskip

\proclaim{\stag{2.6E} Claim}  Assume ${\frak K}$ is weakly nice. \newline
1) If $A < B \in {\Cal K}_\infty$, \underbar{then} the 
following are equivalent:
\medskip
\roster
\item "{(a)}"  $A <_s B$
\item "{(b)}"  for every $m < \omega$ \newline
$1 = \text{ Lim}_n \text{ Prob}_{\mu_n}\biggl($for every embedding $f$ 
of $A$ into ${\Cal M}_n$ there are \newline

$\qquad \qquad \qquad \qquad$ embeddings $g_\ell:B \rightarrow {\Cal M}_n$ 
extending $f$ for \newline

$\qquad \qquad \qquad \qquad \ell < m$ such that 
$\langle g_\ell:\ell < m \rangle$ is disjoint over $A\biggr)$.
\endroster
\medskip

\noindent
(For $(b) \Rightarrow (a),{\frak K}$ weakly nice is not needed.)
\endproclaim
\bigskip

\demo{Proof}  1) The direction $(a) \Rightarrow (b)$ holds as ${\frak K}$
is weakly nice, more elaborately, we prove (b) assuming (a) by induction on 
$m$ where $\langle A_\ell:\ell \le m \rangle$ is a decomposition of $(A,B)$ 
(which exists by \scite{2.6}(5)); for $m=0$ this is trivial and for $m+1$ by 
straight combinatorics.  
Next we prove $\neg(a) \Rightarrow \neg(b)$ even
without using ``${\frak K}$ weakly nice".  So assume $(b) \and \neg(a)$ 
and we shall get a contradiction.  As $\neg(a)$, by \scite{2.6}(4) for some 
$A_1$ we have $A <_i A_1 \le B$, hence by Definition \scite{2.3}(2)(b) 
for some $m^* \in \Bbb N$ we have:
\medskip
\roster
\item "{$(*)_1$}"  $1 = \text{ Lim}_n \text{ Prob}_{\mu_n}\biggl($ 
for every embedding $f$ of $A$ into ${\Cal M}_n$ there are  \newline

$\qquad \qquad \qquad \qquad$ at most $m^*$ extensions of $f$ to an 
embedding of \newline

$\qquad \qquad \qquad \qquad A_1$ into ${\Cal M}_n \biggr)$.
\endroster
\medskip

\noindent
As Clause (b) holds, the limit there is $1$ also for the $m^*$ we have just
chosen.  The contradiction is immediate. \hfill$\square_{\scite{2.6E}}$
\enddemo
\bigskip

\definition{\stag{2.6G} Definition} 1) We say $(A^+,A,B,D)$ is a 
semi-$(k,r)$-good quadruple if:
\medskip

$(*)^{k,r}_{A^+,A,B,D}$  $\quad A \le A^+ \in {\Cal K}_\infty$ and 
$A \le D,B \le D \in {\Cal K}_\infty$ \newline 

$\qquad \qquad \qquad$ and for every random enough ${\Cal M}_n$ we have:
\medskip
\roster
\item "{$(**)$}"  for every embedding $f:A^+ \rightarrow {\Cal M}_n$ 
satisfying \newline
$c \ell^r(f(A),{\Cal M}_n) \subseteq f(A^+)$ \ub{there is} an extension 
$g$ of $f \restriction A$, embedding \newline
$D$ into ${\Cal M}_n$ such that \newline
$c \ell^k(g(B),{\Cal M}_n) = g(c \ell^k(B,D))$.
\endroster
\medskip

\noindent
If $r=k$ we may write $k$ instead of $(k,r)$. \nl
2) We say that ${\frak K}$ is semi-nice \underbar{if} it is weakly nice and
for every $A \in {\Cal K}_\infty$ and $k$ for some $\ell,m,r$ we have:
\medskip
\roster
\item "{$(*)$}"  for every random enough ${\Cal M}_n$, and embedding
$f:A \rightarrow {\Cal M}_n$ and $b \in {\Cal M}_n$ we can find $A_0,A^+,B,D$ 
such that: \nl
[note that we can have finitely many possibilities for $(\ell,m,r)$; does not
matter]
{\roster
\itemitem{ $(\alpha)$ }  $f(A) \le A_0 \le A^+ \le c \ell^m(f(A),
{\Cal M}_n)$ 
\sn
\itemitem{ $(\beta)$ }  $f(A) \cup \{b\} \subseteq B \subseteq D \subseteq 
{\Cal M}_n$
\sn
\itemitem{ $(\gamma)$ }  $|D| \le \ell$
\sn
\itemitem{ $(\delta)$ }  $(A^+,A_0,B,D)$ is semi-$(k,r)$-good
\sn
\itemitem{ $(\varepsilon)$ }  $c \ell^r(A_0,{\Cal M}_n) \subseteq A^+$
\sn
\itemitem{ $(\zeta)$ }    $c \ell^k(B,{\Cal M}_n) \subseteq D$. 
\endroster}
\endroster
\enddefinition
\bigskip

\proclaim{\stag{2.6H} Claim}  1) Assume
\mr
\item "{$(a)$}"  $(A^+,A,B,D)$ is semi-$(k,r)$-good
\sn
\item "{$(b)$}"  $A \le A^+_1 \le_s A^+$
\sn
\item "{$(c)$}"  $B_1 \le B,A \le D_1 \le D,B_1 \le D_1$
\sn
\item "{$(d)$}"  $c \ell^r(A,A^+) \subseteq A^+_1$ (follows from (b))
\sn
\item "{$(e)$}"  $c \ell^k(B_1,D) \subseteq D_1$.
\ermn
\ub{Then} $(A^+_1,A,B_1,D_1)$ is semi-$(k,r)$-good. \nl
2) If $(a)$ of part (1) and $k' \le k,r' \ge r$ and $A^+ \subseteq A^* \in
{\Cal K}_\infty$ satisfies $c \ell^r(A,A^*) \subseteq A^+$, \ub{then}
$(A^*,A,B,D)$ is semi-$(k',r')$-good.  (We can combine parts (1) and (2)). \nl
3) If (a) of part (1) and $r,r_1,m$ satisfies the statement $(*)$
below and $A_1 \le A \le A^+ \le A^+_1 \in {\Cal K}_\infty,
A \subseteq c \ell^m(A_1,A^+_1)$ and $c \ell^r(A,A^+_1) \subseteq A^+$,
\ub{then} $(A^+_1,A_1,B,D)$ is semi-$(k,r_1)$-good where
\medskip

$(*) = (*)^{r_1}_{m,r} \qquad$  if $A' \le B' \in {\Cal K}_\infty$ then

$\qquad \qquad \qquad \quad c \ell^r(c \ell^m(A',B'),B') 
\subseteq c \ell^{r_1}(A',B')$.
\endproclaim
\bigskip

\demo{Proof}  1) Let ${\Cal M}_n$ be random enough and $f:A^+_1 \rightarrow
{\Cal M}_n$ be such that \nl
$c \ell^r(f(A),{\Cal M}_n) \subseteq f(A^+_1)$.
As $A^+_1 \le_s A^+$ and ${\frak K}$ is semi-nice there is an embedding 
$f':A^+ \rightarrow {\Cal M}_n$ 
extending $f$.  By monotonicity, $c \ell^r(f(A),{\Cal M}_n) \subseteq
f'(A^+)$.  As $(A^+,A,B,D)$ is semi-$(k,r)$-good, there is $g:D \rightarrow
{\Cal M}_n$ such that: $g \supseteq f' \restriction A$ and $c \ell^k(g(B),
{\Cal M}_n) = g(c \ell^k(B,D)) \subseteq g(D)$. \nl
But by \scite{2.6}(14) clause (c) always
\mr
\item "{$(*)$}"  $A' \subseteq B' \subseteq C' \in {\Cal K}_\infty,
c \ell^{k'}(A',C') \subseteq B' \Rightarrow c \ell^{k'}(A',C') =
c \ell^{k'}(A',B')$
\ermn
so in our case (see assumption $(e)$)

$$
c \ell^k(g(B_1),{\Cal M}_n) = c \ell^k(g(B_1),g(D)).
$$
\mn
As $g$ embeds $D$ into ${\Cal M}_n$, 
$c \ell^k(g(B_1),g(D)) = g(c \ell^k(B_1,D))$,
and by $(*)$ above and assumption (e) we have 
$c \ell^k(B_1,D) = c \ell^k(B_1,D_1)$.  So together with earlier equality

$$
c \ell^k(g(B_1),{\Cal M}_n) = g(c \ell^k(B_1,D_1))
$$
\mn
as required, that is $g \restriction D_1$ is as required. \nl
2) Easier (by \scite{2.6}(14),clause(e)). \nl
3) Let ${\Cal M}_n$ be random enough and $f:A^+_1 \rightarrow {\Cal M}_n$
be such that $c \ell^{r_1}(f(A_1),{\Cal M}_n) \subseteq f(A^+_1)$.  We know
$A \subseteq c \ell^m(A_1,A^+_1)$.  Now by the assumption on $r_1,m,r$ for
every $A',A^+_1 \le A' \in {\Cal K}^\infty$ we have $c \ell^{r_1}(A_1,A')
\supseteq c \ell^r(A,A')$ hence $c \ell^r(f(A),{\Cal M}_n) \subseteq
c \ell^{r_1}(f(A_1),{\Cal M}_n) \subseteq f(A^+_1)$.  So we can apply the
property ``$(A^+,A,B,D)$ is semi-$(k,r)$-nice".
\hfill$\square_{\scite{2.6H}}$
\enddemo
\bigskip

\proclaim{\stag{2.6I} Claim}  In the Definition of semi-nice, \scite{2.6G}(2), 
we can equivalently omit $\ell,m$ (just $r$ suffice) and replace $(*)$ by
\mr
\item "{$(*)'$}"  for every random enough ${\Cal M}_n$ and $f:A \rightarrow
{\Cal M}_n$ and $b \in {\Cal M}_n$ letting $B = A \cup \{b\}$ we have: \nl
$(c \ell^r(f(A),{\Cal M}_n),f(A),B,f(A) \cup c \ell^k(B,{\Cal M}_n))$ is
semi-$(k,r)$-good.
\endroster
\endproclaim
\bigskip

\demo{Proof}  \ub{Original definition implies new definition}

Let $\ell,m,r$ be as guaranteed by the original definition.  Without loss of
generality $m \ge r$.  For the new definition we choose $r_1 > r,m$ 
such that $A' \le B' \in {\Cal K}_\infty \Rightarrow c \ell^r
(c \ell^m(A',B'),B') \subseteq c \ell^{r_1}(A',B')$, they exist by
\scite{2.6}(14)(e),(g).  Let us check $(*)'$ of the
new definition, so ${\Cal M}_n$ is random enough and $f:A \rightarrow
{\Cal M}_n$.  So by the old definition there are $A_0,A^+,B,D$ satisfying
$(\alpha)-(\zeta)$ of $(*)$.
In particular $(A^+,A_0,B,D)$ is semi-$(k,r)$-good.  
As $A^+ \le c \ell^m(f(A),{\Cal M}_n)$ and $f(A) \le A_0$ by \scite{2.6H}(3)
also the quadruple $(c \ell^{r_1}(f(A),{\Cal M}_n),f(A),B,D)$ is
semi-$(k,r_1)$-good.  (Remember $r_1 \ge m,r$).  Let $B_1 = f(A) \cup
\{b\},D_1 = f(A) \cup c \ell^k(B_1,{\Cal M}_n)$.  Now apply \scite{2.6H}(1)
and get that $(c \ell^{r_1}(f(A),{\Cal M}_n),f(A),f(A) \cup \{b\},f(A) \cup
c \ell^k(f(A) \cup \{b\},{\Cal M}_n))$ is semi-$(k,r_1)$-good as required.
\enddemo
\bn
\ub{New definition implies old definition}

Immediate, letting $m=r$: in $(*)$ of \scite{2.6G}(2) let $A_0 = f(A),
A^+ =$ \nl
$c \ell^r(f(A),{\Cal M}_n),
B=f(A) \cup \{b\}$ and $D = c \ell^k(f(A),{\Cal M}_n)$.  What about $\ell$?
It exists by \scite{2.6}(13).
\hfill$\square_{\scite{2.6I}}$
\bigskip

\demo{\stag{2.6J} Conclusion}  The definition of semi-nice is equivalent
to: \nl
for every $k$ and $\ell$ for some $r$ we have
\mr
\item "{$(*)''$}"  If $A \in {\Cal K}_\infty,|A| \le \ell$ and ${\Cal M}_n$
is random enough and $f:A \rightarrow {\Cal M}_n$ and $b \in {\Cal M}_n$
\ub{then} $(c \ell^r(f(A),{\Cal M}_n),f(A),f(A) \cup \{b\},c \ell^k(f(A) 
\cup \{b\},{\Cal M}_n))$ is semi-$(k,r)$-good.
\endroster
\enddemo
\bigskip

\demo{Proof}  \ub{old $\Rightarrow$ new}.

Let $\{A_i:i < i^*\} \subseteq {\Cal K}_\infty$ list the $A \in 
{\Cal K}_\infty$ with $\le \ell$ elements up to isomorphism.  For each
$A_i$ there is $r_i \in \Bbb N$ as guaranteed in \scite{2.6I}.  Let
$r = \text{ Max}\{r_i:i < i^*\}$. \nl
So let $A \in {\Cal K}_\infty,|A| \le \ell$ be given so for some $i,
A \cong A_i$; if ${\Cal M}_n$ is random enough and $f:A \rightarrow 
{\Cal M}_n$ and $b \in {\Cal M}_n$ and $B = f(A) \cup \{b\}$, then
$(c \ell^{r_i}(f(A),{\Cal M}_n),f(A),B$, \nl
$c \ell^k(B,{\Cal M}_n))$ is semi-$(k,r_i)$-good. \nl
(Why?  By the choice of $r_i$.)  Now by \scite{2.6H}(2) as $r_i \le r$ we know
that $(c \ell^r(f(A),{\Cal M}_n)$, \nl
$f(A),B,f(A) \cup c \ell^k(B,{\Cal M}_n))$
is semi-$(k,r)$-good, as required because $f(A) \subseteq c \ell^k(B,
{\Cal M}_n)$.
\enddemo
\bn
\ub{New $\Rightarrow$ old}

Easier (and not used). \hfill$\square_{\scite{2.6J}}$
\bigskip

\definition{\stag{2.9} Definition}  Let ${\frak K}$ be a 0-1 context.
\newline
1) ${\frak K}$ is complete \underbar{if} for every $A \in {\Cal K}$, 
the sequence

$$
\langle \text{Prob}_{\mu_n} (A \text{ is embeddable into } {\Cal M}_n):
n \in \Bbb N \rangle
$$

\noindent
converges to zero or converges to one. \newline
2)  ${\frak K}$ is weakly complete \underbar{if} 
the sequence above converges. \newline
3)  ${\frak K}$ is very weakly complete \underbar{if} for every 
$A \in {\Cal K}$, the sequence \newline
$\langle[\text{Prob}_{\mu_n}(A \text{ embeddable into }
{\Cal M}_{n+1}) - \text{ Prob}_{\mu_n}(A \text{ embeddable into }
{\Cal M}_n)]:n < \omega \rangle$ converges to zero.
\newline
So if $h(n) = n+1$, we get very weakly complete (similarly in
\scite{2.10}(4)). \nl
4)  ${\frak K}$ is $h$-very weakly complete if for every $A \in {\Cal K}$,
the sequence \newline
$\langle \text{Lim}_n \text{ Max }_{n_1,n_2 \in [n,h(n))}
[\text{Prob}_{\mu_{n_2}}(A \text{ embeddable into }
{\Cal M}_{n_2})$ \newline
$ - \text{ Prob}_{\mu_{n_1}}(A \text{ embeddable into }
{\Cal M}_{n_1})]:n < \omega \rangle$ converges to zero.
\enddefinition
\bigskip

\definition{\stag{2.10} Definition}  Let ${\frak K}$ be a 0-1 context.
\newline
1) ${\frak K}$ satisfies the 0-1 law for the logic ${\Cal L}$ 
\underbar{if} for every sentence $\varphi \in {\Cal L}(\tau)$
(i.e. the logic ${\Cal L}$ with vocabulary $\tau$) the sequence

$$
\langle \text{Prob}_{\mu_n}({\Cal M}_n \models \varphi):n \in \Bbb N \rangle
$$

\noindent
converges to zero or converges to one. \newline
2) ${\frak K}$ satisfies the weak 0-1 law or convergence law for the 
logic ${\Cal L}$ \underbar{if} for every sentence 
$\varphi \in {\Cal L}(\tau)$, the sequence

$$
\langle \text{Prob}_{\mu_n}({\Cal M}_n \models \varphi):n \in \Bbb N \rangle
$$

\noindent
converges. \newline
3)  ${\frak K}$ satisfies the very weak 0-1 law for ${\Cal L}$ \underbar{if}
for every sentence $\varphi \in {\Cal L}(\tau)$ the sequence

$$
\langle \text{Prob}_{\mu_{n+1}}({\Cal M}_{n+1} \models \varphi) -
\text{ Prob}_{\mu_n}({\Cal M}_n \models \varphi):n \in \Bbb N \rangle
$$

\noindent
converges to zero. \newline
4) ${\frak K}$ satisfies the $h$-very weak 0-1 law for ${\Cal L}$
\underbar{if} for every sentence $\varphi \in {\Cal L}(\tau)$, the sequence

$$
\langle \text{max}_{\{n_1,n_2\} \subseteq [n,n+h(n)]} 
|\text{Prob}_{\mu_{n_1}}({\Cal M}_{n_1} \models \varphi) - \text{ Prob}
_{\mu_{n_2}}({\Cal M}_{n_2} \models \varphi)|:n \in \Bbb N \rangle
$$
\medskip

\noindent
converge to zero. \newline
5)  If the logic ${\Cal L}$ is first order logic, we may omit it.
\enddefinition
\bigskip

\fakesubhead{\stag{2.11} Fact}  \endsubhead  1) If ${\frak K}$ is complete, 
\ub{then} it is weakly complete
which implies it is very weakly complete. \newline
2) If $h_1,h_2$ are functions from $\Bbb N$ to $\Bbb N$ and $(\forall n)
(h_1(n) \le h_2(n))$ and ${\frak K}$ is $h_2$-very weakly complete, \ub{then}
${\frak K}$ is $h_1$-very weakly complete. \newline
3)  Similarly for 0-1 laws.
\bigskip

\proclaim{\stag{2.14} Lemma}  1) Assume ${\frak K}$ is semi-nice.  
Modulo the theory $T_\infty$, every formula of the form 
$\psi(x_0,\dotsc,x_{m-1})$ is equivalent to a Boolean combination
of formulas of the form 
$(\exists x_m,\dotsc,x_{k-1})\varphi(x_0,\dotsc,
x_{m-1},x_m,\dotsc,x_{k-1})$, where for some \newline
$A \le_i B \in {\Cal K}_\infty$ we have 
$A = \{a_\ell:\ell < m\}$, $B = \{a_\ell:\ell < k\}$ (so $m \le k$) and

$$
\align
\varphi = \bigwedge \biggl\{ R(\ldots,x_\ell,\dotsc,)_{\ell < k}:&B \models 
R(\ldots,a_\ell,\ldots)_{\ell < k}, \\
  &R \text{ an atomic or negation of atomic formula} \biggr\}.\endalign
$$
\medskip

\noindent
1A)  Another way of saying it, is: there is $k$ computable from $\psi$ such 
that: \newline
for every random enough ${\Cal M}_n$ and $a_0,\dotsc,a_{m-1} \in {\Cal M}_n$,
the truth value of \newline
${\Cal M}_n \models \psi(a_0,\dotsc,a_{m-1})$ is computable from \newline
$({\Cal M}_n \restriction c \ell^k(\{a_0,\dotsc,a_{m-1}\},{\Cal M}_n),
a_0,\dotsc,a_{m-1})/\cong$. \newline
2) If ${\frak K}$ is semi-nice and weakly complete, \underbar{then} the 
weak 0-1 holds (i.e. convergence see Definition \scite{2.10}(2)).
\newline
3)  If ${\frak K}$ is semi-nice and complete (see Definition \scite{2.9}) 
\underbar{then} $T_\infty$ is a
complete theory; and ${\frak K}$ satisfies the $0-1$ law for first order
sentences (see \scite{2.10}(1)). \newline
4)  If $T_\infty$ is a complete theory, \underbar{then} 
${\frak K}$ is complete. \nl
5) The parallel of 2), 3) holds for $h$-very weak.
\endproclaim
\bigskip

\demo{Proof}  1) By (1A). \newline
1A) We prove it by induction on the quantifier depth of $\psi$.  For
$\psi$ atomic, or a conjunction or a disjunction or a negation this should be
clear.  So assume $\psi(x_0,\dotsc,x_{s-1}) = (\exists x_s)\varphi(x_0,\dotsc,
x_s)$, by the induction hypothesis there is a function $F_\varphi$ and
number $k_\varphi$ such that:
\medskip
\roster
\item "{$(*)_{\varphi,F_\varphi}$}"  for every random enough ${\Cal M}_n$ for
every $a_0,\dotsc,a_s \in {\Cal M}_n$ we have: the truth value of
${\Cal M}_n \models \varphi(a_0,\dotsc,a_s)$ is \newline
$F_\varphi\biggl( ({\Cal M}_n \restriction c \ell^{k_\varphi}(\{a_0,\dotsc,
a_s\},{\Cal M}_n),a_0,\dotsc,a_s)/\cong \biggr)$.
\endroster
\medskip

By Definition \scite{2.6G}(2) (of semi-nice) for any 
$A \in {\Cal K}_\infty$ and $k$ there are $\ell = \ell(A,k)$ and 
$m = m(A,k)$ and $r = r(A,k)$ as there.  Let
$m^* = \text{ max}\{m(A,k_\varphi):A \in {\Cal K}_\infty$ and $|A| \le
s + 1\}$ and $\ell^* = \text{ max}\{\ell(A,k_\varphi):A \in 
{\Cal K}_\infty$ and $|A| \le s+1\}$ and see below $(*)_5$(ii) and 
$(*)_6$(ii) and let $r^* = \text{ max}\{r(A,k_\varphi):A \in 
{\Cal K}_\infty$ and $|A| \le s+1\}$.  Let $k = k^*$ be such that 
$A \in {\Cal K}_\infty \and
|A| \le s \Rightarrow c \ell^k(A,{\Cal M}_n) \supseteq 
c \ell^{k_\varphi}(c \ell^{m^*}(A,{\Cal M}_n),{\Cal M}_n)$ 
(such $k$ exists by \scite{2.6}(14(g)).
\medskip

Now for ${\Cal M}_n$ random enough, for any $a_0,\dotsc,a_{s-1} \in
{\Cal M}_n$, we shall prove a sequence of conditions one implying the next
(usually also the inverse), then close the circle thus proving they are all 
equivalent:
\medskip
\roster
\item "{$(*)_1$}"  ${\Cal M}_n \models \psi(a_0,\dotsc,a_{s-1})$
\smallskip
\noindent
\item "{$(*)_2$}"  ${\Cal M}_n \models (\exists x_s) \varphi
(a_0,\dotsc,a_{s-1},x_s)$.
\endroster
\medskip

\noindent
[Clearly $(*)_1 \Leftrightarrow (*)_2$]
\medskip
\roster
\item "{$(*)_3$}"  for some $b \in {\Cal M}_n$ we have
${\Cal M}_n \models \varphi(a_0,\dotsc,a_{s-1},b)$.
\endroster
\medskip

\noindent
[Clearly $(*)_2 \Leftrightarrow (*)_3$]
\medskip
\roster
\item "{$(*)_4(i)$}"  for some $b \in c \ell^{m^*}(\{a_0,\dotsc,
a_{s-1}\},{\Cal M}_n)$ we have \newline
${\Cal M}_n \models \varphi(a_0,\dotsc,a_{s-1},b)$
\endroster
\medskip

\noindent
or
\roster
\item "{${}(ii)$}"  for some $b \in {\Cal M}_n \backslash
c \ell^{m^*}(\{a_0,\dotsc,a_{s-1}\},{\Cal M}_n)$ \newline
we have ${\Cal M}_n \models \varphi(a_0,\dotsc,a_{s-1},b)$.
\endroster
\medskip

\noindent
[Clearly $(*)_3 \Leftrightarrow (*)_4$]
\medskip
\roster
\item "{$(*)_5(i)$}"   letting $N = {\Cal M}_n \restriction
c \ell^k(\{a_0,\dotsc,a_{s-1}\},{\Cal M}_n)$ the following holds:
\newline
for some $b \in c \ell^{m^*}(\{a_0,\dotsc,a_{s-1}\},N) = c \ell^{m^*}
(\{a_0,\dotsc,a_{s-1}\},{\Cal M}_n)$ we have \newline
truth 
$= F_\varphi \biggl( (N \restriction c \ell^{k_\varphi}(\{a_0,\dotsc,
a_{s-1},b\},N),a_0,\dotsc,a_{s-1},b)/\cong \biggr)$.
\endroster
\medskip

\noindent
[Clearly $(*)_4(i) \Leftrightarrow (*)_5(i)$ by the choice of $k$ as: 
$A \subseteq N \subseteq
M \and c \ell^t(A,M) \subseteq N \Rightarrow c \ell^t(A,N) = c \ell^t
(A,M)$, by \scite{2.6}(13)(c) and the induction hypothesis].
\medskip
\roster
\item "{$(*)_5(ii)$}"  letting $N = c \ell^k(\{a_0,\dotsc,a_{s-1}\},
{\Cal M}_n)$ and $A = \{a_0,\dotsc,a_{s-1}\}$ we have: \newline
for some $A_0,A^+$ we have $A \le A_0 \le A^+ \subseteq c 
\ell^{m(A,k_\varphi)}(\{a_0,\dotsc,a_{s-1}\},N)$, and 
$c \ell^{r(A,k_\varphi)}(A_0,{\Cal M}_n) \subseteq A^+$ and there are 
$B^+,b$ such that
$c \ell^{k_\varphi}(A_0 \cup \{b\},{\Cal M}_n) \subseteq B^+ \subseteq 
{\Cal M}_n,|B^+| \le \ell^*$ and $(A^+,A_0,A \cup \{b\},B^+)$ is 
semi-$(k_\varphi,r(A,k_\varphi))$-good and $M \models 
\varphi(a_0,\dotsc,a_{s-1},b)$.
\endroster
\medskip

\noindent
[Clearly $(*)_4(ii) \Rightarrow (*)_5(ii)$ again by \scite{2.6}(14) and
${\frak K}$ being semi-nice (and see \scite{2.6H}(1)) which implies that 
we can use $A \cup \{b\}$ as $B$.].
\newline
Next let
\medskip
\roster
\item "{$(*)_6(ii)$}"    letting $N = c \ell^k(\{a_0,\dotsc,a_{s-1}\},
{\Cal M}_n),A = \{a_0,\dotsc,a_{s-1}\}$ we have: \newline
for some $A^+,A_0$ we have: $A \le A_0 \le A^+ \subseteq c 
\ell^{m(A,k_\varphi)}(\{a_0,\dotsc,a_{s-1}\},{\Cal M}_n)$ and 
$c \ell^{r(A,k_\varphi)}(A_0,{\Cal M}_n) \subseteq A^+$ and there are 
$B^+,b$ such that
$N \restriction (A \cup \{b\}) \le B^+ \in {\Cal K}_\infty,|B^+| \le \ell^*$ 
and $(A^+,A_0,A \cup \{b\},B^+)$ is semi-$(k_\varphi,r(A,k_\varphi))$-good 
and truth $= F_\varphi[B^+ \restriction c \ell^{k_\varphi}(\{a_0,\dotsc,
a_{s-1},b\},B^+),a_0,\dotsc,a_{s-1},b)/\cong]$.
\endroster
\medskip

\noindent
Clearly by the induction hypothesis $(*)_5(ii) \Rightarrow (*)_6(ii)$.

Lastly $(*)_6(ii) \Rightarrow (*)_3$ by Definition \scite{2.6G}(1) + the
induction hypothesis thus we have
equivalence. So $(*)_1 \Leftrightarrow [(*)_5(i) \vee (*)_6(ii)]$, 
but the two later ones depend just on
${\Cal M}_n \restriction c \ell^k(\{a_0,\dotsc,a_{s-1}\},
{\Cal M}_n),a_0,\dotsc,a_{s-1}))/\cong$, thus we have finished. \newline
2) By \scite{2.14}(1) it is enough to prove that the sequence

$$
\langle\text{Prob}_{\mu_n}(A \text{ is embeddable into } {\Cal M}_n):
n < \omega \rangle
$$

\noindent
converge.  This holds by weak completeness. \newline
3),4),5)  Left to the reader. \hfill$\square_{\scite{2.14}}$
\enddemo
\bigskip

\remark{\stag{2.15} Remark}  Note: \newline
If ${\frak K}$ is complete, \underbar{then} $T_\infty$ has a unique (up to
isomorphisms) countable model $M$ such that for some $\langle A_n:n \in
\Bbb N \rangle$ we have: $M = \dsize \bigcup_{n \in \Bbb N} A_n,A_n <_s 
A_{n+1} \in {\Cal K}_\infty$ and every $A \in {\Cal K}_\infty$ 
can be embedded into some $A_n$, and if $n \in \Bbb N,A \le A_n,A \le_s B$, 
then for some $m$ there is an
embedding $f$ of $B$ into $A_m$ such that $f \restriction A = \text{ id}_A$
and $f(B) \le_s A_m$ (see Baldwin, Shelah \cite{BlSh:528}, not used).
\endremark
\bigskip

\proclaim{\stag{2.16} Claim}  1) A sufficient 
condition for ``${\frak K}$ is weakly nice" is
\medskip
\roster
\item "{$(*)$}"  for every $A < B$, if $\neg(A \le_a B)$ \underbar{then} for
some $k < \omega$ we have \newline
$1 = \text{ Lim}_n \text{ Prob}_{\mu_n}\biggl($for every 
embedding $f$ of $A$ into
${\Cal M}_n$ there are no \newline

$\qquad \qquad \qquad \qquad$ embedding $g_\ell:B \rightarrow {\Cal M}_n$
extending $f$ for \newline

$\qquad \qquad \qquad \qquad \ell < k$, disjoint over $f\biggr)$.
\endroster
\endproclaim
\bigskip

\demo{Proof}  Easy.
\enddemo
\newpage

\head {\S2 More accurate measure and drawing again} \endhead  \resetall
\bigskip

We define when a pair of functions $\bar h = (h^d,h^u)$ giving, up to a
factor, the values nu$(f,A,B,M)$ for $M$ random enough, $A <_{\text{pr}} B$
and $f:A \rightarrow M$.  We also define when $\bar h$
obeys $h$ (i.e. $h$ bounds the error factor).  
\relax From $h_{A,B}$ for $A <_{\text{pr}} B$ we define
$h_{A,B}$ also for the case $A <_s B$ and then define a good case when 
the functions are polynomial
in $\|{\Cal M}_n\|$ (see Definition \scite{3.1}, \scite{3.3}).

We then see how large is the factor error for the derived cases and deduce
some natural properties (in \scite{3.4}).

Then we deal with the polynomial case.

Lastly, (\scite{3.7}-\scite{3.12}) we introduce our framework for adding
random relations to random ${\Cal M}_n$.  Reading, you may 
assume ``every $A \in K_\infty$ is embeddable into every random 
enough ${\Cal M}_n"$.
\bigskip

\definition{\stag{3.1} Definition}  1) We say the $0-1$ context ${\frak K}$ 
\underbar{obeys} $\bar h = (h^d,h^u)$ with error $h^e$ where $d$ is for 
down, $u$ is for up and $e$ is for error \ub{if}:
\medskip
\roster
\item "{(a)}"  for $A <_{\text{pr}} B$ we have $h^d_{A,B}$ and 
$h^u_{A,B}$ and $h^e$ are functions from $\dbcu_{n < \omega}{\Cal K}_n$ to 
$\Bbb R^{\ge 0}$
\sn
\item "{(b)}"  for some $\varepsilon \in \Bbb R^{>0}$ for every random 
enough ${\Cal M}_n$ we have \nl
$(h^e[{\Cal M}_n])^\varepsilon \le h^d_{A,B}
[{\Cal M}_n] \le h^u_{A,B}[{\Cal M}_n]$ and $h^e[{\Cal M}_n] \ge 1$
\sn
\item "{(c)}"   for every $\varepsilon \in \Bbb R^{>0}$ we have \nl
$1 = \text{ Lim}_n \text{ Prob}_{\mu_n} \biggl($ 
for every embedding $f$ of $A$ into ${\Cal M}_n$, \newline
we have $h^d_{A,B}[{\Cal M}_n] \times (h^e[{\Cal M}_n])^{- \varepsilon}
\le nu(f,A,B,{\Cal M}_n) \le h^u_{A,B}[{\Cal M}_n] \times 
(h^e[{\Cal M}_n])^\varepsilon \biggr)$ \newline
(see \scite{1.3}(2)).
\ermn
1A) If $h^e$ is identically 1 we may omit it.  If $h^u = h^d$, then we may
write $h^u$ instead of $\bar h$.  If $h^e[{\Cal M}_n] = \|{\Cal M}_n\|$ we may
say ``simply".
\medskip

\noindent
2)  We say $\bar h$ is uniform if $h^x_{A,B}[M]$ (for $x \in \{u,d\}$) 
depends on $\|M\|$ (and $x,A,B$) only but not on $M$ and then we write 
$h^x_{A,B}(\|M\|) = h^x_{A,B}[\|M\|]$.  Similarly for $h^e,h$ used above
and below.  We say $h$ goes to infinity if for every $m$ for every random
enough ${\Cal M}_n,h[{\Cal M}_n] > m$.  \newline
3) We say that $\bar h$ is bounded (or bounded$^+$) by $h$ (for 
$\bar h$ as above) \underbar{if}:
\medskip
\roster
\item "{$(a)$}"  $h$ is a function from $\dsize \bigcup_{n \in \Bbb N}
{\Cal K}_n$ to $\Bbb R^+$ 
(remember that the ${\Cal K}_n$'s are \nl
pairwise disjoint)
\sn
\item "{$(b)$}"   for every random enough ${\Cal M}_n$ we have
$h[{\Cal M}_n] \ge 1$
\sn
\item "{$(c)$}"  for every $\varepsilon > 0$ and $A <_{pr} B$ for every
random enough ${\Cal M}_n$ we have

$$
1 \le h^u_{A,B}[{\Cal M}_n]/h^d_{A,B}[{\Cal M}_n] \le 
(h[{\Cal M}_n])^\varepsilon
$$
\noindent
\item "{$(d)$}"  for every $A <_{pr} B$ and $m \in \Bbb N \backslash
\{0\}$ for some $\varepsilon > 0$ for every random enough 
${\Cal M}_n$ we have

$$
h^d_{A,B}[{\Cal M}_n] > (h[{\Cal M}_n])^\varepsilon \times m
$$
\noindent
\item "{$(e)$}"  for \footnote{this, of course, will not suffice for
0-1 law} every $A \in {\Cal K}_\infty$ for every
$\varepsilon > 0$ for every random enough ${\Cal M}_n$ we have

$$
\text{Prob}_{\mu_n}(A \text{ is embeddable into } {\Cal M}_n) \ge
1/(h[{\Cal M}_n])^\varepsilon.
$$
\endroster
\mn
3A) In the context of (3), let ``${\Cal M}_n$ random enough" mean that for
every $\varepsilon$, the probability of failure is 
$\le 1/(h[{\Cal M}_n])^\varepsilon$.  We say $h$ is standard if for each
$m$, for every random enough ${\Cal M}_n,h[{\Cal M}_n] > m$.  \nl
3B) We \footnote{this is an alternative to part (3), this does not matter
really so we shall use the one \nl
of part (3), the same applies to other cases}
say $\bar h$ is bounded$^-$ by $h$ (for $\bar h$) as above if: 
clauses (a), (b) as above, but in clauses (c),(e) we replace ``every 
$\varepsilon > 0$" by ``for some $m = m(A,B) \in \Bbb N$" in clause (d) 
we replace ``some $\varepsilon > 0$" by ``every $m \in \Bbb N$", i.e.
\medskip
\roster
\item "{$(a)$}"  $h$ is a function from $\dsize \bigcup_{n \in \Bbb N}
{\Cal K}_n$ to $\Bbb R^+$ 
(remember that the ${\Cal K}_n$'s are \nl
pairwise disjoint)
\sn
\item "{$(b)$}"   for every random enough ${\Cal M}_n$ we have
$h[{\Cal M}_n] \ge 1$
\sn
\item "{$(c)$}"  for every $A <_{pr} B$, for some $m = m(A,B) \in \Bbb N
\backslash \{0\}$ for every random enough ${\Cal M}_n$ we have

$$
1 \le h^u_{A,B}[{\Cal M}_n]/h^d_{A,B}[{\Cal M}_n] \le 
(h[{\Cal M}_n])^m
$$
\noindent
\item "{$(d)$}"  for every $A <_{pr} B$ and $m \in \Bbb N$ for some 
$m(A,B) \in \Bbb N \backslash \{0\}$ for every random enough ${\Cal M}_n$ 
we have

$$
h^d_{A,B}[{\Cal M}_n] > (h[{\Cal M}_n])^{m(A,B)},m
$$
\noindent
\item "{$(e)$}"  for \footnote{this, of course, will not suffice for
0-1 law, and though more natural, we shall not follow it here} every 
$A \in {\Cal K}_\infty$ for some $m(A) \in \Bbb N \backslash \{0\}$ 
for every random enough ${\Cal M}_n$ we have

$$
\text{Prob}_{\mu_n}(A \text{ is embeddable into } {\Cal M}_n) \ge
1/(h[{\Cal M}_n])^m
$$
\noindent
(part (3B) is an alternative to \scite{3.1}(3)).
\endroster
\medskip

\noindent
4) Assume ${\frak K}$ obeys $\bar h$. For $A <_s B$ and $M \in \dsize
\bigcup_n {\Cal K}_n$ we let

$$
\align
h^{+u}_{A,B}[M] =: \text{ Max}\biggl\{ \dsize \prod_{\ell < h}
h^u_{A_\ell,A_{\ell + 1}}[M]:&\bar A = \langle A_\ell:\ell \le k \rangle \\
  &\text{is a decomposition of } (A,B) \biggr\}
\endalign
$$

$$
\align
h^{-u}_{A,B}[M] =: \text{ Min}\biggl\{ \dsize \prod_{\ell < h}
h^u_{A_\ell,A_{\ell + 1}}[M]:&\bar A = \langle A_\ell:\ell \le k \rangle \\
  &\text{is a decomposition of } (A,B) \biggr\}
\endalign
$$

$$
\align
h^{-d}_{A,B}[M] =: \text{ Min}\biggl\{ \dsize \prod_{\ell < h}
h^d_{A_\ell,A_{\ell + 1}}[M]:&\bar A = \langle A_\ell:\ell \le k \rangle \\
  &\text{is a decomposition of } (A,B) \biggr\}
\endalign
$$

$$
\align
h^{+d}_{A,B}[M] =: \text{ Max}\biggl\{ \dsize \prod_{\ell < h}
h^d_{A_\ell,A_{\ell + 1}}[M]:&\bar A = \langle A_\ell:\ell \le k \rangle \\
  &\text{is a decomposition of } (A,B) \biggr\}
\endalign
$$

\noindent
Let $h^u_{A,B}[M] = h^{-u}_{A,B}[M]$ and $h^d_{A,B}(M) = h^{+d}_{A,B}(M)$.
\enddefinition
\bigskip

\demo{\stag{3.2} Discussion}  For the semi-nice case, we may consider it 
natural to have the functions $h$ below be $h^{x,k}_{A^+,A,B,D}$ giving 
information on \nl
nu$^k(f,A,B,D,M) = |ex^k(f,A,B,D,M)|$ where $ex^k(f,A,B,D,M) = 
\{g \restriction B:g$ embeds $D$ into $M$, it extends $f$ (which embeds $A$
into $M$) and $c \ell^k(g(B),M) \subseteq g(D)\}$ and we restrict ourselves
to the case that there is an embedding $f^+$ of $A^+$ into $M$ extending $f$
such that $c \ell^r(f(A),M) \subseteq f(A^+)$.  So we may write
$h^{x,k,r}_{A^+,A,B,D}$ and $ex^{k,r}(f,A,A^+,B,D,M)$.  Note that the
variables here of $ex,nu$ are different than in the usual case.
\enddemo
\bigskip

\definition{\stag{3.3} Definition}  1) We say ${\frak K}$ obeys the
polynomial $\bar h$ over (or modulo) $h$ \underbar{if} $\bar h = \langle
h^u,h^d \rangle$ and $h^u = h^d,h$ are functions from $\dbcu_n {\Cal K}_n$ to
$\Bbb R^{\ge 0}$ and $h^u,h$ are uniform (see Definition \scite{3.1}(2)) 
and for every $A <_{pr} B$ a real $\alpha(A,B) \in \Bbb R^{> 0}$ is well 
defined and we have:
\medskip
\roster
\item "{(a)}"  $h:\dbcu_{n < \omega} {\Cal K}_n \rightarrow \Bbb R^+$
\sn
\item "{(b)}"  for every $m$, for random enough ${\Cal M}_n$ we have
$h[{\Cal M}_n] > 1$ and $\|{\Cal M}_n\| \ge m$
\sn
\item "{(c)}"  for every $\varepsilon > 0$ for random enough
${\Cal M}_n$ we have $h[{\Cal M}_n] < \|{\Cal M}\|^\varepsilon$
\sn
\item "{(d)}"  if $A <_{\text{pr}} B$ and $m \in \Bbb N$, \ub{then} 
for every ${\Cal M}_n$ random enough 
\footnote{note this is not as in \scite{3.1}(3)(c)} \nl
$h^d_{A,B}[{\Cal M}_n] = h^u_{A,B}[{\Cal M}_n] =
\|{\Cal M}_n\|^{\alpha(A,B)}$; and if $f$ embeds $A$ into ${\Cal M}_n$
then \nl
$h[{\Cal M}_n]^{-m} h^d_{A,B}[{\Cal M}_n] \le \text{ nu}(f,A,B,{\Cal M}_n)
\le h[{\Cal M}_n]^m h^u_{A,B}[{\Cal M}_n]$
\sn
\item "{$(e)$}"   for some $\varepsilon > 0$ for every ${\Cal M}_n$ random
enough $\|{\Cal M}_n\| > m \times (h[{\Cal M}_n])^\varepsilon$
\sn
\item "{$(f)$}"  if $A \in {\Cal K}_\infty$ then for each $k$ for some $m$ \nl
Prob$_{\mu_n}(A$ is embeddable into ${\Cal M}_n$ assuming
$\|A\| = k) \ge 1/h[{\Cal M}_n]^m$.
\ermn
2) We say $\bar h$ is strictly polynomial if
\mr
\item "{$(a)$}"  if $A <_{pr} B$ then for some $c = c(A,B) \in
\Bbb R^{> 0}$ for some $\varepsilon > 0$, for every random enough
${\Cal M}_n$ and every $f:A \rightarrow {\Cal M}_n$ we have
\ermn
$$
\align
c(A,B)\|{\Cal M}_n\|^{\alpha(A,B)}(1-\|{\Cal M}_n\|^{-\varepsilon}) &\le
h^d_{A,B}({\Cal M}_n) \\
  &\le h^u_{A,B}({\Cal M}_n) \\
  &\le c(A,B)\|{\Cal M}_n\|^{\alpha(A,B)}(1+\|{\Cal M}_n\|^{-\varepsilon}).
\endalign
$$
\mn
3) We say $\bar h$ is a polynomial if $\bar h$ is polynomial over some $h$.
We say ${\frak K}$ is polynomial over $h$ (strictly polynomial), if this holds
for some $\bar h$.
\enddefinition
\bigskip

\demo{\stag{3.4} Fact}  1) Assume ${\frak K}$ obeys $\bar h$ with
error $h^e$. \newline
If $A_0 <_{pr} A_1 <_{pr} < \cdots <_{\text{pr}} A_k$ and $\varepsilon > 0$, 
\underbar{then} every random enough ${\Cal M}_n$ satisfies:
\medskip
\roster
\item "{$(*)$}"  for every embedding $f$ of $A_0$ into ${\Cal M}_n$,
$$
\align
\dsize \prod_{\ell < k} h^d_{A_\ell,A_{\ell + 1}}[{\Cal M}_n] &\times
(h^e[{\Cal M}_n])^{-\varepsilon} \le nu
(f,A_0,A_k,{\Cal M}_n) \\
  &\le \dsize \prod_{\ell < k} h^u_{A_\ell,A_{\ell + 1}}[{\Cal M}_n]
\times (h^e[{\Cal M}_n])^\varepsilon.
\endalign
$$
\endroster
\medskip

\noindent
2) Assume ${\frak K}$ obeys $\bar h$ with error $h^e$, and for 
clause $(\delta)$ assume also that $\bar h$ is bounded by $h$ and
$h \ge h^e$.
\mr
\item "{$(\alpha)$}"  if $A \le_s B$, \ub{then} for every random 
enough ${\Cal M}_n$ we have 
\footnote{if $A$ is embeddable into ${\Cal M}_n$ of course as
otherwise $h^x_{A,B}[{\Cal M}_n]$ is not actually well defined, we tend
to ``forget" to state this}

$$
\align
h^{-d}_{A,B}[{\Cal M}_n] \le h^{+d}_{A,B}[{\Cal M}_n] &= h^d_{A,B}[{\Cal M}_n]
\le h^u_{A,B}[{\Cal M}_n] \\
  &= h^{-u}_{A,B}[{\Cal M}_n] \le h^{+u}_{A,B}[{\Cal M}_n]
\endalign
$$
\sn
\item "{$(\beta)$}"  if $A <_s B$ and $\varepsilon \in \Bbb R^{>0}$, 
\ub{then} for every random enough ${\Cal M}_n$ and embedding \newline
$f:A \rightarrow {\Cal M}_n$ we have: \newline
$h^d_{A,B}[{\Cal M}_n] \times (h^e[{\Cal M}_n])^{- \varepsilon} 
\le nu(f,A,B,{\Cal M}_n) \le h^u_{A,B}[{\Cal M}_n] \times (h^e[{\Cal M}_n])
^\varepsilon$
\smallskip
\noindent
\item "{$(\gamma)$}"  if $A <_{pr} B$, \ub{then} \newline
$h^{+d}_{A,B}[M] = h^{-d}_{A,B}[M] = h^d_{A,B}[M]$ and \newline
$h^{+u}_{A,B}[M] = h^{-u}_{A,B}[M] = h^u_{A,B}[M]$
\smallskip
\noindent
\item "{$(\delta)$}"  if $A <_s B$ \underbar{then} for every $\varepsilon >
0$, for every random enough ${\Cal M}_n$ we have:
\newline
$(h[{\Cal M}_n])^{- \varepsilon} 
\le h^u_{A,B}[{\Cal M}_n] / h^d_{A,B}[{\Cal M}_n] \le 
(h[{\Cal M}_n])^\varepsilon$, moreover \newline
$(h[{\Cal M}_n])^{- \varepsilon} 
\le h^{+u}_{A,B}[{\Cal M}_n] / h^{-d}_{A,B}[{\Cal M}_n] \le (h[M])
^\varepsilon$
\smallskip
\noindent
\item "{$(\varepsilon)$}"  if $A_0 <_s A_1 <_s A_2$ \underbar{then}
for any random enough ${\Cal M}_n$: \newline
$h^d_{A_0,A_1}[{\Cal M}_n] \times h^d_{A_1,A_2}[{\Cal M}_n] \le 
h^d_{A_0,A_2}[{\Cal M}_n]$

$\qquad \qquad \quad h^u_{A_0,A_2}[{\Cal M}_n]
\le h^u_{A_0,A_1}[{\Cal M}_n] \times h^u_{A_1,A_2}[{\Cal M}_n]$ 
\endroster
\enddemo
\bigskip

\demo{Proof}  1)  Easy by induction on $k$. \newline
2) \underbar{Clause $(\alpha)$}: \newline
The first and last inequality holds as Min$(X) \le \text{ Max}(X)$ for
$X \subseteq \Bbb R$ finite non-empty (by \scite{2.6}(5)) as in this case.  
The equalities hold by Definition \scite{3.3}(1).
The middle inequality holds by clause $(\beta)$ below.
\medskip

\noindent
\underbar{Clause $(\beta)$}: \newline
By $(*)$ of \scite{3.4}(1).
\medskip

\noindent
\underbar{Clause $(\gamma)$}: \newline
As $A <_{\text{pr}} B$ implies $(A,B)$ has a unique decomposition.
\medskip

\noindent
\underbar{Clause $(\delta)$}: \newline
Let $\bar A = \langle A_\ell:\ell \le k \rangle$ be a decomposition of
$(A,B)$ and $\varepsilon \in \Bbb R^{>0}$, hence for every random 
enough ${\Cal M}_n$ for every embedding $f$ of $A$ into $M$ \newline

$$
\align
\dsize \prod_{\ell < k} h^d_{A_\ell,A_{\ell + 1}}[{\Cal M}_n] &\times
(h[{\Cal M}_n])^{- \varepsilon} \le
h^d_{A,B}[{\Cal M}_n] \times (h[{\Cal M}_n])^{- \varepsilon} \\
  &\le nu(f,A,B,{\Cal M}_n) \le h^u_{A,B}[{\Cal M}_n] \times 
(h[{\Cal M}_n])^\varepsilon \\
  &\le \dsize \prod_{\ell < k} h^u_{A_\ell,A_{\ell + 1}}[{\Cal M}_n]
\times (h[{\Cal M}_n])^\varepsilon.
\endalign
$$
\medskip

\noindent
This gives the first inequality part of the inequalities.  Let 
$\varepsilon > 0$ be given.

Now for each $\ell < k$ for every random enough ${\Cal M}_n$ we have

$$
1 \le h^u_{A_\ell,A_{\ell + 1}}[{\Cal M}_n]/
h^d_{A_\ell,A_{\ell + 1}}[{\Cal M}_n]
\le (h[{\Cal M}_n])^{\varepsilon/k}.
$$

\noindent
Hence for every random enough ${\Cal M}_n$ we have

$$
\align
h^u_{A,B}[{\Cal M}_n]/h^d_{A,B}[{\Cal M}_n] &\le \biggl(
\dsize \prod_{\ell < k} h^u_{A_\ell,A_{\ell + 1}}[{\Cal M}_n] \biggr)
/ \biggl( \dsize \prod_{\ell < h} 
h^d_{A_\ell,A_{\ell + 1}}[{\Cal M}_n] \biggr) \\
  &= \dsize \prod_{\ell < k} \biggl( h^u_{A_\ell,A_{\ell + 1}}[{\Cal M}_n]
/h^d_{A_\ell,A_{\ell + 1}}[{\Cal M}_n] \biggr) \\
  &\le \dsize \prod_{\ell < k}(h[{\Cal M}_n])^{\varepsilon/k} =
h[{\Cal M}_n]^\varepsilon.
\endalign
$$

\noindent
For the second phrase of clause $(\delta)$ (the moreover) 
note that for every random enough ${\Cal M}_n$ for every
$f:A \rightarrow {\Cal M}_n$ we have: for some decomposition $\bar A$
of $(A,B)$

$$
\align
1 \le h^{+u}_{A,B}[{\Cal M}_n]/nu(f,A,B,{\Cal M}_n) &= \dsize
\prod_{\ell < k} h^u_{A_\ell,A_{\ell +1}}[{\Cal M}_n] / 
nu(f(A,B,{\Cal M}_n)) \\
  &\le \dsize \prod_{\ell < k} h^u_{A_\ell,A_{\ell +1}}[{\Cal M}_n] /
\dsize \prod_{\ell < k} h^d_{A_\ell,A_{\ell + 1}}[{\Cal M}_n] \\
  &\le ((h[{\Cal M}_n])^{\varepsilon/k})^k = (h[{\Cal M}_n])^\varepsilon
\endalign
$$

\noindent
and for possibly other decomposition $\bar A$ of $(A,B)$

$$
\align
1 \le nu(f,A,B,{\Cal M}_n)/h^{-d}_{A,B}[{\Cal M}_n] &= nu(f(A,B,{\Cal M}_n))
/ \dsize \prod_{\ell < k} h^d_{A_\ell,A_{\ell +1}}[{\Cal M}_n] \\
  &\le \dsize \prod_{\ell < k} h^u_{A_\ell,A_{\ell + 1}}[{\Cal M}_n] /
\dsize \prod_{\ell < k} h^d_{A_\ell,A_{\ell + 1}}[{\Cal M}_n] \\
  &\le (h[{\Cal M}_n]^{\varepsilon/k})^k = h[{\Cal M}_n]^\varepsilon.
\endalign
$$  

\noindent
Together we get the desired inequality (well, for $2 \varepsilon$).
\medskip

\noindent
\underbar{Clause $(\varepsilon)$}: \newline
Easy (using \scite{3.1}(4)). \hfill$\square_{\scite{3.4}}$
\enddemo
\bigskip

\proclaim{\stag{3.4A} Claim}  Assume the 0-1 law context ${\frak K}$ 
obeys $\bar h = (h^d,h^u)$ with error $h^e$. \newline
1)  A sufficient condition for ``${\frak K}$ is weakly nice" is
\medskip
\roster
\item "{$(*)_1$}"  if $A <_{pr} B$ and $m^* \in \Bbb N$ and
$\varepsilon > \Bbb R^{>0}$ small enough, \underbar{then} for
every random enough ${\Cal M}_n$ we have 
$h^d_{A,B}[{\Cal M}_n] \times (h^e[{\Cal M}_n])^{- \varepsilon} > m^*$.
\endroster
\medskip

\noindent
2) If $\bar h$ is bounded by $h$ or $\bar h$ is polynomial over $h$, 
\underbar{then} $(*)_1$ from above holds hence ${\frak K}$ is weakly nice. 
\endproclaim
\bigskip

\demo{Proof}  1) Assume $A <_{\text{pr}} B$ and we show that $(*)_1$ of 
Definition \scite{2.6C} holds in this case. \newline
If $|B \backslash A| = 1$, as $A <_{pr} B$, for every $m$, for every random
enough ${\Cal M}_n$, for any $f:A \rightarrow {\Cal M}_n$, by $(*)_1$ of
\scite{3.4A}(1) there
are distinct $g_1,\dotsc,g_m$ which are embedding of $B$ into ${\Cal M}_n$
extending $f$, now they are necessarily pairwise disjoint over $A$ so the
demand in 1.3(1) holds.
\medskip

So assume $|B \backslash A| > 1$ and $m^*$ be given.  
For each $b \in B \backslash A$ by \scite{2.6}(8) we know that 
$(A \cup \{b\}) \le_i B$ but $|A \cup \{ b\}| \le |A| + 1 < |B|$ so 
$A \cup \{b\} <_i B$.  Hence by Definition \scite{2.3}(2)(b) for some
$n_b,m_b \in \Bbb N$ for every $n \ge n_b$ we have $1 - \varepsilon/(|B|
+1) < \text{ Prob}_{\mu_n}({\Cal E}^b_n)$ where ${\Cal E}^b_n$ is the event:
\medskip
\roster
\item "{{}}"  for every embedding $f$ of $(A \cup \{b\})$ into
${\Cal M}_n$ there are at most \newline
$m_b$ extensions of $f$ to an embedding of $B$ into ${\Cal M}_n$.
\endroster
\medskip

\noindent
Let $m^{**} = |B \backslash A| \times |B \backslash A| \times
(\underset {b \in B \backslash A} {}\to {\text{max}} m_b) \times m^*$.
\newline
Also by $(*)_1$ above and clause (c) of \scite{3.1}(1) for some 
$n^* \in \Bbb N$ for every $n \ge n^*$ the event 
${\Cal E}^*_n \cap {\Cal E}^{**}_n$ has probability 
$\ge 1 - \varepsilon/(|B| + 1)$ where ${\Cal E}^{**}_n =
[h^d_{A,B}[{\Cal M}_n] \le nu(f,A,B,{\Cal M}_n)$ for every embedding
$f:A \rightarrow {\Cal M}_n]$ and ${\Cal E}^* = [h^d_{A,B}[{\Cal M}_n] > 
(m^{**} + 1)]$. \newline
Let $n^\otimes = \text{ Max}(\{n_b:b \in B \backslash A\} \cup \{n^*\})$.
\newline
Now suppose $n \ge n^\otimes$, then with probability $\ge 1 - \varepsilon$ all
the events ${\Cal E}^b_n$ for $b \in B \backslash A$ and ${\Cal E}^*_n =
[h^d_{A,B}[{\Cal M}_n] > m^{**}]$ and ${\Cal E}^{**}_n$ occurs for ${\Cal M}_n$.  
It suffices to show that then $(*)_1$ of \scite{2.6C} occurs.  So let $f$ 
be an embedding of $A$ into ${\Cal M}_n$, so as both ${\Cal E}^*_n$ and
${\Cal E}^{**}_n$ occur necessarily there are distinct
extension $g_1,\dotsc,g_{m^{**}}$ of $f$ embedding $B$ into ${\Cal M}_n$.
For $i \in \{1,\dotsc,m^{**}\}$ let $u_i = \{j:j \in \{1,\dotsc,m^{**}\}$ and
Rang$(g_j) \cap \text{ Rang}(g_i) \ne \text{Rang}(f)\}$, and for $b \in B
\backslash A$ and $c \in {\Cal M}_n \backslash \text{Rang}(f)$ let
$v_{b,c} = \{i:i \in \{1,\dotsc,m^{**}\}$ and $g_i(b) = c\}$.  Now clearly
$|v_{b,c}| \le m_b$ as the event ${\Cal E}^*_b$ occurs and

$$
u_i = \dsize \bigcup_{b \in B \backslash A} \qquad \dsize \bigcup_{c \in
\text{ Rang}(g_i) \backslash \text{Rang}(f)} v_{b,c}
$$

\noindent
hence

$$
|u_i| \le |B \backslash A| \times |B \backslash A| \times
\underset {b \in B \backslash A} {}\to {\text{max}} m_b \le m^{**}/m^*.
$$

\noindent
So easily we can find $w \subseteq \{1,\dotsc,m^{**}\}$ such that
$|w| = m^*$ and $i \in w \and j \in w \and i \ne j \Rightarrow j
\notin u_i$.  So $\{g_i:i \in w\}$ is as required. \newline
2) Check (see in particular \scite{3.1}(3)(b),(d)).
\hfill$\square_{\scite{3.4A}}$
\enddemo
\bigskip

\proclaim{\stag{3.4B} Claim}  1) Assume ${\frak K}$ obeys $\bar h$ with
error $h^e$ and:
$A <_s B \le D$ and $A \le C \le_s D$
and $D = B \cup C$.  For every $\varepsilon \in \Bbb R^{>0}$, for every
random enough ${\Cal M}_n$, if $C$ is embeddable into ${\Cal M}_n$, \ub{then}

$$
h^u_{A,B}[{\Cal M}_n] \ge h^d_{C,D}[{\Cal M}_n] \times
(h^e[{\Cal M}_n])^{- \varepsilon}.
$$
\mn
2) If in addition, $\bar h$ is bounded by $h$ and $h \ge h^e$, \ub{then} 
for every
$\varepsilon > 0$, for every random enough ${\Cal M}_n$ and $x \in
\{u,d\}$

$$
h^x_{A,B}[{\Cal M}_n] \ge h^x_{C,D}[{\Cal M}_n] \times
(h[{\Cal M}_n])^\varepsilon
$$
\mn
3) If $A_0 \le A_1 \le \ldots \le A_k,A_0 \le_s A_k$ and $\varepsilon \in
\Bbb R^{>0}$, \ub{then} for every ${\Cal M}_n$ random enough into which 
$A_{k-1}$ is embeddable

$$
h^d_{A_0,A_k}[{\Cal M}_n] \le \Pi\{h^u_{A_i,A_{i+1}}[{\Cal M}_n]:A_i <_s
A_{i+1}\} \times (h^e[{\Cal M}_n])^\varepsilon.
$$
\endproclaim
\bigskip

\definition{\stag{3.4C} Definition}  1) We say that $({\frak K},\bar h,h^e)$ is
semi-nice if
\mr
\item "{$(a)$}"  ${\frak K}$ is a 0-1 context
\sn
\item "{$(b)$}"  ${\frak K}$ obeys $\bar h$ with error $h^e$
\sn
\item "{$(c)$}"   for every $A \in {\Cal K}_\infty$ and $k$ for some $r$ 
we have:
\sn
{\roster
\itemitem{ $(*)$ }  for every random enough ${\Cal M}_n$, and embedding
$f:A \rightarrow {\Cal M}_n$ and $b \in {\Cal M}_n,(c \ell^r(f(A),
{\Cal M}_n),f(A),f(A) \cup \{b\},c \ell^k(f(A) \cup \{b\},{\Cal M}_n))$
is semi$^*-(k,r)$-good for $({\frak K},\bar h,h^e)$, see below.
\endroster}
\item "{$(d)$}"  Condition $(*)_1$ of \scite{3.4A}(1) holds.
\ermn
2) $(A^*,A,B,D)$ is semi$^*-(k,r)$-good for $({\frak K},\bar h,h)$ if
\footnote{note that because of $A \le_s D$, this does not copy the definition
in \S1 even in ``nice" cases} for some $A_0$ we have:
\mr
\item "{$(\alpha)$}"  $A \le A_1 \le A^* \in {\Cal K}_\infty,A_0 \le_s D \in
{\Cal K}_\infty,B \le D \in K_\infty$ and
\sn
\item "{$(\beta)$}"   for every $\varepsilon > 0$, for every
random enough ${\Cal M}_n$, for every embedding $f^*$ of $A^*$ into
${\Cal M}_n$ satisfying $c \ell^r(f^*(A),{\Cal M}_n) \subseteq f^*(A^*)$, we
have, letting $f = f^* \restriction A_1$ the inequality \nl
$(h^e[{\Cal M}_n])^{- \varepsilon} \times
h^d_{A,D}[{\Cal M}_n] \le nu^k(f,f(A_0),B,D,{\Cal M}_n) \le 
h^u_{A,D}[{\Cal M}_n] \times (h^e[{\Cal M}_n])^\varepsilon$ 
(on nu$^k$ see below).
\ermn
2A) We say $(A^*,A,B,D)$ is semi$^{**}$-nice-$(k,r)$-good for
$({\frak K},\bar h,h)$ if: $A_1 = A^*$ in part (2). \nl
3) If $A \le_s D,B \le D,k \in \Bbb N$ and $f:A \rightarrow {\Cal M}$ we let
$nu^k(f,A,B,D,{\Cal M}) = |ex^k(f,A,B,D,{\Cal M}|$ where

$$
\align
\text{ex}^k(f,A,B,D,{\Cal M}) = \bigl\{g:&g \text{ is an embedding of } D 
\text{ into } M \\
  &\text{extending } f \text{ and satisfying } c \ell^k(g(B),M) = 
g(c \ell^k(B,D)) \bigr\}. 
\endalign
$$
\mn
4) We say that $({\frak K},\bar h,h)$ is polynomially semi-nice if a), b),
c), d) of part (1) holds and
\mr
\item "{$(e)$}"  $\bar h$ is polynomial over $h^e$.
\endroster
\sn
We can list some obvious implications.
\enddefinition
\bigskip

\proclaim{\stag{3.4D} Claim}  1) Assume $({\frak K},\bar h,h^e)$ is 
semi-nice and $h^e$ goes to infinity and 
$(A_1,A,B,D)$ is semi$^*-(k,r)$-good for
$({\frak K},\bar h,h)$.  \ub{Then} we can find $B',D',g^*$ such that
\mr
\item "{$(a)$}"  $A \le A_1 \le_s D'$ and $B' \le D' \in {\Cal K}_\infty$
\sn
\item "{$(b)$}"  $g^*$ is an embedding of $D$ into $D'$ 
\sn
\item "{$(c)$}"  $B' = g^*(B),D' = A_1 \cup g^*(D)$ and 
$c \ell^k(g^*(B),D') = g^*(c \ell^k(B,D))$
\sn
\item "{$(d)$}"  for every random enough ${\Cal M}_n$ and $f:A_1 \rightarrow
{\Cal M}_n$ satisfying $c \ell^r(f(A),{\Cal M}_n) \subseteq f(A_1)$ there
is $g':D' \rightarrow {\Cal M}_n$ extending $f$ such that
$c \ell^k(g'(B'),{\Cal M}_n) = g'(c \ell^k(B',D'))$, that is $(A_1,A,B',D')$
is semi$^{**}-(k,r)$-nice
\sn
\item "{$(e)$}"  if $({\frak K},\bar h,h^e)$ is polynomially semi-nice then
\nl
$\alpha(A,D) = \alpha(A_1,D')$.
\ermn
2) Assume ${\frak K}$ obeys $\bar h$ with error $h^e,h^e$ going to infinity.
If $(A,A_0,B,D)$ is semi$^*-(k,r)$-good for $({\frak K},\bar h,
h)$ \ub{then} it is semi-$(k,r)$-good (see Definition \scite{2.6G}(1)). \nl
3) If $({\frak K},\bar h,h^e)$ is semi-nice, \ub{then} ${\frak K}$ is
semi-nice.
\endproclaim
\bigskip

\demo{Proof}  1) Straight by counting. \nl
2) By part (1). \nl
3) By part (2).   \hfill$\square_{\scite{3.4D}}$
\enddemo
\bigskip

\proclaim{\stag{3.4E} Claim}  1) Assume ${\frak K}$ is semi-nice.  \ub{Then} 
for every $A \in {\Cal K}_\infty$ and $k$ and $\ell$, for some $r$ we have:
\mr
\item "{$(*)$}"  for every random enough ${\Cal M}_n$, for every $f:A
\rightarrow {\Cal M}_n$ and $B \le {\Cal M}_n,|B| \le \ell$, we have
\ermn
$$
(c \ell^r(f(A),{\Cal M}_n),f(A),B,c \ell^k(f(A) \cup B,{\Cal M}_n)) 
\text{ is semi-}(k,r)\text{-good}.
$$
\mn
2) Similarly for semi$^*$-nice for $({\frak K},\bar h,h)$.
\endproclaim
\bigskip

\demo{Proof}  1) We prove it by induction on $\ell$.  Now for $\ell = 0$ 
this is trivial by \scite{2.6I}, so let us prove it for $\ell + 1$ 
(assuming we have proved it for $\ell$).  So let 
$A \in {\Cal K}_\infty$ be given.   Let $r(1)$ be such that (exists by 
\scite{2.6J} applied to $k' \ge k$ and $\ell' = \ell + |A|$)
\mr
\item "{$(*)_1$}"  if ${\Cal M}_n$ is random enough, $A' \in {\Cal K}_\infty,
|A'| \le |A| + \ell,A' \le {\Cal M}_n$ and $b \in {\Cal M}_n$, then
$(c \ell^{r(1)}(A',{\Cal M}_n),A',A' \cup \{b\},c \ell^k(A' \cup \{b\},
{\Cal M}_n))$ is semi-$(k,r(1))$-good.
\ermn
Similarly by the induction hypothesis, for some $r(2)$
\mr
\item "{$(*)_2$}"  if ${\Cal M}_n$ is random enough, $A' \in {\Cal K}_\infty,
|A'| \le |A|,B' \le {\Cal M}_n,|B'| \le \ell$ then $(c \ell^{r(2)}(A',
{\Cal M}_n),A',A' \cup B',c \ell^{r(1)}(A \cup B',{\Cal M}_n))$ is
semi-$(r(1),r(2))$-good.
\ermn
We shall show that $r(2)$ is as required.  So let ${\Cal M}_n$ be 
random enough and $f:A \rightarrow {\Cal M}_n$ and
$B \subseteq {\Cal M}_n,|B| = \ell + 1$.  Let $B = B_0 \cup \{b\},|B_0| \le
\ell$.  So by $(*)_1$, the quadruple
$(c \ell^{r(1)}(f(A) \cup B_0,{\Cal M}_n),f(A) \cup B_0,B_0 \cup \{b\},
c \ell^k(A \cup B_0 \cup \{b\},{\Cal M}_n))$ is semi-$(k,r(1))$-good.  \nl
Similarly by $(*)_2$ the quadruple

$$
(c \ell^{r(2)}(f(A),{\Cal M}_n),f(A),B_0,c \ell^{r(1)}(f(A) \cup B_0,
{\Cal M}_n))
$$
\mn
is semi-$(r(1),r(2))$-good.
\sn

By ``transitivity" of the property easily the quadruple

$$
(c \ell^{r(2)}(f(A),{\Cal M}_n),f(A),B,c \ell^k(A \cup B,{\Cal M}_n))
$$
\mn
is $(k,r(2))$-good. \nl
2) Similar to the proof of part (1), using Definition \scite{3.4C} instead
of \scite{2.6J}. \nl
${{}}$  \hfill$\square_{\scite{3.4E}}$
\enddemo
\bigskip

\centerline {$* \qquad * \qquad *$}
\bn
We now turn to ``redrawing".
\definition{\stag{3.7} Definition}  Assume ${\frak K},{\frak K}^+$ are 0-1 
contexts. \newline
1)  We say ${\frak K}^+$ expands ${\frak K}$ if:
\medskip
\roster
\item "{(a)}"   $\tau^+$ a vocabulary extending $\tau$, (hence consisting 
of predicates only, $\tau^+$ locally finite, of course)
\sn
\item "{(b)}"   ${\Cal K}^+$ is the family of $\tau^+$-models satisfying $M^+ 
\restriction \tau \in {\Cal K}$ and \newline
${\Cal K}^+_n = \{ M^+ \in {\Cal K}^+:M^+ \restriction \tau \in {\Cal K}_n\}$
\sn
\item "{(c)}"   we let $\tau_{[\ell]} = 
\tau \cup \{R:R \in \tau^+ \text{ has } \ell \text{ places}\}$ and
$\tau_{< \ell} = \dsize \bigcup_{m < \ell} \tau_{[m]} \cup \tau$
\sn
\item "{(d)}"  for $M \in {\Cal K}_n$ we have \newline
$\mu_n(M) = \sum\{ \mu^+_n(M^+):M^+ \in {\Cal K}^+_n$ and 
$M^+ \restriction \tau = M\}$.
\endroster
\medskip

\noindent
For simplicity $\tau,\tau^+$ are irreflexive and $M \in {\Cal K}_n
\Rightarrow \mu_n(M)(= \mu_n(\{M\})) > 0$. \newline
2)  For ${\Cal M}_n \in {\Cal K}_n$ we define ${}^n \mu^+_{{\Cal M}_n}$, a
distribution on \newline
${\Cal K}^+_n[{\Cal M}_n] = \{ {\Cal M}^+_n \in {\Cal K}^+_n:
{\Cal M}^+_n \restriction \tau = {\Cal M}_n \}$ by \newline
$({}^n \mu^+_{{\Cal M}_n})
({\Cal M}^+_n) = \mu^+_n({\Cal M}^+_n)/\mu_n({\Cal M}_n)$, we write
$\mu^+_{{\Cal M}_n}$ when $n$ is clear from context (this is even formally
clear when the ${\Cal K}_n$'s are pairwise disjoint).
\medskip

We will be mostly interested in the case ${\Cal M}_n^+$ is drawn as in
Definition \scite{3.9}(2) below, but first define less general cases.
\enddefinition
\bigskip

\definition{\stag{3.8} Definition}  1) We define when 0-1 context 
${\frak K}^+$ is independently derived from ${\frak K}$ by the 
function $p$ (everything related to ${\frak K}^+$ has superscript $+$, 
below $\tau^+,{\Cal K}^+,{\Cal K}^+_n$ are as in 2.1(1), and for 
$x \in \{a,i,s,pr\},\le_x$ is defined by \scite{2.3}(2)).
\smallskip

\noindent
The crux of the matter is defining $\mu^+_n$; it suffices to define each
$\mu^+_{{\Cal M}_n}$.  We can think of it as choosing
a $\mu^+_n$-random model ${\Cal M}^+_n$ by expanding ${\Cal M}_n$,
defining ${\Cal M}^+_n \restriction \tau_{< \ell}$ by induction on $\ell$ by
fliping coins: for $\ell = 0$, ${\Cal M}^+_n \restriction \tau_{< 0}$
is chosen $\mu_n$-randomly from ${\Cal K}_n$ (i.e. is ${\Cal M}_n$).  
By induction on $\ell$, for each set $A \in [{\Cal M}_n]^\ell$ (i.e.
$A \subseteq {\Cal M}_n,|A| = \ell$): we choose
$A^+_\ell = ({\Cal M}^+_n \restriction \tau_{\le \ell}) \restriction A$, each
possibility $A^+_\ell$ has probability $p_{A^+_\ell,{\Cal M}^+_n \restriction
\tau_{< \ell}} = p_{A^+_\ell}[{\Cal M}^+_n \restriction \tau_{< \ell}]$ 
depending on ${\Cal M}^+_n \restriction \tau_{< \ell}$ and 
$({\Cal M}^+_n \restriction A) \restriction \tau_{< \ell}$ (not just on 
the isomorphism type), note that the
second one, $({\Cal M}_n^+ \restriction A) \restriction \tau_{< \ell}$ is
determined by $A^+_\ell$ as $A^+_\ell \restriction \tau_{< \ell}$.
Lastly, the drawings above (in stage $\ell$) are done independently for 
distinct $A$ (for each ${\Cal M}^+_n \restriction \tau_{< \ell}$).

2) We say ${\frak K}^+$ is derived uniformly
and independently if in addition $p_{A^+_\ell}[{\Cal M}^+_n \restriction
\tau_{< \ell}]$ depends on ${\Cal M}_n,A^+_\ell / \cong$ only and is 
derived very uniformly if it depends on $A^+ / \cong$ and $\|M_n\|$ 
(and $n$) only.
\enddefinition
\bigskip

\definition{\stag{3.9} Definition}  1) Suppose the 0-1 context 
${\frak K}^+$ is
independently derived from ${\frak K}$ by the function $p$ (see Definition 
\scite{3.8}).  We say $p$ has uniform bounds $\bar p$ if:
\medskip
\roster
\item "{(a)}"  $\bar p = (p^d,p^u)$
\sn
\item "{(b)}"  for $A^+ \in {\Cal K}^+$ and $\ell = |A^+| \in \Bbb N,
p^d_{A^+ \restriction \tau_{\le \ell}},p^u_{A^+ \restriction 
\tau_{\le \ell}}$ are functions from \nl
${\Cal K}_n$ to $[0,1]_{\Bbb R}$ (depending only on 
$A^+ \restriction \tau_{\le \ell}$ up to isomorphism)
such that for every random enough ${\Cal M}_n$:
\medskip
\noindent
{\roster
\itemitem{ $(*)$ }  for every embedding $f$ of $A = A^+ \restriction \tau$ 
into ${\Cal M}^+_n \restriction \tau$, \newline
letting $B^+ = f''(A^+)$ and $A^+_* = f(A^+)$ we have
\newline
$p^d_{A^+ \restriction \tau_{\le \ell}}[{\Cal M}_n] \le
p_{B^+ \restriction \tau_{\le \ell}}[{\Cal M}_n \restriction \tau_{< \ell}] 
\le p^u_{A^+ \restriction \tau_{\le \ell}}[{\Cal M}_n]$ \newline
and $p^x_{A^+ \restriction \tau_{\le \ell}}[{\Cal M}_n]$ depends on 
$(A^+ \restriction \tau_{\le \ell})/\cong$ only, (and, of course,
${\Cal M}_n$ and $n$); independently for the relevant distinct $A^+_*$'s. 
\newline
So we can write $p^x_{A^+ \restriction \tau_{\le \ell}}$ (for ${\Cal M}_n$) 
for $p^x_{A^+_* \restriction \tau_{\le \ell}}[{\Cal M}_n$].
\endroster}
\ermn
(Note that essentially $A^+ = A^+ \restriction \tau_{\le \ell}$ as
$|A^+| = \ell$ and the relation are assumed to be irreflexive, so we can waive
the $\restriction \tau_{\le \ell}$ abusing notation.
\nl
2)  Suppose the 0-1 context ${\frak K}^+$ is an expansion of ${\frak K}$.
We say that the drawing of ${\Cal M}_n^+$ (or for ${\frak K}^+$) 
\underbar{obeys} (the pair of functions) $\bar p = (p^d,p^u)$ with error
$h^e$ over ${\frak K}$ if for every $\varepsilon \in \Bbb R^{> 0}$ and
$A^+ \in {\Cal K}^+$, letting $\ell = |A^+|, 
A = A^+ \restriction \tau$ and given ${\Cal M}_n \in {\Cal K}_n$, and an
embedding
$f:A \rightarrow {\Cal M}_n$, assuming ${\Cal M}^+_n \restriction B$ was
already drawn for every $B \in [{\Cal M}_n]^{\le \ell}$ such that 
$B \ne f(A)$, and $f$ is an embedding of $A^+ \restriction \tau_{< \ell}$ 
to ${\Cal M}^+_n$ \ub{then} the probability (by $\mu^+_{{\Cal M}_n}$, modulo
the assumptions above) that $f$ embeds $A^+ \restriction \tau_{\le \ell}$ 
into ${\Cal M}^+_n$ is at least $p^d_{A^+ \restriction \tau_{\le \ell}}
[{\Cal M}_n] \times (h^e[{\Cal M}_n])^{- \varepsilon}$ and at most 
$p^u_{A^+ \restriction \tau_{\le \ell}}[{\Cal M}_n] \times (h^e[{\Cal M}_n])
^\varepsilon$ (so we assume always $p^d_{A^+ \restriction \tau_{\le \ell}}
[{\Cal M}_n] \le p^u_{A^+ \restriction \tau_{\le \ell}}[{\Cal M}_n]
\times (h^e[{\Cal M}_n])^\varepsilon$, at least for random enough 
${\Cal M}_n$.) \nl
3)  If $p^x_{A^+ \restriction \tau_{\le \ell}}[{\Cal M}_n]$ depends only on 
$\| {\Cal M}_n \|$ and $n$, we may write 
$p^x_{A^+ \restriction \tau_{\le \ell}}[\| M_n \|,n]$; if 
$\| {\Cal M}_n \|$ determines $n$ we may omit the latter (and when the 
intention is clear from context also in other cases). \newline
4)  For $A^+ \le^+ B^+ \in {\Cal K}^+_\infty$ and $x \in \{ d,u\}$ we let:

$$
p^x_{A^+,B^+}[{\Cal M}_n] =: \Pi \{ p^x_C[{\Cal M}_n]:C \le B^+ \text { and }
C \nsubseteq A^+\}.
$$
\mn
5) We omit $h^e$ if $h^e=1$, we say simply if $h^e[{\Cal M}_n] = \|
{\Cal M}_n\|$. \nl
6) Let $h_1,h_2,h^e$ be functions from $\dbcu_n {\Cal K}_n$ to 
$\Bbb R^{>0},\Bbb R^{\ge 1}$, respectively.  We say $h_1 
\underset h {}\to \sim h_2$ if for
every random enough ${\Cal M}_n,(h^e[{\Cal M}_n])^{-\varepsilon} \le h_1
[{\Cal M}_n]/h_2[{\Cal M}_n] \le (h^e[{\Cal M}_n])^\varepsilon$.
\enddefinition
\bigskip

\remark{\stag{3.10} Remark}  1) ``Obeys" (Definition \scite{3.9}(2)) 
means we have independence but only approximately, so we shall be 
able to give later other distributions in which the drawing are independent 
and which give lower and upper bounds to the situation for ${\frak K}$. \nl
2) Among those variants we use mainly Definition \scite{3.9}(2) and even more
the polynomial case.
\endremark
\bigskip

\definition{\stag{3.11} Definition}  In Definition \scite{3.9} 
we say that $\bar p$ is \underbar{polynomial} over $h$ if:
\medskip
\roster
\item "{(a)}"  $h$ is a function from $\dsize \bigcup_{n < \omega}
{\Cal K}_n$ to $\Bbb N$ converging to infinity
\sn
\item "{(b)}"   for every $\varepsilon \in \Bbb R^{>0}$ for every $\zeta
\in \Bbb R^{>0}$ for ${\Cal M}_n$ random enough \nl
$\zeta > h[{\Cal M}_n]/\|{\Cal M}_n\|^\varepsilon$
\sn
\item "{(c)}"  for every $A^+ = B^+ \restriction \tau_{\le \ell},\ell =
|B|,B^+ \in {\Cal K}^+$, for some $\beta(A^+) \in \Bbb R$ we have:
{\roster
\itemitem{ $(*)$ }  there are constants $c^d_{A^+},c^u_{A^+}$ such that:
\endroster}
\endroster

$$
p^d_{A^+}(M) = (c^d_{A^+})\|M\|^{\beta(A^+)}/h[M]
$$

$$
p^u_{A^+}(M) = (c^u_{A^+})\|M\|^{\beta(A^+)}/h[M].
$$
\mn
so this is not necessarily the very uniform case.
\enddefinition
\bigskip

\remark{Remark}  Of course, we can replace ``constant" by any slow enough
function.
\endremark
\bigskip

\proclaim{\stag{3.12} Claim}  1) Definition \scite{3.8}(1) is a particular 
case of Definition \scite{3.9}(2).  Also Definition \scite{3.8}(2), is a 
particular case of Definition \scite{3.9}(1), (with $p^u = p^d$), and 
all of them are particular cases of Definition \scite{3.7}(1).  
Also, if we have \scite{3.8}(1) + \scite{3.9}(2),then we have 
\scite{3.9}(1).  \newline
2)  In \scite{3.7}(1) necessarily
\medskip
\roster
\item "{(a)}"  ${\Cal K}_\infty = \{ A^+ \restriction \tau:A^+ \in 
{\Cal K}^+_\infty \}$
\sn
\item "{(b)}"  if $A^+ \subseteq B^+ \in {\Cal K}^+_\infty$ and $(A^+ 
\restriction \tau) \le_i (B^+ \restriction \tau)$ \underbar{then} $A^+ \le_i
B^+$
\sn
\item "{(c)}"  if $A^+ \subseteq B^+ \in {\Cal K}^+_\infty$ and $A^+ \le_s
B^+$ \underbar{then} $(A^+ \restriction \tau) \le_s (B^+ \restriction \tau)$.
\endroster
\endproclaim
\bigskip

\demo{Proof}  Check.
\enddemo


\newpage

\head {\S3 Regaining the context for ${\frak K}^+$: reduction} \endhead  \resetall
\bigskip

We write down the expected values for nu$(f,A^+,B^+,{\Cal M}^+_n)$ (see
Definition \scite{1.3}(2) and \scite{4.2} below) then we define $\le^+_b$ 
as what will be a variant of $\le^+_a$ if things are close
enough to the expected value and derive $\le^+_j$ (the parallel to
$\le^+_i$), $\le^+_t$ (the parallel to $\le_s$) and $\le^+_{\text{qr}}$ (the
parallel to $\le^+_{\text{pr}}$) (all done in Definition \scite{4.3}).  We
phrase a natural sufficient condition to (the probability condition part
for) $A^+ \le^+_a B^+$ (in \sciteu{4.5}) and show that when it is equivalent
to a natural strengthening then $A^+ <^+_b B \Rightarrow A^+ <^+_a B$ and
moreover $\le^+_j = \le^+_i,\le^+_t = \le^+_s$ (and some obvious fact, all 
this in \scite{4.6}).  We then prove that $\le^+_b$, has the formal 
properties of $\le^+_a,\dotsc,$ (in \scite{4.7}).

We then concentrate on the ``polynomial" case, ending with sufficient
conditions for ${\frak K}$ being semi-nice (\scite{4.8} - \scite{4.12}).
\bigskip

\demo{\stag{4.1} Context}  ${\frak K},\bar h$ as in \scite{3.1} and 
${\frak K}^+,\bar p$ as in \scite{3.9}(2) and we use \scite{3.9}(4).
\enddemo
\bigskip

\definition{\stag{4.1A} Definition}  Assume $f:A \rightarrow {\Cal M}_n$ (an
embedding),
$$
A_0 <_{pr} A_1 <_{pr} \cdots <_{pr} A_k \text{ and } A = A_0,B = A_k 
\text{ and } \bar A = \langle A_\ell:\ell = 0,\dotsc,k \rangle.
$$

\noindent
Let for $i \le k$:

$$
T^{[\ell]} = T^{[\ell]}(f,A,B,\bar A,{\Cal M}_n) =: \{g:g
\text{ an embedding of } A_\ell \text{ into } {\Cal M}_n \text{ extending }
f\}
$$

$$
T(f,A,B,\bar A,{\Cal M}_n) = \dsize \bigcup_{\ell \le k} T^{[\ell]}
(f,A,B,\bar A,{\Cal M}_n).
$$
\medskip

\noindent
in fact, we can omit $A,B$ as they are determined by $\bar A$.
\enddefinition
\bigskip

\proclaim{\stag{4.2} Claim}  Assume
\mr
\item "{$(\alpha)$}"  ${\frak K}$ obeys $\bar h$ with error $h^e_1$
\sn
\item "{$(\beta)$}"  the drawing of ${\Cal M}^+_n$ obeys $\bar p$ with error
$h^e_2$ (see \scite{3.9}(2))
\sn
\item "{$(\gamma)$}"  $\bar A = \langle A_\ell:\ell \le k \rangle$ is
a decomposition of $A <_s B$ (for ${\frak K}$)
\sn
\item "{$(\delta)$}"  $\varepsilon \in \Bbb R^{>0}$ and ${\Cal M}_n$ is
random enough
\sn
\item "{$(\varepsilon)$}"   $f:A \rightarrow {\Cal M}_n \in {\Cal K}_n$ an 
embedding, and $A^+ \restriction \tau = A,B^+ \restriction \tau = B,
A^+ \le B^+$ and $A^+_\ell =: B^+ \restriction A_\ell$
\sn
\item "{$(\zeta)$}"  $\varepsilon \in \Bbb R^{>0}$.
\ermn
\underbar{Then} the expected value of
nu$(f,A^+,B^+,{\Cal M}^+_n)$ under the distribution $\mu^+_{{\Cal M}_n}$
(see \scite{2.4}(1)) and the assumption ``$f$ will be an embedding of 
$A^+$ into ${\Cal M}^+_n$ and ${\Cal M}_n$ random enough" is:

$$
\text{at least } \dsize \prod_{\ell < k} \left( p^d_{A^+_\ell,A^+_{\ell + 1}}
[{\Cal M}_n] \times h^d_{A_\ell,A_{\ell + 1}}[{\Cal M}_n] \right) \times
(h[{\Cal M}_n])^{- \varepsilon}
$$

$$
\text{and at most }
\dsize \prod_{\ell < k} \left( p^u_{A^+_\ell,A^+_{\ell + 1}}
[{\Cal M}_n] \times h^u_{A_\ell,A_{\ell + 1}}[{\Cal M}_n] \right)
\times (h[{\Cal M}_n])^\varepsilon.
$$
\endproclaim
\bigskip

\demo{Proof}  Straight. \hfill$\square_{\scite{4.2}}$
\enddemo
\bigskip

\noindent
So by Claim \scite{4.2}, if the upper and lower bounds are close enough,
we can show that other decompositions of $(A_0,A_k)$ give similar results.
\newline
In the following, for the interesting case (here), $\le^+_j,\le^+_t,
\le^+_{\text{qr}}$ will be proved to be equal to $\le^{{\frak K}^+}_i,
\le^{{\frak K}^+}_s,\le^{{\frak K}^+}_{\text{pr}}$ respectively
(but $\le^+_b$ will not be $\le^{{\frak K}^+}_b$).
\definition{\stag{4.3} Definition}  1)  ${\Cal K}^\oplus_\infty = 
\{ A^+:A^+ \in {\Cal K}^+$ and $A^+ \restriction \tau \in 
{\Cal K}_\infty \}$ and, of course, $K^+_\infty = \{A^+:A^+ \in {\Cal K}^+$
and $0 < \lim \sup_n \text{ Prob}_{\mu^+_n}(A^+$ embeddable into
${\Cal M}^+_n)\}$. \newline
2)  For $A^+,B^+ \in {\Cal K}^\oplus_\infty$ let $A^+ \le^+_b B^+$ holds
\underbar{iff} $(A^+ \restriction \tau) \le_a (B^+ \restriction \tau)$ or
\newline
$(A^+ \restriction \tau) \le_s (B^+ \restriction \tau)$ and 
\medskip

\noindent
$\otimes^0_{A^+,B^+} \qquad \quad$  \underbar{for some} $k \in \Bbb N$, 
\ub{for every} $\zeta >0$ we have

$$
\align
1 = &\text{ Lim}_n \text{ Prob}_{\mu_n} \biggl( \text{we have } \zeta 
\text{ is larger than the probability that for some embedding} \\
  &f \text{ of } (A^+ \restriction \tau) \text{ into } {\Cal M}_n,
\text{ the number of extensions } g \text{ of} \\
  &f \text{ to embedding of } B^+ \text{ into } {\Cal M}^+_n 
\text{ is } \ge k, \text{ by the distribution } \mu^+_n[{\Cal M}^+_n] \\
  &\text{ under the assumption that } f \text{ embeds } A^+ \text{ into } 
{\Cal M}^+_n \biggr).
\endalign
$$
\medskip

\noindent
3)  $A^+ \le^+_j B^+$ \underbar{if} for every $A^+_1$, we have 
$A^+ \le^+ A^+_1 <^+ B \Rightarrow A^+_1 <^+_b B$. \newline
4)  For $A^+,B^+ \in {\Cal K}^\oplus_\infty$ let $A^+ \le^+_t B^+$ 
\underbar{if} $A^+ \le^+ B^+$ and for no $C^+ \in K^\oplus_\infty$ do we have
$A^+ <^+_j C^+ \le^+ B^+$. \newline
5)  For $A^+,B^+ \in {\Cal K}^\oplus_\infty$ let $A^+ \le^+_{qr} B^+$ if
$A^+ <^+_t B^+$ but for no $C^+$ do we have $A^+ <^+_t C^+ <^+_t B^+$.
\newline
6)  We say $\bar A^+$ is a $\le^+_{\text{qr}}$ decomposition of $A^+ <^+_s 
B^+$ if $\bar A^+ = \langle A^+_\ell:\ell \le k \rangle,A^+_\ell 
<^+_{\text{qr}} A^+_{\ell + 1},A^+_0 = A^+,A^+_k = B^+$. 
\newline
7)  ${\Cal K}^\otimes_\infty = \biggl\{ A^+:A^+ \in {\Cal K}^+$ and for some
$\bar A = \langle A_\ell:\ell \le k \rangle$ \newline

$\qquad \qquad \qquad$ and 
$\bar A^+ = \langle A^+_\ell:\ell \le k \rangle$ we have: \newline

$\qquad \qquad \qquad A_\ell = A^+_\ell \restriction
\tau,A_\ell <_{\text{qr}} A_{\ell + 1},\emptyset \le_i A_0$ and for some
$\varepsilon \in \Bbb R^{> 0}$ we have \newline

$\qquad \qquad \qquad 0 < \text{ Lim sup}_n \text{ Prob}
_{\mu_n} \biggl( \varepsilon < \dsize \prod_{\ell < k} (p^u_{A^+_\ell,
A^+_{\ell + 1}}[{\Cal M}_n] \times h^u_{A^+_\ell,A^+_{\ell + 1}}
[{\Cal M}_n]) \biggr) \biggr\}$.
\enddefinition
\bigskip

\proclaim{\stag{4.4A} Claim}  Assume $A^+ \subseteq B^+ \subseteq D^+,
A^+ \subseteq C^+ \subseteq D^+$ are in $K^+_\infty$ and $D^+ = B^+ \cup C^+$.
\nl
1) If $A^+ \le^+_b B^+$, \ub{then} $C^+ \le^+_b D^+$. \nl
2) If $A^+ \le^+_j B^+$, \ub{then} $C^+ \le^+_j D^+$.
\endproclaim
\bigskip

\demo{Proof}  1) Reflect. \nl
2) Follows from part (1).  \hfill$\square_{\scite{4.4A}}$
\enddemo
\bigskip

\proclaim{\stag{4.7} Claim}  1) ${\Cal K}^+_\infty \subseteq 
{\Cal K}^\otimes \subseteq {\Cal K}^\oplus_\infty \subseteq {\Cal K}^+$ are
closed under submodels and isomorphisms \newline
[why? reread Definition \scite{2.3}(1),\scite{4.3}(1), \scite{4.3}(7)]. \nl
2) If $A^+_1 \le^+_j A^+_2 \le^+ C^+ \in {\Cal K}^+_\infty$ and 
$B^+ \le^+ C^+ \in {\Cal K}^+_\infty$ and $A^+_1 \subseteq B^+$ 
\underbar{then} \newline
$B^+ \le^+_j B^+ \cup A^+_2$.  If $A^+ \le^+ B^+ \le^+ C^+,
A^+ \le^+_j C^+$ \underbar{then} $B^+ \le^+_j C^+$ and $B^+ \le^+_b C^+$
\newline
[why?  by Definition \scite{4.3}(3) and \scite{4.4A}]. \newline
3) On ${\Cal K}^\oplus_\infty$ the relation $\le^+_b$ is a partial order and
also the relation $\le^+_j$ is a partial order \nl
[why?  first we prove that $\le^+_b$ is a partial order so assume
$A^+_0 \le^+_b A^+_1 \le^+_b A^+_2 \in {\Cal K}^+_\infty$.  For $\ell = 1,2$
let $k_\ell \in \Bbb N^{>0}$ be such that the condition in
$\otimes^0_{A^+_\ell,A^+_{\ell+1}}$ of \scite{4.3}(2) holds.  Let $k =: k_0 +
k_1 \in \Bbb N$.  So now check. \nl
Suppose that $A^+ \le^+_j B^+ \le^+_j C^+$ but assume toward
contradiction that $\neg(A^+ \le^+_j C^+)$, so by the definition of $\le_j$
for some $D^+,A^+ \le^+ D^+ <^+ C^+,\neg(D^+ <^+_b C^+)$.  Now if $B^+ \le^+
D^+$ we get contradiction to $B^+ \le_j C^+$, so assume $\neg(B^+ \le D^+)$.
By monotonicity $A^+ \le^+_j B^+$ implies $D^+ <^+_j D^+ \cup B^+$ hence by
the definition of $\le^+_j$ we have $D^+ \le^+_b D^+ \cup B^+$.  Also as
$B^+ \le^+_j C^+$, necessarily $D^+ \cup B^+ \le^+_j C^+$, so as
$\le^+_b$ is transitive clearly $D^+ \le^+_b C^+$ as required.] \nl
4) If $A^+ \le^+ C^+$ are in ${\Cal K}^+_\infty$ \underbar{then} 
for one and only one $B^+ \in {\Cal K}^+_\infty$ we have \newline
$A^+ \le^+_j B^+ \le^+_t C^+$. \newline
[why?   let $B^+ \le^+ C^+$ be maximal 
such that $A^+ \le^+_j B^+$, it exists as $C^+$ is finite and 
$A^+ \le^+_j A^+$ (because $A \le_i A$ where $A = A^+ \restriction \tau$), 
now $B^+ \le^+_t C^+$ by \scite{4.7}(3) + Definition \scite{4.3}(4).  
Hence at least one $B^+$ exists, so suppose
$A^+ \le^+_j B^+_\ell \le^+_t C^+$ for $\ell = 1,2$ and $B^+_1 \ne B^+_2$ 
so without loss
of generality $B^+_2 \backslash B^+_1 \ne \emptyset$, by \scite{4.7}(2), 
$B^+_1 \le^+_j B^+_1 \cup B^+_2$ hence $B^+_1 <^+_j B^+_1 \cup B^+_2 
\le^+ C^+$, but this contradicts $B^+_1 \le^+_t C^+$ (see Definition 
\scite{4.3}(3))]. \newline
5)  If $A^+ <^+_t B^+$ \underbar{then} there is a 
$<_{\text{qr}}$-decomposition $\bar A^+$ of $A^+ <^+_t B^+$ [see 
Definition \scite{4.3}(5), remembering $B$ is finite]. \newline
6)  If $A^+ \le^+_t B^+$ and $C^+ \le^+ B^+$ \underbar{then} 
$C^+ \cap A^+ \le^+_t C^+$ \newline
[why? otherwise for some $C^+_1,C^+ \cap A^+ <^+_j C^+_1 \le^+ C^+$, 
\underbar{then} by \scite{4.7}(2) we have $A^+ <^+_j A^+ \cup C^+_1$, 
contradiction to $A^+ \le^+_t B^+$]. \newline
7)  The relations $\le^+_b,\le^+_j,\le^+_t,\le^+_{\text{qr}}$ are 
preserved by isomorphisms. \newline
8)  If $A^+ <^+_{\text{qr}} B^+$ \underbar{then} for every $b \in B^+ 
\backslash A^+$ we have $(A^+ \cup \{a\}) \le^+_j B^+$ \newline
[why? if not, then for some $C^+,(A^+ \cup \{a\}) \le^+ C^+ <^+_t B^+$, 
but \newline $A^+<^+_{\text{qr}} B^+ \Rightarrow A^+ <^+_t B^+ \Rightarrow 
A^+ <^+_t C^+$ 
(by Definition \scite{4.3}(5),\scite{4.7}(6) respectively)
so $A^+ <^+_t C^+ <^+_t B^+$ contradicting $A^+ <^+_{\text{qr}} B^+$]. \nl
9)  $T^+_\infty$ is a consistent (first order) theory.
\endproclaim
\bigskip

\proclaim{\stag{4.6} Claim}  Assume that ${\frak K}^+$ obeys $\bar p$ 
(over ${\frak K}$), ${\frak K}$ obeys $\bar h$ and
\mr
\item "{$(*)$}"  if $A^+,B^+ \in {\Cal K}^\oplus_\infty,(A^+ \restriction
\tau) <_s (B^+ \restriction \tau),\bar A = \langle A_\ell:\ell \le n \rangle$,
a decomposition so $A \restriction \tau = A_0 \le_s \ldots \le_s
A_n = B,A^+_\ell = B^+ \restriction A_\ell$ then $\bigotimes^0_{A^+,B^+}$ of
\scite{4.3}(2) and $\bigotimes^1_{A^+,B^+,\bar A}$ and $\bigotimes^2
_{A^+,B^+,\bar A}$ below are equivalent where
\mn

$\otimes^1_{A^+,B^+,\bar A}$ $\qquad \quad$  for some $\varepsilon > 0$ 
we have

$$
1 = \text{ Lim}_n \, \text{Prob}_{\mu_n}
\biggl( \varepsilon > \dsize \prod_{\ell < k} 
\left( p^u_{A^+_\ell,A^+_{\ell + 1}}[{\Cal M}_n] 
\times h^u_{A_\ell,A_{\ell + 1}}[{\Cal M}_n] \right) \biggr).
$$
\mn

$\otimes^2_{A^+,B^+,\bar A}$ $\qquad \quad$  for some $\varepsilon \in 
\Bbb R^{>0}$ we have \footnote{of course, we can think of cases that 
there are few copies
of $A^+_\ell$, then $\|{\Cal M}_n\|$ can be replaced by such upper bounds;
this has no influence in the polynomial case}

$$
1 = \text{ Lim}_n \, \text{Prob}_{\mu_n}
\biggl( \|{\Cal M}_n\|^{- \varepsilon} > \dsize \prod_{\ell < k}
(p^u_{A^+_\ell,A^+_{\ell + 1}}[{\Cal M}_n] \times h^u_{A_\ell,
A_{\ell + 1}}[{\Cal M}_n]) \biggr).
$$
\ermn
1) If $A^+ <^+_b B^+$ (in ${\Cal K}^\oplus_\infty$) 
\underbar{then} for some $m \in \Bbb N$ we have:

$$
\align
1 = \text{ Lim}_n \, \text{Prob}_{\mu^+_n} \biggl( &\text{ there is no }
f:A^+ \rightarrow {\Cal M}^+_n \text{ and } g_\ell:B^+ 
\rightarrow {\Cal M}^+_n \text{ for } \ell < m \\
  &\text{such that } \langle g_\ell:\ell < m \rangle \text{ is a disjoint
sequence of extensions of } f \biggr).
\endalign
$$
\medskip

\noindent
2) If $A^+ <^+_j B^+$ (in $K^\oplus_\infty$) \underbar{then} 
for some $k \in \Bbb N$ we have:

$$
\align
1 = \text{ Lim}_n \, \text{Prob}_{\mu^+_n} &\text{( there is no }
f:A^+ \rightarrow {\Cal M}^+_n \text{ and } g_\ell:B^+ \rightarrow {\Cal M}^+
\text{ for } \ell < k \\
  &\text{such that } \langle g_\ell:\ell < k \rangle \text{ is a sequence of 
distinct extensions of } f).
\endalign
$$
\medskip

\noindent
3) If  
$A^+ \le^+_b B^+ \in {\Cal K}^+_\infty$, \underbar{then} $A^+ \le^+_a B^+$. 
\newline
4) If $A^+ \le^+_j B^+ \in {\Cal K}^+_\infty$, 
\underbar{then} $A^+ \le^+_i B^+$. \newline
5)  If $(A^+ \restriction \tau) \le_a (B^+ \restriction \tau)$ and
$A^+ \le^+ B^+ \in {\Cal K}^+_\infty$, \ub{then} $A^+ \le^+_b B^+$. \nl
6) If $(A^+ \restriction \tau) \le_i (B^+ \restriction \tau)$ and
$A^+ \le B^+ \in {\Cal K}^+_\infty$, \underbar{then} $A^+ \le^+_j B^+$.
\endproclaim
\bigskip

\remark{\stag{4.6A} Remark}  1) Are there cases we may be interested in 
which are not covered by this claim?  If ${\frak K}$ obeys $\bar h$, and

$$
h^u_{A,B}(n) \sim n^{1/(\text{log} n)^{1/2}} \sim h^d_{A,B}(n).
$$

\noindent
2) We may rephrase the assumption in \scite{4.6} to cover 
those cases. \newline
3) Note: if all is polynomial, then $\otimes^2_{A^+,B^+,\bar A}$ 
is equivalent to $\otimes^1_{A^+,B^+,\bar A}$. \newline
4) We can use $h$ as in \scite{3.1}(3),\scite{3.1}(3A); see latter.
\endremark
\bigskip

\demo{Proof}  1) Let $\varepsilon > 0$ be as in $\otimes^2_{A^+,B^+,\bar A}$, 
and let $m \in \Bbb N$ be such that $\varepsilon m > |A^+|$ so 
$\zeta =: \varepsilon k - |A^+| \in
\Bbb R^{>0}$.  Let $\langle A_\ell:\ell \le k \rangle,\langle A^+_\ell:\ell
\le k \rangle$, be as in \sciteu{4.5} and let $A = A^+ \restriction \tau$; 
considering our aim, without loss of generality $A_0 = A$ 
(see Definition \scite{4.3}(1)).

So for every $M \in {\Cal K},\|M\|^{- \varepsilon m} \times |\{f:f \text{ an
embedding of } A \text{ into } M\}|$ is \newline
$\le \|M\|^{-(\varepsilon m - |A|)} = \|M\|^{- \zeta}$, 
hence it suffices to prove:
\mr
\item "{$(*)$}"  if ${\Cal M}_n$ is random enough and $f$ is an embedding
of $A$ into ${\Cal M}_n$ then \newline
$\|{\Cal M}_n\|^{- \varepsilon m} \ge \text{ Prob}_{\mu_n[M_n]} \biggl(
\text{there are a sequence of } k \text{ disjoint extensions of } f$\newline

$\qquad \qquad \qquad \qquad \qquad\text{   to an embedding of } 
B^+ \text{ into } {\Cal M}^+_n$, \newline

$\qquad \qquad \qquad \qquad \qquad \text{   under the assumption } f
\text{ is an embedding of}$ \newline

$\qquad \qquad \qquad \qquad \qquad \,A^+ \text{ into } {\Cal M}^+_n \biggr)$.
\endroster
\medskip

Now if ${\Cal M}_n$ is random enough then by \scite{3.6}(1) we know \newline
nu$(f,A,B,{\Cal M}_n) \le \dsize \prod_{\ell < k} h^u_{A_\ell,A_{\ell +1}}
[{\Cal M}_n]$, i.e. ex$(f,A,B,{\Cal M}_n)$ has $\le \dsize \prod_{\ell < k}
h^u_{A_\ell,A_{\ell + 1}}[{\Cal M}_n]$ members.  Let 
$F = F^m = F^m[{\Cal M}_n] = \{ \bar f:\bar f = \langle f_\ell:\ell < m 
\rangle$ and \nl
$f_\ell \in \text{ ex}(f,A,B,{\Cal M}_n)$ and $\bar f$ is disjoint 
over $A\}$. \newline
So $|F_m| \le |\text{ex}(f,A,B,{\Cal M}_n)|^m \le \biggl( \dsize
\prod_{\ell < k} h^u_{A^+_\ell,A^+_{\ell + 1}}[{\Cal M}_n] \biggr)^m$.
\newline
Hence, if we draw ${\Cal M}^+_n$ (by the distribution $\mu_n[{\Cal M}_n]$)
under the assumption ``$f$ is an embedding of $A^+$ into ${\Cal M}^+_n$", 
then the expected value of \newline
$|\{\bar f \in F^m[{\Cal M}_n]:\text{for } \ell < m,f_\ell
\text{ is an embedding of } B^+ \text{ into } {\Cal M}^+_n\}|$ is

$$
\align
\le \dsize \sum_{\bar f \in F} \text{ Prob}_{\mu^+_n[{\Cal M}_n]} \biggl(
  &\text{for } \ell < m,f_\ell \text{ is an embedding of } B^+ \text{ into }
{\Cal M}^+_n \biggl | \\
  &f \text{ is an embedding of } A^+ \text{ into } {\Cal M}_n \biggr)\\
  &\le \dsize \sum_{\bar f \in F} \biggl( \dsize \prod_{\ell < k}
p^u_{A^+_\ell,A^+_{\ell + 1}} [{\Cal M}_n] \biggr)^m \\
  &= |F| \times \biggl( \dsize \prod_{\ell < k}
p^u_{A^+_\ell,A^+_{\ell + 1}} [{\Cal M}_n] \biggr)^m \\
  &\le \biggl( \dsize \prod_{\ell < k}
h^u_{A^+_\ell,A^+_{\ell + 1}} [{\Cal M}_n] \biggr)^m \times
\biggl( \dsize \prod_{\ell < k}
p^u_{A^+_\ell,A^+_\ell} [{\Cal M}_n] \biggr)^m \\
  &= \biggl( \dsize \prod_{\ell < k} \biggl(
p^u_{A^+_\ell,A^+_{\ell + 1}} [{\Cal M}_n] \times
h^u_{A^+_\ell,A^+_{\ell + 1}} [{\Cal M}_n] \biggr) \biggr)^m
\endalign
$$
\medskip

\noindent
but if ${\Cal M}_n$ is random enough $\otimes^2_{A^+,B^+,\bar A}$ (i.e. by the
assumption of \scite{4.6} and $\otimes^0_{A^+,B^+}$ which holds as we are
assuming $A^+ \le^+_j B^+$) this last number is \nl
$\le \biggl( \| {\Cal M}_n\|^{- \varepsilon} \biggr)^m =
\|{\Cal M}_n\|^{- \varepsilon m}$. \newline

As said above, this suffices.  \newline
2) By the $\triangle$-system lemma and \scite{4.6}(1) (and the 
definitions). \newline
3) Follows by \scite{4.6}(1).  \newline
4) Follows from \scite{4.6}(2). \newline
5) Read the definitions.  \hfill$\square_{\scite{4.6}}$
\enddemo
\bigskip

\proclaim{\stag{4.6B} Claim}  1) A sufficient condition for $\le^+_i =
\le^+_j \restriction {\Cal K}^+_\infty$ is:
\mr
\item "{$(*)$}"  like $(*)$ of \scite{4.6} but we add another equivalent
condition
\mn
$\otimes^{1,d}_{A^+,B^+,\bar A}$ $\qquad \quad$  for some $\varepsilon > 0$ 
we have

$$
1 = \text{ Lim}_n \, \text{Prob}_{\mu_n}
\biggl( \varepsilon > \dsize \prod_{\ell < k} 
\left( p^d_{A^+_\ell,A^+_{\ell + 1}}[{\Cal M}_n] 
\times h^d_{A_\ell,A_{\ell + 1}}[{\Cal M}_n] \right) \biggr).
$$
\ermn
2) Note that $\bigotimes^2_{A^+,B^+,\bar A} \Rightarrow \bigotimes^1
_{A^+,B^+,A^+} \Rightarrow \bigotimes^{1,u}_{A^+,B^+,\bar A}$ and \nl
$\bigotimes^2_{A^+,B^+,A^+} \Rightarrow \bigotimes^0_{A^+,B^+,\bar A}
\Rightarrow \bigotimes^{1,u}_{A^+,B^+,\bar A}$.
\endproclaim
\bigskip

\demo{Proof}  By \scite{4.6}(4), it suffices to prove, assuming
$A^+ \le^+ B^+ \and \neg(A^+ \le^+_j B^+)$, that $\neg(A^+ \le^+_i B^+)$.
As $\neg(A^+ \le^+_j B^+)$ necessarily there is $A^+_1$ such that
$A^+ \le^+ A^+_1 \le^+ B^+,\neg(A^+_1 \le^+_b B^+)$.

The rest is easy, too.  \hfill$\square_{\scite{4.6B}}$
\enddemo
\bigskip

\definition{\stag{4.6C} Definition}  1) Assume
\mr
\item "{$(*)$}"  ${\frak K}$ obeys $\bar h$ with error $h^e,{\frak K}^+$ is
drawing obeys $\bar p$ with error $h^e_1$.
\ermn
We say that $(A^+,A^+_0,B^+,D^+)$ is a pretender to semi$^*-(k,r)$-good
for \nl
$({\frak K}^+,{\frak K},\bar h,h^e_0,\bar p,h^e_1)$ if: for some witness
$A^+_1$,
\mr
\item "{$(a)$}"  $A^+_0 \le^+_j A^+_1 \le^+_j A^+,
A^+_1 \le^+_t D^+,B^+ \le^+ D^+$
\sn
\item "{$(b)$}"  $(A^+ \restriction \tau,A^+_0 \restriction \tau,B^+
\restriction \tau,D^+ \restriction \tau)$ is semi$^*-(k,r)$-good for
$({\frak K},\bar h,h)$ witnessed by $A^+_1 \restriction \tau$
\sn
\item "{$(c)$}"  if $D^+ <^+ D^+_1 \in {\Cal K}^\oplus_\infty,A^+_1 <_t
D^+_1$ and $D^+ \restriction \tau <_i D^+_1 \restriction \tau$ and
$C^+_1 <^+_j C^+_2 \le D^+_1$ and $C^+_1 \le^+ A^+_1$ and $|C^+_2| \le k$
and $h^+_{A^+_0,D^+_1} \underset h^+ {}\to \sim h^+_{A^+_0,D^+}$, \ub{then}
$C^+_2 \le^+ D^+$.
\ermn
2) We in part one write semi$^{**}$ if $A^+_1 = A^+$.
\enddefinition
\bigskip

\proclaim{\stag{4.6D} Claim}  1) Assume ${\frak K}$ obeys $\bar h$ which is
bounded by $h,K^+$ drawn by $p$ which obeys $\bar p$ which is bounded by
$h$.  Sufficient conditions for $(*)_1 + (*)_2 + (*)_3$ are
$\otimes_1 + \otimes_2$ where
\mr
\item "{$(*)_1$}"   $({\Cal K}^+_\infty,\le^+_j \restriction 
{\Cal K}^+_\infty,\le^+_t \restriction 
{\Cal K}_\infty,\le^+_{qr} \restriction {\Cal K}_\infty) = 
({\Cal K}^+_\infty,\le^+_i,\le^+_s,\le^+_{pr})$
\sn
\item "{$(*)_2$}"  ${\frak K}^+$ is weakly nice
\sn
\item "{$(*)_3$}"  ${\Cal K}^+$ obeys $\bar h^+$ with error $h^e$ 
\ermn
and 
\mr
\item "{${}$}"  $\bigotimes_1 \quad \bar h^+ = 
\langle {}^+h^u,{}^+h^d \rangle$, if $A^+ <^+_{qr} B^+$ then 
{\roster
\itemitem{ ${}$ }  $(\alpha) \quad {}^+h^u_{A^+,B^+}[{\Cal M}^+_n] =$ \nl 

$\qquad \qquad {}^+h^u_{A^+,B^+}[{\Cal M}_n] = h^u_{A^+,B^+}[{\Cal M}^+_n 
\restriction \tau] \times p^u_{A^+,B^+}[{\Cal M}^+_n \restriction \tau]$ 
\sn
\itemitem{ ${}$ } $(\beta) \quad {}^+h^d_{A^+,B^+}[{\Cal M}^+_n] = 
{}^+h^d_{A^+,B^+}[{\Cal M}_n] =$ \nl

$\qquad \qquad h^d_{A^+,B^+}[{\Cal M}^+_n \restriction
\tau] \times p^d_{A^+,B^+}[{\Cal M}^+_n \restriction \tau]$
\sn
\itemitem{ ${}$ } $(\gamma) \quad h^e[{\Cal M}^+_n] \ge h_1[{\Cal M}^+_n 
\restriction \tau] \times h_2[{\Cal M}^+_n \restriction \tau]$ 
goes to infinity
\endroster}
\item "{${}$}"  $\bigotimes_2 \quad$ for 
$A^+ <_{qr} B^+$ and $\varepsilon \in \Bbb R^{>0}$, for 
random enough ${\Cal M}_n$ we have \nl

$\quad$ if $f$ embeds $A^+ \restriction \tau$ into ${\Cal M}_n$, under the
assumption that $f$ embeds \nl

$\quad A^+$ into ${\Cal M}^+_n$, the probability
that nu$(f,A^+,B^+,{\Cal M}^+_n)$ is not in the \nl

$\quad$ interval
$({}^+h^d_{A^+,B^+}[{\Cal M}_n] \times (h^e[{\Cal M}_n])^{- \varepsilon},
{}^+h^d_{A^+,B^+}[{\Cal M}^+_n] \times (h^e([{\Cal M}_n])^\varepsilon)$ \nl

$\quad$ is $< 1/\|{\Cal M}_n\|^{1/\varepsilon}$
\sn
\item "{${}$}"  $\bigotimes_3 \quad$  for $A^+ <^+_{qr} B^+,\varepsilon \in
\Bbb R^{>0}$ and $m \in \Bbb N$ for random enough ${\Cal M}_n$ and \nl

$\qquad f:A^+ \restriction \tau \rightarrow {\Cal M}_n$ we have

$$
{}^+h^d_{A^+,B^+}[{\Cal M}_n]/(h^e_3[{\Cal M}_n])^\varepsilon \ge m.
$$
\ermn 
2) In part (1) if in addition $\otimes^+_2$ then in addition $(*)_4 + (*)_5$
\mr
\item "{$(*)_4$}"  if $(A^+,A^+_0,B^+,D^+)$ is a pretender to being
semi$^*-(k,r)$-good for \nl
$({\frak K}^+,{\frak K},\bar h,h^e_0,\bar p,h^e_1)$, 
\ub{then} it is semi$^*-(k,r)$-good for $({\frak K}^+,\bar h^+,h^+)$
\sn
\item "{$(*)_5$}"  ${\Cal K}^+$ is semi-nice, morever even 
$({\Cal K}^+,\bar h^+,h^+)$ is semi-nice, ${\Cal M}^+_n$ random enough
\ub{and}
\sn
{\roster
\itemitem{ $\bigotimes^+_2$ }  if $(A^+,A^+_0,B^+,D^+)$ is a pretender to
semi$^*-(k,r)$-good, \ub{then} it is $f:A^+ \rightarrow {\Cal M}^+_n$ and
$[C^+_1 <^+_j C^+_2 \le {\Cal M}^+_n \and C^+_1 \le f(A) \and |C^+_2| \le r
\Rightarrow C^+_2 \subseteq f(A^+)]$, \ub{then} the inequalities in
$\bigotimes_2$ holds for nu$^k(f,A^+,B^+,D^+)$.
\endroster}
\endroster
\endproclaim
\bigskip

\demo{Proof}  Straight.
\enddemo
\bigskip

\centerline {$* \qquad * \qquad *$}
\bigskip

\noindent
Now we deal with the polynomial case.
\definition{\stag{4.8} Definition}  Assume
\medskip
\roster
\item "{$(*)_{\text{ap}}(a)$}"  ${\frak K}$ is a 0-1 law context
\sn
\item "{$(b)$}"     ${\frak K}$ obeys $\bar h = (h^d,h^u)$ (see Definition
\scite{3.1})
\sn
\item "{$(c)$}"     $\bar h$ is polynomial over $h_1$ with 
the function $\alpha(-,-)$ (see Definition \scite{3.3})
\sn
\item "{$(d)$}"     the 0-1 law context ${\frak K}^+$ is an expansion of
${\frak K}$ obeying the pairs of functions $\bar p = (p^d,p^u)$ (see
Definition \scite{3.9}(2))
\sn
\item "{$(e)$}"     $\bar p$ is polynomial over $h_2$ (see Definition 
\scite{3.11}) by the function $\beta(-)$. \newline
Here let $A^+ \le^- B^+$ mean $A^+ \subseteq B^+ \in {\Cal K}^\otimes_\infty$.
\sn
\item "{$(f)$}"  for $\varepsilon > 0$, for random enough 
${\Cal M}_n,h_1[{\Cal M}_n],h_2[{\Cal M}_n] \le \|{\Cal M}_n\|^\varepsilon$.
\endroster
\medskip

\noindent
1) For $A^+ \le^+ B^+$ satisfying \newline
$(A^+ \restriction \tau) <_s (B^+ \restriction \tau)$ we define \newline
$\beta(A^+,B^+) =: \alpha(A^+ \restriction \tau,B^+ \restriction \tau)
+ \sum\{\beta(C^+):C^+ \subseteq B^+,C^+ \nsubseteq A^+\}$. \newline
\sn
2)  For $A^+ \le^+ B^+$ let $A^+ \le^{+p}_b B^+$ mean: \newline
for some $A^+_1$ we have:
\medskip
\roster
\widestnumber\item{$(iii)$}
\item "{$(i)$}"  $A^+ \le^+ A^+_1 \le^+ B^+$
\sn
\item "{$(ii)$}"  $(A^+ \restriction \tau) \le_i (A^+_1 \restriction \tau)
\le_s (B^+ \restriction \tau)$
\sn
\item "{$(iii)$}"  $A^+ \ne A^+_1$ \underbar{or} for some $A^+_2$ we
have $A^+_1 <^+ A^+_2 \le^+ B^+$ and \newline
$\beta(A^+_1,A^+_2) < 0$.
\endroster
\medskip

\noindent
3)  For $A^+ \le^+ B$ let $A^+ \le^{+p}_j B^+$ mean: \newline
for every $A^+_1$ we have $A^+ \le^+ A^+_1 <^+ B^+ \Rightarrow A^+_1 <^{+p}_b
B^+$. \newline
\sn
4)  For $A^+ \le^+ B$ let $A^+ \le^{+p}_t B^+$ mean: \newline
for no $A^+_1$ do we have $A^+ <^{+p}_j A^+_1 \le^+ B$. \newline
\sn
5)  For $A^+ \le^+ B$ let $A^+ \le^{+p}_{qr} B^+$ mean: \newline
$A^+ \le^+_t B^+$ but for no $A^+_1$ do we have $A^+ <^{+p}_t A^+_1 <^{+p}_t
B^+$. \newline
\sn
6)  ${\Cal K}^{+p}_\infty = \biggl\{ A^+ \in {\Cal K}^+:\text{letting }
(A^+ \restriction \tau \restriction \emptyset) \le_i A_1 <_s A^+ \restriction
\tau \text{ and } A^+_1 =: A^+ \restriction A_1$, \newline

$\qquad \qquad \qquad \qquad \,\text{ \,\,we have: } [A^+_2 \subseteq A^+_1
\Rightarrow \beta(A^+_2) = 0] \text{ and } \beta(A^+_1,A^+) \ge 0$ and
\newline

$\qquad \qquad \qquad \qquad \,\,\,
\,\,0 < \text{ Lim inf}_n \text{ Prob}_{\mu_n}$ 
(we can embed $A^+_1$ into ${\Cal M}_n)\biggr\}$.
\enddefinition
\bigskip

\noindent
\underbar{\stag{4.9A} Discussion}:  Note that there can be $A^+$ 
such that the sequence \newline
$\langle \text{Prob}_{\mu_n}(A^+ \text{ embeddable into } {\Cal M}^+_n):
n \in \Bbb N \rangle$ does not converge, but essentially this occurs only when
$((A^+ \restriction \emptyset) \restriction \tau) \le_i (A^+ \restriction 
\tau)$.
More exactly if $A^+ \restriction \emptyset \restriction \tau \le_i 
(A^+_1 \restriction \tau) \le_s (A^+ \restriction \tau)$ assuming the answer 
for $(A^+_1$ embeddable into ${\Cal M}_n)$, we almost surely know if $A^+$ 
can be embedded into ${\Cal M}^+_n$.
\bigskip

\noindent
We first note
\demo{\stag{4.10} Fact}  Assume $(*)_{\text{ap}}$ of \scite{4.8} and
\medskip
\roster
\item "{$\otimes_1$}"  the irrationality assumption for $({\frak K},\bar h,
{\frak K}^+,\bar p)$ which means: if $A^+ <^{+p}_{qr} B^+$, \newline
$(A^+ \restriction \tau) <_s (B^+ \restriction \tau)$ \underbar{then}
$\beta(A^+,B^+) \ne 0$.
\endroster
\medskip

\noindent
1) If $A^+ \subseteq B^+ \subseteq C^+$ and $(A^+ \restriction \tau) <_s
(B^+ \restriction \tau) <_s (C^+ \restriction \tau)$ then \newline
$\beta(A^+,B^+) + \beta(B^+,C^+) = \beta(A^+,C^+)$, (of course, $\otimes_1$
is not used). \newline
2) If $A^+ \le^+ B^+ \in {\Cal K}^\oplus_\infty$ and
$(A^+ \restriction \tau) <_s (B^+ \restriction \tau)$ and 
$\beta(A^+,B^+) < 0$ \underbar{then} for some $m$ for every random
enough ${\Cal M}^+_n$ and any embedding $f:A^+ \rightarrow {\Cal M}_n^+$,
there are no $m$ disjoint extensions $g:B^+ \rightarrow {\Cal M}_n^+$ of $f$. 
\newline
2A) The four conditions in $(*)$ of \scite{4.6B} (see also $(*)$ of
\scite{4.6}) are equivalent (for any $A^+,B^+,\bar A$ as there). \nl
3) $\le^{+p}_x = \le^+_x$ for $x \in \{b,j,t,qr\}$ and ${\Cal K}^+_\infty =
{\Cal K}^\otimes_\infty = {\Cal K}^{\otimes_p}_\infty$. \newline
4) Assume $A^+ \subseteq B^+ \subseteq D^+,A^+ \subseteq C^+ \subseteq D^+,
D^+ = B^+ \cup C^+,B^+ \cap C^+ = A^+,A^+ <^{+p}_t B^+$ and
$(C^+ \restriction \tau) \le_s (D^+ \restriction \tau)$ \underbar{then}
$\beta(A^+,B^+) \ge \beta(C^+,D^+)$. 
\enddemo
\bigskip

\demo{Proof}  1) \newline

$$
\align
\beta(A^+,C^+) &= \alpha(A^+,C^+) + \sum \{\beta(D^+):D^+ \subseteq
C^+ \text{ and } D^+ \nsubseteq A^+\} \\
  &= \alpha(A^+,B^+) + \alpha(B^+,C^+) + \sum \{\beta(D^+):D^+ \subseteq
C^+ \text{ and } D^+ \nsubseteq A^+\} \\
  &= \alpha(A^+,B^+) + \alpha(B^+,C^+) + \sum \{\beta(D^+):D^+ \subseteq
B^+,D^+ \nsubseteq A^+\} \\
  &+ \sum \{\beta(D^+):D^+ \subseteq C^+,D^+ \nsubseteq B^+\} \\
  &= \beta(A^+,B^+) + \beta(B^+,C^+).
\endalign
$$

\noindent
2) For ${\Cal M}_n$ which is random enough and $f:(A^+ \restriction \tau)
\rightarrow {\Cal M}_n$, the set \newline
$Y = ex(f,A,B,{\Cal M}_n)$ has
$\le h^u_{A,B}({\Cal M}_n)$ members (see Definition \scite{3.1}).
Hence the set
$Y_m = \{\bar g:\bar g = \langle g_\ell:\ell < m \rangle \text{ is a 
sequence of members of } Y \text{ disjoint over } f\}$ has \newline
$\le (h^u_{A,B}({\Cal M}_n))^m$ members which is 
$\le (h_1[{\Cal M}_n])^{t_1} \cdot \|M_n\|^{\alpha(A,B)m}$.  Now under 
the assumption $f:A^+ \rightarrow {\Cal M}_n^+$, for each $\bar g \in Y_m$, 
the probability that $\dsize \bigwedge_{\ell < m}$ $(g_\ell$ embeds $B^+$ 
into ${\Cal M}_n^+$) for appropriate $t_2 \in \Bbb N$ is:

$$
\le \bigl( h_2[{\Cal M}_n]^{t_2} \cdot \|{\Cal M}_n\|
^{\Sigma\{\beta(C^+):C^+ \subseteq B^+,C^+ \nsubseteq A^+\}} \bigr)^m.
$$
\mn
So the expected value is

$$
\align 
\le h_1[{\Cal M}_n]^{t_1} &\times h_2[{\Cal M}_n]^{mt_2} \\
  &\times \bigl(\|{\Cal M}_n\|^{\alpha(A^+ \restriction \tau,B^+ \restriction
\tau) + \Sigma\{\beta(C^+):C^+ \subseteq B^+,C^+ \nsubseteq A^+\}} \bigr)^m \\
  &= h_1[{\Cal M}_n]^{t_1} \times h_2[{\Cal M}_n]^{mt_2} \times 
\|{\Cal M}_n\|^{m \beta(A^+,B^+)}.
\endalign
$$
\mn
So as $\beta(A^+,B^+) < 0$, and the assumptions on $h_1,h_2$ we have:
for $m$ large 
enough, this probability is $< \|{\Cal M}_n\|^{|A|+1}$ so the conclusion 
follows. \newline
2A)  Straight. \nl
3)  Assume $A^+ \le^- B^+$ and we shall prove that $A^+ \le^+_b B^+
\Leftrightarrow A^+ \le^{+p}_b B^+$. \newline
Let $A^+_1$ be such that $A^+ \le^+ A^+_1 \le^+ B^+$ and
$(A^+ \restriction \tau) \le_i (A^+_1 \restriction \tau) \le_s
(B^+ \restriction \tau)$.  Now if $A^+ \ne A^+_1$ then both 
$A^+ \le^+_b B^+,A^+ \le^{+p}_b B^+$ hold and we are done so assume
$A^+ = A^+_1$.  Now compare condition $\otimes^0_{A^+,B^+}$ of \scite{4.3}(2) 
and $\beta(A^+,B^+) < 0$, (as in the proof of part (2)) they are the same.
So $\le^{+ p}_b = \le^+_b$, and the other equalities follow. \newline
4)  Let ${\Cal M}_n$ be random enough but such that some $f_1$ embeds
$C^+ \restriction \tau$ into it (note $C^+ \in {\Cal K}^{+p}_\infty$ and
see Definition \scite{4.8}(7)).  Assume $f_1$ embeds $C^+$ into 
${\Cal M}^+_n$ and
let $f_0 = f \restriction A^+$; now the expected value of the number of
$g_0:B^+ \rightarrow {\Cal M}^+_n$ is $\approx \|{\Cal M}_n\|^{\beta(C^+,
D^+)}$, so the desired inequality follows. 
\hfill$\square_{\scite{4.10}}$
\enddemo
\bigskip

\proclaim{\stag{4.11} Claim}  1) Assume $(*)_{\text{ap}}$ of Definition 
\scite{4.8}. \newline
A sufficient condition for $(*)_1 + (*)_2 + (*)_3$ below is 
$\otimes_1 + \otimes_2 + \otimes_3$ below where:
\medskip
\roster
\item "{$(*)_1$}"  $({\Cal K}^{+p}_\infty,
\le^{+p}_j \restriction {\Cal K}^+_\infty,\le^{+p}_t \restriction 
{\Cal K}_\infty,\le^{+p}_{qr} \restriction {\Cal K}_\infty)$ \newline
$= ({\Cal K}^+_\infty,\le^+_i,\le^+_s,\le^+_{pr})$
\smallskip
\noindent
\item "{$(*)_2$}"  ${\frak K}^+$ is weakly nice
\smallskip
\noindent
\item "{$(*)_3$}"  ${\frak K}^+$ obeys a pair $\bar h$ which is 
polynomial over some $h^e$ \newline
with the function $(A^+,B^+) \mapsto \beta(A^+,B^+)$ when $A^+ <^+_s B^+$.
{\roster
\itemitem{ $\otimes_1$ }   for some function 
$h_3:{\Cal K} \rightarrow \Bbb N$ satisfying
$\varepsilon \in \Bbb R^+ \Rightarrow$ for random enough \nl

$\qquad {\Cal M}_n$ we have 
$[h_3({\Cal M}_n)/\|{\Cal M}_n\|^\varepsilon] = 0$, 
($h_3$ somewhat above) \nl

$\qquad h_1[{\Cal M}_n] \times h_2({\Cal M}_n]$
\sn
\itemitem{ $\otimes_2$ }  the irrationality assumption: for 
$({\frak K},\bar h,{\frak K}^+,\bar p)$: if $A^+ <^{+p}_{qr} B^+$, \newline
$(A^+ \restriction \tau) <_s (B^+ \restriction \tau)$ 
\underbar{then} $\beta(A^+,B^+) \ne 0$
\smallskip
\noindent
\itemitem{ $\otimes_3$ }  if $A^+ <_{qr} B^+$, \ub{then} for every 
${\Cal M}^+_n$ random enough, for every \nl
$f_0:A^+ \rightarrow {\Cal M}^+_n$ we have: \newline
$\|{\Cal M}_n\|^{\beta(A^+,B^+)}/h_3[\|{\Cal M}_n\|] \le 
|\{f_1:f_1 \text{ is an embedding of } B^+ \text{ into}$ \nl

$\qquad {\Cal M}^+_n \text{ extending } f_0\}| \le 
\|{\Cal M}_n\|^{\beta(A^+,B^+)} \times h_3[\|{\Cal M}_n\|]$.
\endroster}
\endroster
\mn
2) In part (1) we add $(*)_4, (*)_5$ if we assume also $\otimes^+_3$
\mr
\item "{$(*)_4$}"  if $(A^+,A^+_0,B^+,D^+)$ is a pretender to 
semi$^*-(k,r)$-good, \ub{then} it is semi-$(k,r)$-good
\sn
\item "{$(*)_5$}"  $({\frak K}^+,\bar h^+,h^+)$ is semi-$(k,r)$-good
\ermn
where
\mr
\item "{$\bigotimes^+_3$}"  if $(A^+,A^+_0,B^+,D^+)$ is a pretender to
semi-$(k,r)$-good, \ub{then} for every random enough ${\Cal M}_n$ and
$f:A^+ \rightarrow {\Cal M}_n$ letting $f_0 = f \restriction A^+_0$ we
have
$$
\align
\|{\Cal M}_n\|^{\beta(A^+_0,D^+)}/h_3({\Cal M}_n) &\le \text{ nu}^k(f_0,
A_0,B,D,{\Cal M}_n) \\
  &\le \|{\Cal M}_n\|^{\beta(A^+_0,D)} \times h_3({\Cal M}_n).
\endalign
$$
\endroster
\endproclaim
\bigskip

\demo{Proof}  Straightforward.
\enddemo
\bigskip

\demo{\stag{4.17} Conclusion}  Assume $(*)_{ap}$ of \scite{4.8} and
$\bigotimes_1$ of \scite{4.10} (the irrationality) and for simplicity let
$h^e(M) = \|M\|$. \nl
\ub{Then}
\mr
\item "{$(a)$}"  the conditions of \scite{4.11}(1), \scite{4.11}(2) hold
hence their conclusions
\sn
\item "{$(b)$}"   for $A^+ \le B^+ \in {\Cal K}^\oplus_\infty$, we
have 
{\roster
\itemitem{ (i) }  $A^+ \le B^+ \in {\Cal K}^+_\infty$ and
\sn
\itemitem{ (ii) }  $A^+ \le^+_b B^+ \Rightarrow A^+ \restriction \tau
\le_a B^+ \restriction \tau$
\sn
\itemitem{ (iii) }  $A^+ \le^+_j B^+ \Leftrightarrow 
A^+ \le^+_i B^+ \Leftarrow A^+ \restriction \tau
\le_i B^+ \restriction \tau$ 
\sn
\itemitem{ (iv) }  $A^+ \le^+_t B^+ \Leftrightarrow A^+ \le^+_s B^+
\Leftarrow A^+ \restriction \tau \le_s B^+ \restriction \tau$ 
\sn
\itemitem{ (v) }  $A^+ \le^+_b B^+ \Leftrightarrow A^+ \le^+_a B^+
\Leftarrow A^+ \restriction \tau \le_a B^+ \restriction \tau$
\endroster}
\item "{$(c)$}"  if $A^+ <^+_{qr} B^+$, then for every random enough
${\Cal M}_n$ for every \nl
$f^+:A^+ \rightarrow {\Cal M}^+_n$, we have, for some constant $\bold c$

$$
\align
h^d_{A^+ \restriction \tau,B^+ \restriction \tau}[{\Cal M}_n] \times
p^d_{A^+,B^+}[{\Cal M}_n] &- \sqrt{h^d_{A^+ \restriction \tau,
B^+ \restriction \tau}[{\Cal M}_n] \times p^d_{A^+,B^+}[{\Cal M}_n]} \\
  &\times \bold c \times (\text{log}\|{\Cal M}_n\|) \\
  &\le \text{ nu}(f^*,A^+,B^+,{\Cal M}^+_n) \\
  &\le h^u_{A^+ \restriction \tau,B^+ \restriction \tau}[{\Cal M}_n] \times
p^u_{A^+,B^+}[{\Cal M}_n] \\
  &+ \sqrt{h^d_{A^+ \restriction \tau,B^+ \restriction \tau}[{\Cal M}_n] 
\times p^u_{A^+,B^+}[{\Cal M}_n]} \times \bold c \times 
(\text{log}\|{\Cal M}_n\|)
\endalign
$$

\item "{$(d)$}"  Assume $B^+_0 \le^+ B^+,A^+ \le^+ A^+_0 \le^+_b B^+$.  
\ub{Then} $(A^+,A^+_0,B^+_0,B^+)$ is a pretender to 
$\text{semi}^*-(k,r)$-good for $({\frak K}^+,\bar h,h_0,\bar p,h_1)$ 
\ub{iff} $(A^+,A^+_0,B^+_0,B^+)$ is a semi$^*-(k,r)$-good for 
$({\frak K}^+,\bar h^+,h^e)$
\sn
\item "{$(e)$}"  inequalities for (d) parallel to those in (c).
\endroster
\enddemo
\bigskip

\demo{Proof}  All is reduced to the case of the binomial distribution
by \S4,\S5. \nl
In more details,
\mn
\ub{Stage A}:  Clauses (b)(i)-(v) and (d), the first iff \nl
The first by \scite{4.10}, the second (i.e. clause (d)) left to the reader.
\mn
\ub{Stage B}:  Clauses (c) + (e) \nl

Let $\varepsilon > 0$ and let ${\Cal M}_n$ be random enough.  We need to
consider all $f \in F = \{f:f \text{ an embedding of } A^+ \restriction
\tau \text{ into } {\Cal M}_n\}$, for each of them the appropriate inequality
should hold.  So it is enough if the probability of failure is 
$< \|{\Cal M}_n\|^{1/\varepsilon}/|F|$ so 
$\le \|{\Cal M}_n\|^{-|A|-1-1/\varepsilon}$ suffice.  
Success means that under the hypothesis 
$f:A^+ \rightarrow {\Cal M}^+_n$, we should consider the
candidates $g \in G = \{g:g \text{ an embedding of } B^+ \restriction \tau$
into ${\Cal M}_n$ extending $f\}$, how many of them will be in $G^+ =
\text{ ex}(f,A^+,B^+,{\Cal M}_n)$.  Well the events $g \in G^+$ are not
independent (in $\mu^+_{{\Cal M}_n}$, of course). \nl
\ub{Note}:  We should prove that the probability of deviating from the
expected number of extensions of $f:A^+ \rightarrow {\Cal M}^+_n$, is much
smaller than the number of $f$'s which is $\le \|{\Cal M}_n\|^{|A|}$.  By
\scite{5.10A} we can restrict ourselves to the separated case (see
Definition \scite{5.10}(3)).  By \scite{5.8} we can ignore the case that in
$G[{\Cal M}^+_n]$ (for our $f,A,B$ or $f,A,A_0,B_0,B,k$ as in \scite{5.1}(2))
every connectivity component has $< \bold c$ element for some $\bold c$
depending on $A,B$ only.

Now our problem is a particular case of the context in \scite{6.1} and by the
previous sentence we can deal separately giving an interval into which the
number of components of isomorphism type $\bold t$, for the $\bold t$ with
$< \bold c$ elements (see \scite{6.3}). Actually each one is another
instance of \scite{6.1} only separated is replaced by weakly separated.  So
clearly it is enough to deal with $L_{\bold t^*}[{\Cal M}^*_n]$ ($\bold t^*$
the isomorphism type of a singleton).  With well known estimates,
\scite{6.7}(2) gives ``very low probability" for the value
$L_{t^*}[{\Cal M}^*_n]$ being too small.  The dual estimate is given by
\scite{6.11A} (remember in our case having instances of the new relations 
has low probability so we can use \scite{6.7}(2) to show that $\zeta$ there
is small even for $\alpha < 1$ very near to 1).

For clause (e) just note that the number of extensions violating the desired
conclusion is much smaller (the definition of pretender is just made for 
this).  
\mn
\ub{Stage C}:  Rest.

Straight by now.  \hfill$\square_{\scite{4.17}}$
\enddemo
\newpage

\head {\S4 Clarifying the probability problem} \endhead  \resetall 
\bigskip

We are still in the context \scite{4.1}.  
Our aim is to clarify the probability problem to which we reduce our aim 
in \scite{4.11} (in the polynomial case). \nl
So our aim is to get good enough upper bounds and lower bounds to the
$h_{A^+,B^+}({\Cal M}_n)$.
\demo{\stag{5.1} Hypothesis}  1) ${\Cal M}_n \in {\Cal K}_n,A^+ \le_{qr} B^+,
\bar A = \langle A_\ell:\ell \le k \rangle,A = A^+ \restriction \tau$, \nl
$B = B^+ \restriction \tau,A = A_0 <_{pr} A_1 <_{pr} \cdots <_{pr} A_k = B,
A^+_\ell = A^+ \restriction A_\ell$ \newline
and an embedding $f^*:A \rightarrow {\Cal M}_n$.

We try to approximate nu$(f^*,A^+,B^+,{\Cal M}^+_n)$ under the condition
``$f^*$ embed $A^+$ into ${\Cal M}^+_n$". \newline
2) (a variant) Assume further $(A^+,A^+_0,B^+_0,B^+)$ is a pretender to
being semi-$(k,r)$-good $B^-_1 = B \restriction B_1,f^*_1:A \rightarrow
{\Cal M}_n$.  We try to approximate nu$^k(f^*,A^+,B^+_1,B^+,{\Cal M}^+_n)$
under the assumption $c \ell^r(f(A^+_0),{\Cal M}^+_n) \subseteq A^+$. 
\mn
\ub{Notation}:  Lastly let $x \in \{ d,u \}$.
\enddemo
\bn
We shall speak mainly for \scite{5.1}(1), and then indicate the changes
for \scite{5.1}(2).
\demo{\stag{5.2} Notation}  $T_\ell = T_\ell[f^*,{\Cal M}_n] = 
\text{ ex}(f^*,A,A_\ell,{\Cal M}_n)$, \newline
$T_\ell[f^*,{\Cal M}^+_n] = T_\ell \cap \text{ ex}(f^*,A^+,A^+_\ell,
{\Cal M}^+_n)$, \newline
$T = \dsize \bigcup_{\ell \le k} T_\ell$ and $T[f^*,{\Cal M}^+_n] = \dsize
\bigcup_{\ell \le k} T_\ell[f^*,{\Cal M}^+_n]$.  
For $g \in T_\ell$, let $\text{lev}
(g) = \ell$. \newline
For $g_1,g_2 \in T$, let $g_1 \doublecap g_2$ be 
$g_1 \restriction A_\ell$ where
$\ell \le k$ is maximal such that \newline
$g_1 \restriction A_\ell = g_2 \restriction A_\ell$.  Let

$$
m^d_\ell =: h^d_{A_\ell,A_{\ell + 1}}[{\Cal M}_n]
$$

$$
m^u_\ell =: h^u_{A_\ell,A_{\ell + 1}}[{\Cal M}_n]
$$

$$
p^d_\ell =: p^d_{A^+_\ell,A^+_{\ell + 1}}[{\Cal M}_n]
$$

$$
p^u_\ell =: p^u_{A^+_\ell,A^+_{\ell + 1}}[{\Cal M}_n].
$$
\medskip

\noindent
Lastly let

$$
\align
{\Cal K}^+_{{\Cal M}_n,f^*} =: \{ {\Cal M}^+_n \in {\Cal K}^+_n:&\,{\Cal M}
^+_n \text{ expand } {\Cal M}_n \text{ and } f^* \text{ is an embedding of} \\
  &\,A^+ \text{ into } {\Cal M}^+_n \}.
\endalign
$$
\medskip

\noindent
Let $\mu^+_n[f^*,{\Cal M}_n]$ be the distribution that $\mu^+_n[{\Cal M}_n]$
induce on ${\Cal K}^+_{{\Cal M}_n,f^*}$.
\enddemo
\bigskip

\demo{\stag{5.3} Hypothesis}  ${\Cal M}_n$ is random enough, so that: 
for $\ell < k,f \in T_\ell$ the number of $g,f \subseteq g \in T_{\ell + 1}$, 
is in the interval $[m^d_\ell,m^u_\ell]$ where $m^x_\ell = h^x_{A_\ell,
A_{\ell + 1}}[{\Cal M}_n]$.
\enddemo
\bigskip

\demo{\stag{5.4} Observation}  For ${\Cal M}^+_n$ random enough and
$f^*:A^+ \rightarrow {\Cal M}^+_n$:
\medskip
\roster
\item "{$(a)$}"  each $f \in T_\ell[{\Cal M}_n]$ with $\ell < k$, has at
least $m^d_\ell$ immediate successors and at most $m^u_\ell$ immediate 
successors
\sn
\item "{$(b)$}"  if $f \in T_\ell[{\Cal M}_n],x \in {\Cal M}_n$
\underbar{then} \newline
$|\{g \in T_{\ell + 1}:f \subseteq g,x \in \text{ Rang }
g \backslash \text{ Rang } f\}| \le \bold c^r_j$ \newline
(where $\bold c^r_j$ is a constant depending on $(A_\ell,A_{\ell + 1})$ only)
\sn
\item "{$(c)$}"  the set of immediate successors of $f$ can be represented
as $\dsize \bigcup_{i < \bold c} T_{\ell,f,i}$, \newline
($\bold c$ depending on
$(A_\ell,A_{\ell +1})$ only and) $\langle \text{Rang } g \backslash
\text{Rang } f:g \in T_{\ell,f,i} \rangle$ are pairwise disjoint.
\endroster
\enddemo
\bigskip

\demo{\stag{5.5} Observation}:  So when ${\Cal M}_n$ is random enough if 
$f \in T_\ell$ embed $A^+_\ell$ into ${\Cal M}^+_n$, \underbar{then} the 
number of $f' \in T_{\ell + 1}$
extending $f$ which embed $A^+_{\ell + 1}$ into ${\Cal M}^+_n$ is
in the expected case in the interval 
$[p^d_\ell \, m^d_\ell,p^u_\ell \, m^u_\ell]$ except for the case $\ell =
k-1$ (then this interval is $\subseteq [0,1)_{\Bbb R}$, so the number
of $f'$ has a bound, but its expected value is $< 1$).
\newline
Of course, for $\ell < k-1$ we expect that for various $f$'s the number will
deviate (from the expected value), 
but for ${\Cal M}^+_n$ random enough none will deviate too much.
\enddemo
\bigskip

\definition{\stag{5.6} Definition}  For 
${\Cal M}^+_n \in {\Cal K}^+_{{\Cal M}_n,f^*}$, we define a graph

$$
G[{\Cal M}^+_n] = G[f^*,{\Cal M}^+_n] = G[f^*,{\Cal M}^+_n,T].
$$

\noindent
Its nodes are $G[{\Cal M}^+_n] = \{f \in T_k:f$ embed $B^+$ into
${\Cal M}^+_n\}$. \newline
Its set of edges is

$$
\align
R[{\Cal M}^+_n] = \biggl\{ 
\{ g_1,g_2\}:&g_1 \in G[{\Cal M}^+_n],g_2 \in
G[{\Cal M}^+_n], \\
  &\text{and Rang}(g_1) \cap \text{ Rang}(g_2) \ne \text{ Rang}
(g_1 \doublecap g_2) \biggr\}.
\endalign
$$
\medskip

\noindent
If $\bold{\frak c} = \{ f_1,\dotsc,f_m\}$ is a component, its domain is
$\dsize \bigcup^m_{\ell =1} \text{ Rang } f_\ell \backslash \text{Rang }
f^*$.
\enddefinition
\bigskip

\proclaim{\stag{5.8} Claim}  Assume that the tuple 
$({\frak K},{\frak K}^+,\bar h,\bar p)$ satisfies $(*)$ of \scite{4.8} and
the irrationality inequality, i.e. $\otimes_1$ of \scite{4.11}. \newline
For any $c \in \Bbb R^+$, for some $m^\oplus(c) \in \Bbb N$,
if ${\Cal M}^+_n$ is random enough, \underbar{then} 
every component of $G[{\Cal M}^+_n]$
has $\le m^\otimes(c)$ members, with the probability of failure 
$\le \| {\Cal M}^+_n \|^{-c}$.
\endproclaim
\bigskip

\demo{Proof}  Choose $\varepsilon > 0$ such that $(\varepsilon < 1$ and)
\mr
\item "{$(*)_1$}"  if $A^+_0 \le C^+ < A^+_k,C^+ \restriction \tau <_s A_k$
then $\beta(C^+,A^+_k) \le - \varepsilon$
\sn
\item "{$(*)_2$}"  $p(A^+_0,A^+_k) \le - \varepsilon$
\ermn
(as $A^+_0 <_{\text{qr}} A^+_k$ each $\beta(C^+,A^+_k) < 0$ by the
irrationality condition and $p(A^+_0,A^+_k) < 0$ as otherwise
$A^+_0 <_t A^+_1 <_t A_2 <_t \ldots <_t \ldots$; so $\varepsilon$ has just to
be below finitely many reals which are $> 0$). 
\sn
Next we choose $m^0$ such that
\mr
\item "{$(*)_3$}"  $\alpha(c \ell(\emptyset,A_0),A_k) - m^0 \times
\varepsilon < -e$
\ermn
(clearly possible).  Next choose $m^1$ such that
\mr
\item "{$(*)_4$}"  if $B \in K$ and $f_\ell:A_k \rightarrow B$ for $\ell <
m^0-1$ are embeddings $f_\ell \restriction A_0 = f_0 \restriction A_0$ and
$C = c \ell^{|A_k|}(\dbcu_\ell f_\ell(A_k),B)$, \ub{then}
$|C|^{|A_k \backslash A_0|} < m^1$ \nl
(actually, somewhat less is needed).
\ermn
So
\mr
\item "{$(*)_5$}"   if $\{f_\ell:\ell < m^1\} \subseteq G[{\Cal M}^+_n]$ is
connected then reordering we have: \nl
for every $\ell \in (0,m^0]$ one of the following occurs:
{\roster
\itemitem{ $(a)$ }  Rang$(f_\ell \restriction (A_k \backslash A_0)) \cap
\dbcu_{m < \ell} \text{ Rang}(f_m) \ne \emptyset$ but \nl
Rang$(f_\ell) \nsubseteq c \ell^{|A_k|}(\dbcu_{m < \ell} \text{ Rang}
(f_m),{\Cal M}_n)$
\sn
\itemitem{ $(b)$ }  Rang$(f_\ell \restriction (A_k \backslash A_0)) \cap
\dbcu_{m < \ell} \text{ Rang}(f_m) = \emptyset$ but \nl
Rang$(f_\ell) \cap c \ell^{|A_k|}(\dbcu_{m < \ell} \text{ Rang}(f_m),
{\Cal M}_n) \ne \emptyset$.
\endroster}
\ermn
[Why?  Suppose this holds for $\ell \in (0,m')$, with $m'$ maximal and assume
that $m' < m^0$ and we shall derive a contradiction; note
that for $m'=0$ this holds trivially.  So if we can find $\ell \in (m',m^1)$
satisfying (a) or (b) we get contradiction to ``$m'$ maximal" so there is no
such $\ell$. \nl
Let $S =: \{\ell < m':\text{Rang}(f_\ell \restriction (A_k \backslash A_0))
\cap \dbcu_{m < m'} \text{ Rang}(f_m) \ne \emptyset\}$. \nl
Note that
\mr
\item "{$(\alpha)$}"  $\ell \le m' \Rightarrow \ell \in S \Rightarrow
\text{ Rang}(f_\ell \restriction (A_k \backslash A_0)) \subseteq 
\dbcu_{m < m'} \text{ Rang}(f_m) \Rightarrow \text{ Rang}(f_\ell \restriction
(A_k \backslash A_0)) \subseteq c \ell^{|A_k|}(\dbcu_{m<m'} \text{ Rang}
(f_m),{\Cal M}_n)$ and
\sn
\item "{$(\beta)$}"  $\ell \in S \and \ell > m' \Rightarrow \ell$ fails
clause (a) $\Rightarrow \text{ Rang}(f_\ell) \subseteq 
c \ell^{|A_k|}(\dbcu_{m \le m'} \text{ Rang}(f_m),{\Cal M}_n)$.
\ermn
By $(*)_4$ we have $|S| < m^1$, so $S \ne \{\ell:\ell < m\}$.  Also
$0 < m'$ so $0 \in S$, hence $S \ne \emptyset \and S \ne \{\ell:\ell < 
m'\}$.  So by the connectivity of $\{f_\ell:\ell < m\}$, i.e. of the graph
$G[{\Cal M}^+_n]$, for some $\ell_1 \in S,\ell_2 < m^1,\ell_2 \notin S$ we
have Rang$(f_{\ell_1} \restriction (A_k \backslash A_0)) \cap \text{ Rang}
(f_{\ell_2} \restriction (A_k \backslash A_0)) \ne \emptyset$.  Now by
$(\alpha) + (\beta)$ above Rang$(f_1)$ is $\subseteq c \ell^{|A_k|}
(\dbcu_{m \le m'} \text{ Rang}(f_m),{\Cal M}_n)$ hence Rang$(f_2 \restriction
(A_k \backslash A_0))$ has an element in this set, but as $\ell_2 \notin S$
it has no element in $\dbcu_{m \le m'} \text{ Rang}(f_m)$, so $\ell_2$
satisfies clause (b) above.  But this contradicts the maximality of $m'$.
So we are done.]
\mn
Now it is enough to fix the isomorphism type of
$({\Cal M}_n \restriction \dbcu_{\ell \le m^0} \text{ Rang}(f_\ell),
f_\ell(d))_{\ell < m^0,d \in A_k}$, call it ${\frak t}$ (as their number is
fixed not depending on $n$).  Let for $g:A_0 \rightarrow {\Cal M}_n$

$$
\align
F_{\frak t}({\Cal M}_n,g) = \biggl\{ \langle f_\ell:\ell \le m^0 
\rangle:&f_\ell \text{ embed } A_k \text{ into } {\Cal M}_n, \\
  &f_\ell \text{ extends } g \text{ and } {\frak t} \text{ is the isomorphism
type of} \\
  &({\Cal M}_n \restriction \dbcu_{\ell \le m^0} \text{ Rang}(f_\ell),
f_\ell(d))_{\ell < m^0,d \in A_k} \biggr\}.
\endalign
$$
\mn
(Note: the ``$\le m^0$" rather than ``$< m^0$" is intentional). \nl
Let $F_{\frak t}[{\Cal M}^+_n,g] = \{ \bar f \in F[{\Cal M}_n,g]:\text{ for }
\ell \le m^0,f_\ell \text{ embeds } A^+_k \text{ into } {\Cal M}^+_n\}$.
We will show that for each ${\frak t}$, for random enough ${\Cal M}_n$,
the expected value of $|F_{\frak t}[{\Cal M}^+_n,g]|$ under the assumption
that $g$ embeds $A^+_0$ into ${\Cal M}_n$ is $\le 
\|{\Cal M}_n\|^{-m^0 \varepsilon} \times \|{\Cal M}_n\|^{\alpha(c \ell
(\emptyset,A_0),A_k)} < \|{\Cal M}_n\|^{-e}$ this clearly suffices. 
\enddemo
\bn
The rest is straight, still we first note \nl
\ub{\stag{5.8A} Observation}:  1) Assume 
$A_0 \le A_1 \le B_1,A_0 \le B_0 \le B_1,A_1 \cup B_0 = B_1$.  \ub{Then} 
$\alpha(c \ell(A_0,A_1),A_1) \ge \alpha (c \ell(B_0,B_1),B_1)$.
\bigskip

\demo{Proof}  1) Let $\varepsilon > 0$ and let ${\Cal M}_n$ be random enough.
We can find $g:B_1 \rightarrow {\Cal M}_n$ and so let

$$
F_1 = \{f:f \text{ embed } A_1 \text{ into } {\Cal M}_n,f \restriction A_0
\subseteq g\}
$$

$$
F_2 = \{f:f \text{ embed } A_1 \text{ into } {\Cal M}_n,f \restriction 
c \ell(A_0,A_1) \subseteq g\}
$$

$$
F_3 = \{f:f \text{ embed } B_1 \text{ into } {\Cal M}_n,f \restriction B_0
\subseteq g\}
$$

$$
F_4 = \{f:f \text{ embed } B_1 \text{ into } {\Cal M}_n,f \restriction 
c \ell(B_0,B_1) \subseteq g\}.
$$
\mn
Clearly
\mr
\item "{$(a)$}" $|F_2| \le \|{\Cal M}_n\|^{\alpha(c \ell(A_0,A_1),A_1) +
\varepsilon}$
\sn
\item "{$(b)$}" $|F_4| \ge \|{\Cal M}_n\|^{\alpha(c \ell(B_0,B_1),B_1) -
\varepsilon}$
\sn
\item "{$(c)$}" $|F_3| \le |F_1|$
\sn
\item "{$(d)$}" $|F_1| \le \bold c|F_2|$ where $\bold c > 0$ is a real
depending on $A_0,A_1$ only
\sn
\item "{$(e)$}" $|F_4| \le |F_3|$.
\ermn
Together (if $\|{\Cal M}_n\|^\varepsilon > \bold c$) we get
$(\alpha(c \ell(A_0,A_1),A_1) + \varepsilon \ge \alpha(c \ell(B_0,B_1),
B_1) - 2 \varepsilon$ but $\varepsilon$ was any positive real so we are 
done.  \hfill$\square_{\scite{5.8A}}$
\enddemo
\bn
\ub{Continuation of the proof of \scite{5.8}}:

$$
B^+_\ell = \dbcu_{m \le \ell} \text{ Rang}(f_m) \text{ and } B_\ell = B_\ell
\restriction \tau, \text{ so } B_0 = A = A_0.
$$
\mn
Let $D_\ell = c \ell(B_\ell,B_{\ell +1})$ and $C^+_\ell = 
\{a \in B^+:f_\ell(a) \in D^+_\ell\}$ and $C_\ell = C^+_\ell \restriction
\tau$.
\mr
\item "{$(*)$}"  $\alpha(D_\ell \restriction B^+_{\ell +1}) \le 
\alpha(C_\ell,B)$ \nl
[why?  apply observation \scite{5.8A} with $A_0,A_1,B_0,B_1$ there 
standing for $(\text{Rang }f_{\ell +1}) \cap D_\ell,(\text{Rang }f_{\ell +1}),
D_\ell,B_{\ell +1}$ (note: $(\text{Rang }f_\ell) \cap D_\ell \subseteq
D_\ell \subseteq B_{\ell +1},(\text{Rang }f_\ell) \cap D_\ell \subseteq
\text{ Rang}(f_\ell) \subseteq B_{\ell +1}$ and $B_{\ell +1} = D_\ell \cup
(\text{Rang }f_{\ell +1})$, \nl
$B_{\ell +1} \in {\Cal K}_\infty$, so the assumption
of observation).  So $\alpha(D_\ell \cap \text{ Rang }f_{\ell +1}$, \nl
$\text{Rang }f_{\ell +1}) \ge \alpha(c \ell(D_\ell,B_{\ell +1}),
B_{\ell +1})$.  But $c \ell(D_\ell,B_{\ell +1}) = D_\ell$ by the
definition of $D_\ell$ and $f_{\ell +1}$ is an embedding mapping $B$ onto
$\text{Rang }f_{\ell +1}$ and $C_\ell$ onto $D_\ell \cap \text{ Rang }
f_{\ell +1}$ so $\alpha(D_\ell \cap \text{ Rang }f_{\ell +1},
\text{Rang }f_{\ell +1}) = \alpha(C_\ell,B)$.  Together we get $(*)$.]
\ermn
Now
$$
\align
\beta(A^+,B^+_{\ell +1}) &= \alpha(A,B_{\ell +1}) + \sum\{\beta(C^+):C^+
\subseteq B^+_{\ell +1},C^+ \nsubseteq A^+\} \\
  &= \alpha(A,D_\ell) + \alpha(D_\ell,B_{\ell +1}) + 
\sum\{\beta(C^+):C^+ \subseteq B^+_{\ell +1} \text{ and } C^+ \nsubseteq
A^+\}\\
  &\le \alpha(A,B_\ell) + \alpha(D_\ell,B_{\ell +1}) + 
\sum\{\beta(C^+):C^+ \subseteq B^+_{\ell +1} \text{ and } C^+ \ne A^+\}\\
  &= \alpha(A,B_\ell) + \sum\{\beta(C):C \subseteq B_\ell,C \nsubseteq A\} \\
  &+ \alpha(D_\ell,B_{\ell +1}) + \sum\{\beta(C^+):C^+ \subseteq
B^+_{\ell +1},C^+ \nsubseteq B^+_\ell\} \\
  &\,\,\,\beta(A^+,B^+_\ell) + \alpha(B^+_\ell,B^+_{\ell +1}) +
\sum\{\beta(C^+):C^+ \subseteq B^+_{\ell +1},C^+ \nsubseteq B^+_\ell\} \\
  &\le \beta(A^+,B^+_\ell) + (\alpha(C_\ell,B) + \sum\{\beta(C^+):C^+ 
\subseteq B^+_{\ell +1},C^+ \nsubseteq B^+_\ell\} \\
  &\qquad \qquad \qquad \qquad \qquad \text{[by } (*) \text{ above]} \\
  &\le \beta(A^+,B^+_\ell) + \beta(C^+_\ell,B^+) + \sum
\{\beta(C^+):C^+ \subseteq B^+_{\ell +1},C^+ \subseteq B^+_\ell \\
  &\text{ \ub{but} } C^+ \subseteq f_\ell(C^+_\ell) \vee C^+ \nsubseteq
f_\ell(B)\}.
\endalign
$$
\bn
\ub{Case 1}:  $C^+_\ell \ne B^+$.

Hence $\beta(A^+,B^+_{\ell +1}) \le \beta(A^+,B^+_\ell) + \beta
(C^+_\ell,B^+) \le \beta(A^+,B^+_\ell) - \varepsilon$ \nl
(why?  first inequality as the third summand above is a sum of reals $\le 0$,
the second inequality by definition of $<^+_{\text{qr}}$ and
choice of $\varepsilon$).
\bn
\ub{Case 2}:  $C^+_\ell = B^+$.

So $\beta(C^+_\ell,B^+) = 0$, so we get \nl
$\beta(A^+,B^+_{\ell +1}) \le \beta(A,B_\ell) + \sum\{\beta(C^+):
C^+ \subseteq B^+_{\ell +1},C^+ \nsubseteq B^+_\ell$ \ub{but} $C^+ \subseteq
f_\ell(C^+_\ell) \vee C^+ \subseteq f_\ell(B^+)\} \le \beta(A^+,B^+_\ell) 
- \varepsilon$.  \hfill$\square_{\scite{5.8}}$
\bn
The following definition points to the fact that there may be quite a
different situation in spite of our treating them together as they are
similar enough for our aims. 

\definition{\stag{5.9} Definition}  1) We say $T$ is simple of the first kind 
\underbar{if}:
\newline
$g_1,g_2 \in T_i \Rightarrow \text{ Rang}(g_1) \cap \text{ Rang}(g_2) =
\text{ Rang}(g_1 \doublecap g_2)$. \newline
2)  We say $T$ is simple of the second kind \underbar{if}: \newline
for every $g_1,g_2 \in T_\ell$, we have

$$
\{ g'_1 \restriction (A_{\ell + 1} \backslash A_\ell):g_1 \subseteq g'_1
\in T_{\ell + 1} \} =
\{ g'_2 \restriction (A_{\ell + 1} \backslash A_\ell):g_2 \subseteq g'_2
\in T_{\ell + 1} \}.
$$

\noindent
3)  $T$ is separated when:

$$
\{g_1,g_2 \} \subseteq T,\{y_1,y_2\} \subseteq A_k,g_1(y_1) = g_2(y_2)
\Rightarrow y_1 = y_2.
$$
\medskip

\noindent
4)  $T$ is locally disjoint if for $f \in T_\ell,f \subseteq g_1 \in
T_{\ell + 1},f \subseteq g_2 \in T_{\ell + 1},g_1 \ne g_2$ we have
Rang$(g_1) \cap \text{ Rang}(g_2) = \text{ Rang}(f)$.
\enddefinition
\bigskip

\remark{\stag{5.10} Remark}  Note the separated 
case assumption helps as it gives monotonicity in the probability. \newline
The following indicates that we can assume $T$ is separated.
\endremark
\bigskip 

\proclaim{\stag{5.10A} Claim}  Each problem (i.e. from \scite{5.1}) we can
split $T_k[f^*,{\Cal M}_n]$ to \nl
$\le (\text{log}_2\|{\Cal M}_n\|)^{|B^+|+1}$
sets, each of them separable.  (So if our estimates can absorb the
inaccuracy involved we have reduced our problem to a separable one).
\endproclaim
\bigskip

\demo{Proof}  Just choose for each $c \in {\Cal M}_n$ a sequence $\rho_c$
of zeroes and ones of length 
$\le \, \rbrack \text{ log}_2\|{\Cal M}_n\| \lbrack$ such that $c \ne d
\Rightarrow \rho_c \ne \rho_d$.  For each $g \in T_k[f^*,{\Cal M}_n]$, let \nl
$w_g = \{\text{Min}\{i:\rho_c(i) \ne \rho_d(i)\}:c \ne d \in
\text{ Rang } g \backslash \text{Rang }f^*\}$, and let  \nl
$v_g = \{\langle c,\rho_{g(c)}(i) \rangle:
c \in B^+ \backslash A^+,i \in w_g\}$.  Define
an equivalence relation $\bold e$ on $T_k[f^*,{\Cal M}_n]$ as follows:
$g_1 \, \bold e \, g_2$ iff $w_{g_1} = w_{g_2} \and v_{g_1} = v_{g_2}$.  Now
$\bold e$ gives a division as required.  \hfill$\square_{\scite{5.10A}}$
\enddemo
\bigskip

\proclaim{\stag{5.11} Discussion}  Now \scite{5.10A} is sufficient for our
aims, but we can get better: division to constant number, and preserving the
order of magnitude of the splitting of the tree.
\endproclaim
\bigskip

\proclaim{\stag{5.11A} Claim}  Consider the probability space of all 
$c$ where $c$ is a function from $A^* = \bigcup \{ \text{Rang}(f):f \in T_k\}
\backslash \text{Rang}(f^*)$ to $A_k \backslash A_0$; all of the 
same probability.  Let the space measure be called $\mu_{sep}$.  Let 
$T^{[c]} = \{ \eta:c \circ f_\eta$ is the identity on $A_{\ell g(\eta)}\}$ 
and $T^{[c]}_\ell = T^{[c]} \cap T_\ell$.
\newline
Then 
\medskip
\roster
\item   The probability that all the splitting in $T^{[c]}$ are nearly the
expected value (meaning if the expected value is $v$, error is 
$\le v^{{\frac 12}+ \varepsilon}$)
is very near to $1$ assuming the \underbar{largeness condition}, (e.g. in
the polynomial case)
\sn
\item   Let $e^* = ((|A_k \backslash A_0|)!)
|A^*_{k^*} \backslash A^*_0|^{-|A^*_{k^*} \backslash
A^*_0|} \in (0,1)_{\Bbb R}$ and $a_{dn}$ (or $a_{up}$) be the natural lower 
(or upper bound) of the expected value of $|T[{\Cal M}^+_n]|$,
\underbar{then} for any $\varepsilon > 0$ we have $\mu_{sep}$ - almost surely
the chosen $c$ satisfies \newline
$\text{Prob}_{\mu_n} \left( |T_k \cap T[{\Cal M}^+_n]|) \le a_{dn} -
a^{{\frac 12} + \varepsilon}_{dn} \right) \le$ \newline  
$\text{Prob}_{\mu_n} \left( |T^{[c]}_k \cap T[{\Cal M}^+_n]| \le e^*a_{dn} -
a^{{\frac 12} + \varepsilon} \right) + \| {\Cal M}_n\|^{- \varepsilon}$.
\sn
\item  $\text{Prob}_\mu \left( |T_k \cap T[{\Cal M}^+_n]| \ge a_{up} +
a^{{\frac 12} + \varepsilon}) \right) \le$ \newline  
$\text{Prob} \left( T_k^{[c]} \cap T[{\Cal M}^+_n] \le e^*a_{up} +
a^{{\frac 12} + \varepsilon}_{up} \right) + \| {\Cal M}_n \|^{- \varepsilon}$.
\endroster
\endproclaim
\bigskip

\demo{Proof}  If we first draw ${\Cal M}^+_n$ then ignoring an event with 
probability $\| {\Cal M}_n \|^{-e}$, the components of
$T[{\Cal M}^+_n]$ are all of size $\le m^\otimes(e)(\in \Bbb N$, from
\scite{5.8}).  So
the number we get after drawing a $\mu_{sep}$-random $c$, behave by a
multinomial distribution, so almost surely for $c$, we get the expected
number with small error.  By commutativity of probability, this implies the
conclusion. \hfill$\square_{\scite{5.11A}}$
\enddemo
\newpage

\head {\S5 The probability arguments} \endhead  \resetall
\bigskip

We relax our framework, forget about the tree (from \S4), and 
just have a family $F$
of one-to-one functions from $[m]$ to $[n]$ (thinking $n,|F|$ are much larger
than $m$), $F$ seperative for simplicity (i.e. $f_1(\ell) = f_2(\ell_2)
\Rightarrow \ell_1 = \ell_2$), and $A^*$ a $\tau^*$-model with set of
elements $[m]$ with vocabulary $\tau^+$.  Now we draw some relations 
on $[n]$ to get ${\Cal M}^*_n$ independently and want to know enough on the 
number ${\Cal L}$ of $f \in F$ such that all
appropriate relations were chosen.  The easiest case is when $f_1 \ne f_2 \in
F \Rightarrow \text{ Rang}(f_1) \cap \text{ Rang}(f_2) = \emptyset$, then we
get a binomial distribution.  Still we are interested in the case that for 
every ``successful" $f \in F$, the number of succsesful $f' \in F \backslash
\{ f\}$ not disjoint to $f$ is small; i.e. the expected number is $\ll 1$.
So we define components of the set of successful $f$, and look what is their
number.  We first show that for $L$ which is larger than the expected value
the probability of having ``the number sucesses is $L$"
to decrease with
$L$ as in the binomial distributions.  Then we get a slightly worse statement
on what occurs for $L$ smaller than the expected value; the ``error" term"
comes from the number of $f \in F$ such that in $f([m])$ there is no relation
but if we change ${\Cal M}^*_n$ such that $f$ is succsesful, it is not a 
singleton.
But by the first argument (or direct checking in the cases from \S4) this is
small.  Of course, we have larger components, but for each isomorphic type
the problem of the distribution of their number is like the original one 
except being only weakly separative.  Clearly this framework is wide 
enough to include what is needed in \S4. \nl
\ub{Note}:  Clearly the higher components contribute little but we do not
elaborate as there is no need: we may restrict ourselves to finitely many
components.
\bigskip

\definition{\stag{6.1} Definition}  1) We say $\bar y = (m,\bar p^n,n,F^n) 
= (m,\bar p,n,F) = (m^{\bar y},p^{\bar y},n^{\bar y},F^{\bar y})$ is a 
system (or $m$-system or $(m,n)$-system) if:
\medskip
\roster
\item "{$(a)$}"  $F$ is a family of one-to-one functions with domain 
$[m] = \{1,\dotsc,m\}$ into the set $[n] = \{1,\dotsc,n\}$
\sn
\item "{$(b)$}"  ${\Cal P} \subseteq \{ \bold u:\bold u \subseteq [m],\bold u
\ne \emptyset\},\bold e$ an equivalent relation on ${\Cal P}$,
\sn
\item "{$(c)$}"  $\bar p^n = \langle p_{\bold u / \bold e}:\bold u / \bold e 
\in {\Cal P} / \bold e \rangle$, where $p_{\bold u / \bold e}$ is
a probability and let $p_\bold u = p_{\bold u / \bold e}$.
Let
$$
R_\bold u = \{ f(\bold u):f \in F\},R_{\bold u / \bold e} = \bigcup
\{ R_{\bold u'}:\bold u' \in \bold u / \bold e\}
$$
\noindent
(we can look at them as symmetric relations)
\sn
\item "{$(d)$}"  we choose for each $\bold u \in {\Cal P},
{\Cal R}_{\bold u/\bold e} \subseteq R_{\bold u/\bold e}$ by 
drawing for each $v \in R_{\bold u/\bold e}$ a decision 
for $v \in {\Cal R}_{\bold u/\bold e}$, independently (for distinct
$(\bold u,v)$) with probability $p_{\bold u/\bold e}$.  The distribution (on
${\Cal K}^*_n$ (see below)) is called $\mu^*_n$.
\endroster
\medskip

\noindent
2) We call ${\Cal M}^*_n = ([n],\dotsc,R_{\bold u/\bold e},
{\Cal R}_{\bold u/\bold e},\ldots)_{\bold u \in {\Cal P}}$ a 
$\mu_n[\bar y]$-random model; we may omit $\bar y$.  
Let ${\Cal K}^*_n = {\Cal K}^*_n[\bar y] = 
K_{\bar y}$ be the set of all possible ${\Cal M}^*_n$.  Note that ${\Cal P},
\bold e$ can be defined from $\bar p^n$, so we write ${\Cal P} = {\Cal P}
^{\bar y},\bold e = \bold e^{\bar y}$ or ${\Cal P}^n,\bar e^n$.  
Note that from $M_n \in {\Cal K}^*_{\bar y}$ we can reconstruct 
$m^{\bar y},n^{\bar y}$ and $F^{\bar y}$
(though not $\bar p^{\bar y}$) and from $F^{\bar y}$ we can reconstruct
$m^{\bar y},n^{\bar y}$ if $[m] = \bigcup\{\bold u:\bold u \in {\Cal P}\}$.
\enddefinition 
\bigskip

\centerline {$* \qquad * \qquad *$}
\bigskip

It is natural to demand $F$ is separative (see Definition \scite{6.2}(1) 
below), as we can reduce the general case to this one (though increasing 
the ``error" term see \scite{5.10A}, \scite{5.11A}).  But why do we consider
``weakly separative" (see Definition \scite{6.2}(2) below)?  The main 
arguments here gives reasonable estimates if we have estimated the number
of occurances of non-trivial components, so we need to estimate them.  In
order to bound the ``error" gotten by making the estimation of the number
of occurances of a component, we weaken the Definition to include this
(we could somewhat further weaken ``weakly separative", but do not, just as
it gives no reasonable gain now).
\bigskip

\definition{\stag{6.2} Definition}  1) We say $F = F^n$ (or 
$\bar y = (m,\bar p^n,n,F^n)$) is separative \underbar{if}: \newline
$\{ f_1,f_2\} \subseteq F_n \and \{ \ell_1,\ell_2\} \in [m] \and f_1(\ell_1)
= f_2(\ell_2) \Rightarrow \ell_1 = \ell_2$. \newline
2)  We call $\bar y$ semi-separative if: \newline
there is an equivalence relation $\bold e^*$ on $[m]$ such that
\mr
\widestnumber\item{$(iii)$}
\item "{$(i)$}"  $\bold u \in {\Cal P} \Rightarrow (\forall \ell,k)
(\ell \in \bold u \and k \in \bold u \and \ell \ne k 
\rightarrow \neg \ell \bold e^* k)$ and
\sn
\item "{$(ii)$}"  $(\forall \bold u_1,\bold u_2 \in {\Cal P})
((\forall k)[\bold u_1 \cap (k/\bold e^*) \ne \emptyset \leftrightarrow 
\bold u_2(k \cap \bold e^*)] \rightarrow \bold u_1 \, \bold e \, 
\bold u_2)$; \nl
i.e. $\bold e$ refines
$\bold e^{**} =: \{(\bold u_1,\bold u_2):\{\bold u_1,\bold u_2\} \subseteq 
{\Cal P}$ and \nl
$(\forall k)[\bold u_1 \cap (k/\bold e^*) \ne \emptyset 
\leftrightarrow \bold u_2 \cap (k/\bold e^*) \ne 0]\}$
\sn
\item "{$(iii)$}"   there is an equivalence relation $\bold e'$ on $[n]$
such that: \nl
if $\{f_1,f_2\} \subseteq F,\{m_1,m_2\} \subseteq [m]$, \ub{then} \nl
$m_1 \, \bold e^* \, m_2 \Leftrightarrow f_1(m_1)\bold e' \, f_2(m_2)$
\sn
\item "{$(iv)$}"  if 
$\{ f_1,f_2\} \subseteq F^n,\{ \bold u_1,\bold u_2\} \subseteq {\Cal P}$ and
$f_1(\bold u_1) = f_2(\bold u_2)$ \underbar{then} \newline
$\bold u_1 \, \bold e \, \bold u_2 \and
f_1 \restriction \bold u_1 = f_2 \restriction \bold u_2$.
\ermn
3)  We say $\bar y$ is weakly separative if $(i), (ii), (iii)$ of 
part (2) holds. \nl
4)  For $X \subseteq {\Cal P}$ let

$$
q_X = \dsize \prod_{\bold u \in X} p_{\bold u} \times 
\dsize \prod_{\bold u \in {\Cal P} \backslash X} (1-p_\bold u).
$$
\enddefinition
\bigskip

\remark{Remark}  Note:  we are thinking of the 
cases the $p_{\bold u}$'s are small.  If some are essentially 
constant, we treat them separately.
\endremark
\bigskip

\definition{\stag{6.3} Definition}  1) For a system 
$\bar y = (m,\bar p,n,F)$ and $M^*_n \in {\Cal K}^*_{\bar y}$, we 
define a graph $G[M^*_n] = G_{\bar y}[M^*_n] = (F[M^*_n],E[M^*_n])$. \newline
Its set of nodes is

$$
F[M^*_n] = \{f \in F:f(\bold u) \in R_{\bold u/\bold e}^{M^*_n} 
\text{ for every } \bold u \in {\Cal P}\}
$$

$$
E[M^*_n] = E[F] \restriction F[M_n]
$$

$$
\text{where } E[F] = \{(f_1,f_2):f_1 \in F,f_2 \in F,f_1 \ne f_2 
\text{ and Rang}(f_1) \cap \text{ Rang}(f_2) \ne \emptyset\}.
$$
\medskip

\noindent
For a component ${\Cal C}$ of $G[M^*_n]$ let $V[{\Cal C}] = 
\dsize \bigcup_{f \in {\Cal C}} \text{ Rang}(f)$.
\smallskip
\noindent
2) We say that two components ${\Cal C}_1,{\Cal C}_2$ of $G[M^*_n]$ 
are isomorphic if there is a \newline
$\bold f:{\Cal C}_1 \overset{1-1}\to{\underset \text{onto}\to \longrightarrow} 
{\Cal C}_2$ such that for any $a_1,a_2 \in [m]$ and $f_1,f_2 \in C_1$ we have:

$$
(\bold f(f_1))(a_1) = (\bold f(f_2))(a_2) \Leftrightarrow f_1(a_1) =
f_2(a_2)
$$ 
\medskip

\noindent
2A) We say that $\bold f$ is an embedding of a possible component ${\Cal C}_1$
into a possible component ${\Cal C}_2$ (possible means it is a component in 
some $M^*_n$) \newline
if: $\bold f$ is a one to one function from ${\Cal C}_1$ into ${\Cal C}_2$ 
satisfying the demand in (2).  The isomorphism type ${\Cal C}/\cong$ 
of a component ${\Cal C}$, is naturally defined. \newline
3) Let $\bold T_k = \bold T^m_k = 
\{ {\Cal C}/\cong:\text{for some } \bar y \text{ an } m
\text{-system and some } M^*_n \in {\Cal K}^*_{\bar y},{\Cal C} 
\subseteq F^{\bar y}$ is \newline

$\qquad \qquad \qquad \qquad \qquad \quad$ connected component in the graph 
$G[M^*_n] \text{ and}$ \newline

$\qquad \qquad \qquad \qquad \qquad \quad |{\Cal C}| = k\}$. \newline
4)  Assume $\bar y$ is semi-separative, $\bold f_\ell:[m] \rightarrow 
[m^*]$ for $\{\ell=1,\dotsc,k\}$. \nl
Let $\bold T = \bold T^m = \dsize \bigcup_k \bold T^m_k$.  
Let $\bold t^*$ be the isomorphic type of singleton.  Normally $m$ is constant
so we may omit it. \newline
5) Let $L_{\bold t}[M^*_n]$ be $|{\Cal L}_{\bold t}[M^*_n]|$ where
\newline
${\Cal L}_{\bold t}[M^*_n] = \{ {\Cal C}:{\Cal C} \text{ a component of }
G[M^*_n]$ such that the isomorphic type of ${\Cal C}$ is $\bold t\}$. \nl
6) $\bar L = \bar L[M^*_n] = \langle L_{\bold t}[M^*_n]:\bold t \in 
\dsize \bigcup_k \bold T_k \rangle$. \newline
7) ${\Cal K}^*_n[F,\bar L] = \{ M_n \in {\Cal K}^*_n[F]:\bar L[M^*_n] 
= \bar L\}$. \newline
8) $F_X[M_n] = \{ f \in F:\bold u \in {\Cal P} \Rightarrow [f(\bold u)
\in {\Cal R}_\bold u \equiv \bold u \in X]\}$ for $X \subseteq {\Cal P}$. 
\newline
9) $F_*[M^*_n] = F^*_\emptyset[M^*_n]$ where for $X \subseteq
{\Cal P}$ we let: 

$$
\align
F^*_X[M^*_n] = \biggl\{ f \in F:&f \in F_X[M^*_n] 
\text{ and for no } f' \in F \backslash \{f\} \text{ do we have:} \\
  &\bold u \in {\Cal P} \and 
f'(\bold u) \nsubseteq f([m]) \Rightarrow f'(\bold u) \in 
{\Cal R}_{\bold u/\bold e}[M^*_n] \biggr\}.
\endalign
$$
\medskip

\noindent
10) $F_\otimes[M^*_n] = F^\otimes_\emptyset[M^*_n]$ where for 
$X \subset {\Cal P}$ we let $F^\otimes_X[M^*_n] =: 
\{f \in F:f \in {\Cal L}_{t^*}[M^*_n]\}$.
\enddefinition
\bigskip

\proclaim{\stag{6.4} Claim}  1) Separative implies semi-separative which
implies weakly separative. \newline
2) If $\bar y$ is weakly separative, $\bold f_\ell:[m] \rightarrow [m^*]$
for $\ell = 1,\dotsc,k,\{\bold f_\ell:\ell \in [k]\}$ is weakly separative,
$\dsize \bigcup_{\ell \in [k]} \text{ Rang}(\bold f_\ell) = [m^*]$ and
\newline
$F^* = \{g:g:[m^*] \rightarrow [n] \text{ is one to one, } \ell \in [k] 
\Rightarrow g \circ \bold f_\ell \in F^{\bar y}\},
p^*_{\bold f_\ell(u)} = p_u$,\newline
${\Cal P}^* = \{ \bold f_\ell(u):u \in {\Cal P}\}$ \underbar{then} 
$\bar y^* = (m^*,\bar p^*,n,F^*)$ is weakly separative. \nl
3) Similarly with semi-separative instead weakly separative.
\endproclaim
\bigskip

\demo{Proof}  Straight.
\enddemo
\bigskip

\remark{Remark}  The following claim says that above 
the expected value, the probability goes down fast enough.
\endremark
\bigskip

\proclaim{\stag{6.5} Claim}  Assume $F$ is weakly separative and
$\bar L^\ell = \langle L^\ell_{\bold t}:\bold t \in T \rangle$, 
for $\ell = 1,2$ and 
\medskip

\noindent
$$
\bar L^2_{\bold t} = \cases L^1_{\bold t} + 1 &\text{\underbar{if} }
\bold t = \bold t^* \text{ (is singleton)}\\
  L^1_{\bold t} &\text{ otherwise}. \endcases
$$
\medskip

\noindent
\underbar{Then} ($q_\emptyset,q_{\Cal P}$ are defined in \scite{6.2}(4)):

$$
\text{Prob}_{\mu^*_n}(\bar L[{\Cal M}^*_n] = \bar L^1) \ge q_\emptyset
\left( \frac{L^2_{\bold t^*}}{q_{\Cal P}|F|} \right) \text{ Prob}_{\mu^*_n}
\biggl( \bar L[{\Cal M}^*_n] = \bar L^2 \biggr).
$$
\endproclaim
\bigskip

\remark{\stag{6.5A} Remark}  1) On $F_{\Cal P}$ see Definition \scite{6.3}(7).
\nl
2)  Note: Exp$(F_{\Cal P}[{\Cal M}^*_n]) = q_{\Cal P}|F|$ and
$|F_{\Cal P}[{\Cal M}^*_n]| = \dsize \sum_{\bold t \in \bold T} |\bold t| 
L_{\bold t}[{\Cal M}^*_n]$ where $|\bold t|$ is the number of $f \in C$ 
for any $C$ of isomorphism type $\bold t$.
\endremark
\bigskip

\demo{Proof}  Let us consider for $X \subset {\Cal P}$,

$$
\align
W_X =: \biggl\{ (M^1,M^2):&M^1 \in {\Cal K}^*_n[F,\bar L^1],M^2 \in 
{\Cal K}^*_n[F,\bar L^2],\\
  &\text{for } \bold t \in T \backslash \{ \bold t^*\},{\Cal L}_{\bold t}
[{\Cal M}^1] = {\Cal L}_{\bold t}[{\Cal M}^2] \text{ and} \\
  &{\Cal L}_{\bold t^*}[{\Cal M}^1] \subseteq {\Cal L}_{\bold t^*}[{\Cal M}^2]
 \text{ and if } \{f\} = {\Cal L}_{\bold t^*}[{\Cal M}^2] \backslash 
{\Cal L}_{\bold t^*} [{\Cal M}^1] \\
  &\text{(note that necessarily there is one and only one such } f) \\
  &\text{\ub{then} } f(\bold u) \in {\Cal R}^{M_1}_{\bold u/\bold e} 
\Leftrightarrow \bold u \in X \text{ for each } \bold u \in {\Cal P} \\
  &\text{and: if } f' \in F,\bold u \in {\Cal P}^+,f'(\bold u) \nsubseteq
\text{ Rang}(f) \\
  &\text{then } \bold u \in {\Cal R}^{M_1}_{\bold u/\bold e} 
\Leftrightarrow \bold u \in {\Cal R}^{M_2}_{\bold u/\bold e} \biggr\}
\endalign
$$

$\quad W =: \dsize \bigcup_{X \subset {\Cal P}} W_X$.
\medskip

\noindent
Clearly
\medskip
\roster
\item "{$(*)_1$}"  for every $X \subset {\Cal P}$ and $M^2 \in {\Cal K}^*_n
[F,\bar L^2]$ the set \newline
$\{M^1:(M^1,M^2) \in W_X\}$ has exactly $L^2_{\bold t^*}$ members
\endroster
\medskip

\noindent
(here we use ``$F$ is weakly separative") and
\medskip
\roster
\item "{$(*)_2$}"  for every $M^1 \in {\Cal K}^*_n[F,\bar L^1]$, the number of
$\{M^2:(M^1,M^2) \in W_X\}$ has at most 
$|F| - \dsize \sum_{\bold t \in \bold T}|\bold t| 
\cdot L^1_{\bold t} \le |F|$ members \newline
[note: this is a quite crude bound as it does not take into account that
not only each $f \in C \in {\Cal L}^1_{\bold t}$ is disqualified as the
possible member of ${\Cal L}^2_{\bold t} \backslash {\Cal L}^1_{\bold t}$, but
also each $f' \in F$ such that $(\exists \bold u \in {\Cal P})(f(\bold u) =
f'(\bold u))$; but at present the effect does not disturb us].
\medskip
\noindent
\item "{$(*)_3$}"  for every $(M^1,M^2) \in W_X$ and 
$X \subseteq {\Cal P}$, we have (see Definition \scite{6.2}(4)) \newline
Prob$_{\mu_n}({\Cal M}^*_n = M^1)/q_X = \text{ Prob}_{\mu_n}({\Cal M}^*_n 
= M^2)/q_{\Cal P}$.
\endroster
\mn
Now
\bigskip

\noindent
$$
\multline
|F| \text{ Prob}_{\mu_n}(\bar L[{\Cal M}^*_n] = \bar L^1)/q_\emptyset \\
  \text{[by the definition of } {\Cal K}_n[F,\bar L^1]]
\endmultline
$$
\smallskip

\noindent
$$
\multline
= \dsize \sum_{M^1 \in {\Cal K}_n[F,\bar L^1]} \qquad \quad
|F| \text{ Prob}_{\mu_n}({\Cal M}^*_n = M^1)/q_\emptyset \\
  \text{[by } (*)_2]
\endmultline
$$
\smallskip

\noindent
$$
\multline
\ge \dsize \sum_{M^1 \in {\Cal K}_n[F,\bar L^1]} \qquad \quad 
\dsize \sum_{M^2 \text{ satisfies }(M^1,M^2) \in W_\emptyset} \qquad \quad
\text{ Prob}_{\mu_n}({\Cal M}^*_n = M^1)/q_\emptyset \\
  \text{[by interchanging sums]}
\endmultline
$$
\smallskip

\noindent
$$
\multline
= \dsize \sum_{M^2 \in {\Cal K}_n[F,\bar L^2]} \qquad \quad 
\dsize \sum_{M^1 \text{ satisfies }(M^1,M^2) \in W_\emptyset} \qquad \quad
\text{ Prob}_{\mu_n}({\Cal M}^*_n = M^1)/q_\emptyset \\
  \text{[by } (*)_3]
\endmultline
$$
\smallskip

$$
\multline
= \dsize \sum_{M^2 \in {\Cal K}_n[F,\bar L^2]} \qquad \quad
\dsize \sum_{M^1 \text{ satisfies }(M^1,M^2) \in W_\emptyset} \qquad \quad
\text{ Prob}_{\mu_n}({\Cal M}^*_n = M^2)/q_{\Cal P} \\
  \text{[by } (*)_1]
\endmultline
$$
\smallskip

$$
\multline
= \dsize \sum_{M^2 \in {\Cal K}_n[F,\bar L^2]} \qquad \quad L^2_{\bold t^*} 
\times \text{ Prob}_{\mu_n}({\Cal M}^*_n = M^2)/q_{\Cal P} \\
  \text{[by the definition of } {\Cal K}_n[F,\bar L^2]]
\endmultline
$$
\smallskip

$\,\,= L^2_{\bold t^*} \times \text{ Prob}(\bar L_{\bold t^*}[{\Cal M}^*_n] = 
\bar L^2)/q_{\Cal P}$.
\bigskip

\noindent
Now the conclusion follows. \hfill$\square_{\scite{6.5}}$
\enddemo
\bigskip

\proclaim{\stag{6.6} Claim}  Assume $F$ is weakly separative.  
If $L^1 + 1 = L^2$, \ub{then}

$$
\text{Prob}_{\mu^*_n} \biggl( L_{\bold t^*}[{\Cal M}^*_n] = L^1 \biggr) 
\ge q_\emptyset \biggl( {\frac{L^2}{q_{\Cal P} \times |F|}}
\biggr) \times \text{ Prob}_{\mu^*_n} \biggl( L_{\bold t^*}
[{\Cal M}^*_n] = L^2 \biggr).
$$
\endproclaim
\bigskip

\demo{Proof}  By \scite{6.5}, dividing the event to cases.
\hfill$\square_{\scite{6.6}}$
\enddemo
\bigskip

\demo{\stag{6.7} Conclusion}  Assume $F$ is weakly separative. \newline
1) If $L^* \ge 2(q_{\Cal P}|F|/q_\emptyset)$ \underbar{then}
Prob$_{\mu_n}(L_{\bold t^*}[{\Cal M}_n] > L^*) \le \text{ Prob}_{\mu_n}
(L_{t^*}[{\Cal M}_n] = L^*)$. \newline
2) If $L^* \ge [q_{\Cal P}|F|/q_\emptyset]$ \underbar{then} \newline

$\qquad \qquad$ 
Prob$_{\mu_n}(L_{\bold t^*}[{\Cal M}^*_n] = L^*) \le \dsize \prod^{L^*-1}
_{L = [q_{\Cal P}|F|/q_\emptyset]} \quad \frac{L}{[q_{\Cal P}
|F|/q_\emptyset]}$.
\enddemo
\bigskip

\demo{Proof}  Iterate \scite{6.5}, i.e. by induction on $L^*$. 
\hfill$\square_{\scite{6.7}}$
\enddemo
\bigskip

\remark{Remark}  1) In the cases we are interested in, $q_\emptyset$ is
near to one or at least above a constant $> 0$ and $q_{\Cal P}$ is small, so
$(q_{\Cal P}|F|)/q_\emptyset$ is near the expected value of
$|F[{\Cal M}_n]|$. \newline
2) This is enough for the bounds in the case ``above the expected value" for
\S4 modulo separability as we do it for each component. 
\endremark
\bigskip

\proclaim{\stag{6.8} Claim}  Let $F$ be weakly separative.  Assume 
$\bar L^\ell$ for $\ell =1,2$ are as in \scite{6.5} and (see 
Definition \scite{6.3}(8))
\medskip
\roster
\item "{$\bigotimes$}"  $\alpha \in (0,1)_{\Bbb R}$ and $\zeta \in 
(0,1)_{\Bbb R}$ and \newline
$\zeta \ge \text{ Prob}_{\mu^*_n} \biggl( \bar L[{\Cal M}^*_n] = \bar L^1$ 
and $\bigl| F_*[{\Cal M}^*_n] \bigr| < \alpha|F| \biggr)$
\endroster
\medskip

\noindent
(on $F_*[{\Cal M}^*_n]$ see Definition \scite{6.3}(8)). \underbar{Then}

$$
\text{Prob}_{\mu_n}(\bar L[{\Cal M}^*_n] = \bar L^1) \le q_\emptyset
\frac{L^2_{\bold t^*}}{q_{\Cal P}\alpha|F|} \, \text{ Prob}_{\mu^*_n}
(\bar L[{\Cal M}^*_n] = \bar L^2) + \zeta.
$$
\endproclaim
\bigskip

\remark{Remark}  In the cases we have in mind, for $n$ going to infinity,
$\alpha$ goes to 1, $\zeta$ goes to 0, very fast indeed.
\endremark
\bigskip

\demo{Proof}  Start as in the proof of \scite{6.5} getting $(*)_1,
(*)_2,(*)_3$ but then:
\newpage

$$
\multline
\alpha|F| \text{ Prob}_{\mu_n}(\bar L[{\Cal M}^*_n] = \bar L^1)/q_\emptyset \\
  \text{[by } \bigotimes 
\text{ in the assumption and by the Definition of } {\Cal K}_n[F,\bar L^1]]
\endmultline
$$
\smallskip

$$
\multline
\le \alpha \times |F| \times \zeta + 
\dsize \sum_{M^1 \in {\Cal K}_n[F,\bar L^1],|F_*[{\Cal M}^*_n]| 
\ge \alpha|F|} \qquad \quad
\alpha|F| \text{ Prob}_{\mu_n}({\Cal M}^*_n = M^1)/q_\emptyset \\
  \text{[by } (*)_2]
\endmultline
$$
\smallskip

$$
\multline
\le \alpha \times |F| \times \zeta + 
\dsize \sum_{M^1 \in {\Cal K}_n[F,\bar L^1]} \qquad \quad 
\dsize \sum_{M^2 \text{ satisfies }(M^1,M^2) \in W_\emptyset} \qquad \quad
\text{ Prob}_{\mu_n}({\Cal M}^*_n = M^1)/q_\emptyset \\
  \text{[by interchanging sums]}
\endmultline
$$
\smallskip

$$
\multline
= \alpha \times |F| \times \zeta + 
\dsize \sum_{M^2 \in {\Cal K}_n[F,\bar L^2]} \qquad \quad 
\dsize \sum_{M^1 \text{ satisfies }(M^1,M^2) \in W_\emptyset} \qquad \quad
\text{ Prob}_{\mu_n}({\Cal M}^*_n = M^1)/q_\emptyset \\
 \text{[by } (*)_3]
\endmultline
$$
\smallskip

$$
\multline
= \alpha \times |F| \times \zeta + 
\dsize \sum_{M^2 \in {\Cal K}_n[F,\bar L^2]} \qquad \quad
\dsize \sum_{M^1 \text{ satisfies }(M^1,M^2) \in W_\emptyset} \qquad \quad
\text{ Prob}_{\mu_n}({\Cal M}^*_n = M^2)/q_{\Cal P} \\
  \text{[by } (*)_1]
\endmultline
$$
\smallskip

$$
\multline
= \alpha \times |F| \times \zeta + 
\dsize \sum_{M^2 \in {\Cal K}_n[F,\bar L^2]} \qquad \quad 
L^2_{\bold t^*} \times \text{ Prob}_{\mu_n}({\Cal M}^*_n = M^2)/q_{\Cal P} \\
  \text{[by the definition of } {\Cal K}_n[F,\bar L^2]]
\endmultline
$$
\smallskip

$\,\,= \alpha \times |F| \times \zeta + 
L^2_{\bold t^*} \times \text{ Prob}(\bar L_{\bold t^*}
({\Cal M}_n) = \bar L^2)/q_{\Cal P}$. 
\medskip

\noindent
Dividing by $\alpha|F|/q_\emptyset$ we get the desired conclusion.
\hfill$\square_{\scite{6.8}}$
\enddemo
\bigskip

\centerline {$* \qquad * \qquad *$}
\bigskip
\noindent

\proclaim{\stag{6.11A} Claim}  Assume $F$ is weakly separative.
Assume $L^2 = L^1 + 1$ are given and
\medskip
\roster
\item "{$(*)$}"  $\alpha \in (0,1)_{\Bbb R}$ and $\zeta \in (0,1)_{\Bbb R}$
and
\smallskip
$$
\zeta \ge \text{ Prob}_{\mu^*_n} \biggl( L^1_{\bold t^*}[{\Cal M}^*_n] =
L^1 \text{ and } \biggl| F_*[{\Cal M}^*_n] \biggr| < \alpha |F| \biggr).
$$
\endroster
\medskip

\noindent
\underbar{Then}

$$
\text{Prob}_{\mu^*_n} \biggl( L_t[{\Cal M}^*_n] = L^1 \biggr) \le
{\frac{\zeta}{\bold t^*}} + q_\emptyset {\frac{L^2_{\bold t^*}}
{q_{\Cal P} \alpha|F|}} \text{ Prob}_{\mu^*_n} \biggl( L_{\bold t^*}
[{\Cal M}^*_n] = L^2 \biggr).
$$
\endproclaim
\bigskip

\demo{Proof}  By \scite{6.8}, dividing the even to case, noting that:

$$
\align
\text{Prob}_{\mu^*_n}\bigl( L_{t^*}[{\Cal M}^*_n] &= L_1 \text{ and }
\bigl| F_*[{\Cal M}^*_n] \bigr| < \alpha|F| \bigr) \\
  &= \Sigma \bigl\{ \text{Prob}_{\mu^*_n} \bigl( \bar L[{\Cal M}^*_n] = \bar L
\text{ and } \bigl| F_*[{\Cal M}^*_n] \bigr| < \alpha|F| \bigr):\bar L \\
  &\qquad \qquad \qquad \qquad = \langle 
L_{\bold t}:\bold t \rangle,L_{\bold t^*} = L_1 \bigr\}.
\endalign
$$
\sn
${}$ \hfill$\square_{\scite{6.11A}}$
\enddemo
\newpage

\head {\S6 Free Amalgamation} \endhead  \resetall
\bigskip

We like to axiomatize ``free amalgamation" in its connection to 0-1 laws
(in previous cases the ``edgeless disjoint amalgamation" serves).
\bigskip

\centerline {$* \qquad * \qquad *$}
\bigskip

We first define a context having ``free amalgamation".  The idea is that it
is not necessarily a ``disjoint amalgamation with no additional relations" as
we may allow say a two-place relation with probability $\frac 12$, so cases
of this relation has no influence on the amalgamation being free.
\bigskip

\definition{\stag{2.6A} Definition}  1) We say $({\frak K},\nonfork{}{}_{})$ 
is a 0-1 context [or weakly 0-1 context] 
(no contradiction to \scite{2.1}(1)) \underbar{if} ${\Cal K}$ satisfies 
(a),(b),(c) as in \scite{2.1}(1) and $({\Cal K}_\infty,<_s$ are defined as
in \scite{2.3}(1),\scite{2.3}(2)(c) above and): \newline
(d)  $\nonfork{}{}_{}$ is a four-place relation on ${\Cal K}_\infty$ 
written as $\nonforkin{B}{C}_{A}^{D}$ or 
$\nonfork{}{}_{} \, (A,B,C,D)$.  This relation is preserved under 
isomorphism and we say: 
$B,C$ are $\nonfork{}{}_{}$-freely amalgamated over $A$ inside $D$ 
(and omit $D$ if clear from the context)  \newline
(e)  $\nonforkin{B}{C}_{A}^{D}$ implies $A \le_s B \subseteq D,A \subseteq
C \subseteq D,B \cap C = A$. \newline
(f) (base increasing): if $\nonforkin{B}{C}_{A}^{D}$ and $A \le C_1 \le C$
\underbar{then} $\nonforkin{B \cup C_1}{C}_{C_1}^{D}$. \newline
(g)  monotonicity: $A \subseteq B' \le_s B,A \subseteq C' \subseteq C$
and $\nonforkin{B}{C}_{A}^{D}$ implies $\nonforkin{B'}{C'}_{A}^{D}$. \newline
\smallskip
\noindent
(h)  monotonicity: if $B \subseteq D' \subseteq D,C \subseteq D' \subseteq
D$ then $\nonforkin{B}{C}_{A}^{D} \Leftrightarrow \nonforkin{B}{C}_{A}^{D'}$
but in the ``weakly" version for $\Leftarrow$ we add the assumption 
$\boxtimes_{D',D}$ where we let $\boxtimes_{A_1,A_2}$ be defined as
\medskip
\roster
\item "{$\boxtimes_{A_1,A_2}$}"  $\qquad A_1 \subseteq A_2 \in K_\infty$ 
and if $\emptyset \le_i A' \le_s A_2$ then \footnote{of course, if
$A \in {\Cal K}_\infty \Rightarrow \emptyset \le_s A$ this condition
holds trivially.  We expect that in ``reasonable" cases such assumption
can be removed (in clauses (g) and (h))} $A' \subseteq A_1$
\endroster
\medskip

\noindent
(i) existence: if $A \le_s B$ and $A \subseteq C$ \underbar{then} for some $D$
and $f$ we have: $C \subseteq D, f$ is an embedding of $B$ into $D$ over $A$
and $\nonfork{}{}_{}(A,f(B),C,D)$ (but in the ``weakly" version this is
omitted) \newline
\smallskip
\noindent
(j) Right Transitivity: if $\nonforkin{B_0}{A_1}_{A_0}^{B_2}$ and 
$\nonforkin{B_1}{A_2}_{A_1}^{B_2}$ and $B_0 \subseteq B_1$
then $\nonforkin{B_0}{A_2}_{A_0}^{B_2}$ \newline
\smallskip
\noindent
(k) Left Transitivity: If $\nonforkin{A_1}{B_0}_{A_0}^{B_2}$ and 
$\nonforkin{A_2}{B_1}_{A_1}^{B_2}$ and $B_0 \subseteq B_1$
\underbar{then} \newline
$\nonforkin{A_2}{B_0}_{A_0}^{B_2}$ \newline
\smallskip
\noindent
(l) $({\frak K},\nonfork{}{}_{})$ has symmetry which means
\medskip
\roster
\item "{$(*)$}"  if $\nonforkin{B}{C}_{A}^{D}$ and $A \le_s C$ then
$\nonforkin{C}{B}_{A}^{D}$
\endroster
\medskip

\noindent
(m)  Smoothness:  $\nonforkin{B}{C}_{A}^{D}$ implies 
$C \le_s B \cup C$ (of course, $(d)+(f)$ implies this).
\medskip
\noindent
2) We say $({\frak K},\nonfork{}{}_{})$ has the strong finite basis property
when
\smallskip
\noindent
(n) for every $\ell$ for some $m$
\medskip
\roster
\item "{$(*)$}"  if $B \le D,C \le D,|B| \le \ell$ then for some
$A \le C,|A| \le m$ and \newline
$\nonforkin{B \cup A}{C}_{A}^{D}$.
\endroster
\medskip

\noindent
or at least \newline
(n)$^-$  for every $\ell,r$ for some $m$
\medskip
\roster
\item "{$(*)$}"  if $|B| \le \ell,B \le D \in {\Cal K}_\infty,A_i \le D$
for $i \le m$ such that $A_i \subseteq A_{i+1}$, \underbar{then} for some $i$
and $C \le A_i$ we have:

$$
B \cap A_{i+1} \subseteq C \subseteq A_i \text{ and }
\nonforkin{B \cup C}{c \ell^r(C,A_{i+1})}_{C}^{D}
$$
\endroster
\medskip

\noindent
3) We say $\nonfork{}{}_{}$ (or $({\frak K},\nonfork{}{}_{})$) has the
uniqueness if the $\nonfork{}{}_{}$-free amalgamation is unique, that is:
suppose that for $\ell = 1,2$ \,\, $\nonforkin{B_\ell}{C_\ell}_{A_\ell}
^{D_\ell}$ and $D_\ell = B_\ell \cup C_\ell$ and $f$ is an isomorphism from
$B_1$ onto $B_2$ and $g$ is an isomorphism from $C_1$ onto $C_2$ and
$A_1 = A_2,f \restriction A_1 = g \restriction A_2$, 
\underbar{then} $f \cup g$ is an isomorphism fom $D_1$ onto $D_2$. \newline
4) We say $\nonfork{}{}_{}$ (or $({\frak K},\nonfork{}{}_{})$ has dual
transitivity when:
\medskip
\roster
\item "{$(*)$}"  $\nonforkin{A_1}{C_0}_{A_0}^{C_1}$,
$\nonforkin{A_2}{C_1}_{A_0}^{C_2}$ then $\nonforkin{A_2}{C_0}_{A_0}^{C_2}$ but
in the weak case assume $\boxtimes_{C_1,C_2}$.
\endroster
\enddefinition
\bigskip

\fakesubhead{\stag{2.6B} Fact} \endsubhead
1) If ${\frak K}$ is a 0-1 context (see \scite{2.1}, with $<_s$ as in
Definition \scite{2.3}(2)(c)) and $\nonfork{}{}_{} = \{(A,B,C,D):A \le_s
B \subseteq D,A \subseteq C \subseteq D,B \cap C = A$ and the quadruple is
freely amalgamated in the sense of \scite{1.2}(6) (no new instances of the 
relations)$\}$ \underbar{then}
$({\frak K},\nonfork{}{}_{})$ is a 0-1 context except possibly 
$(f)$(base increasing),$(i)$(existence), \newline
$(m)$(smoothness); with uniqueness (see \scite{2.6A}(3)). \newline
2) Assume $({\frak K},\nonfork{}{}_{})$ is a 0-1 context
\medskip
\roster
\item "{$(a)$}"  if $\nonforkin{B}{C}_{A}^{D},D = B \cup C,
B' \le B,A \le A^+ \le C$, \newline
$r = k +|A|,c \ell^r(A,C) \subseteq A^+$ \underbar{then} 
$c \ell^k(B',D) \subseteq B \cup A^+$
\smallskip
\noindent
\item "{$(b)$}"  $\nonforkin{B}{C}_{A}^{D},A \le_i D' \le D,
D' \subseteq B \cup C$ \underbar{then} $A \le_i D' \cap C$.
\endroster
\bigskip

\demo{Proof}  1) Straight. \newline
\underbar{clauses (a)-(e)}:  By the definitions.
\medskip

\noindent
\underbar{clause (g)}:  Check the definitions, and by
\scite{2.6}(6) we have $A \le_s B'$.
\medskip

\noindent
\underbar{clauses (h),(j)}:  Check.
\medskip

\noindent
\underbar{clause (k)}:  The point is that $A_0 \le_s A_1 \le_s A_2
\Rightarrow A_0 \le_s A_2$ by \scite{2.6}(10).
\medskip

\noindent
\underbar{clause (l)}:  Read the definitions.
\medskip

\noindent
\underbar{uniqueness}:  Reflect on the meaning of $\nonfork{}{}_{}$.
\medskip

\noindent
2) \underbar{clause (a)}:  By monotonicity of $\nonfork{}{}_{}$ (see
Definition \scite{2.6A}(1)(g) we may assume that $B = B'$.  Assume toward
contradiction that the conclusion fails.  So there are $C' \le D,|C'| \le k,
d \in C' \backslash B \backslash A^+ \subseteq C \backslash A^+$ and 
$C' \cap B <_i C'$, hence (see \scite{2.6}(3)), $B <_i B \cup C'$.  Let 
$C_1 = A \cup (C' \cap C)$, and let $C_0$ be such that $A \le_i C_0 \le_s C_1$
(exists by \scite{2.6}(4)).
By clauses (f) + (g) of Definition \scite{2.6A}(1) we have
$\nonforkin{B \cup C_0}{C_1}_{C_0}^{D}$, by symmetry (i.e. \scite{2.6A}(1)
(l)) we have $\nonforkin{C_1}{B \cup C_0}_{C_0}^{D}$ hence by smoothness
\scite{2.6A}(1)(m) $B \cup C_0 \le_s B \cup C_1$.  
So as $C' \cap B <_i C' \subseteq B \cup C_1$ (as $C_1 = A \cup (C' \cap C)$
and $D = B \cup C)$ by \scite{2.6}(3) we necessarily have $C' \subseteq B \cup
C_0$, so as $d \in C' \backslash B$ necessarily $d \in C_0$, but
$|C_0| \le |A| + |C' \cap C \backslash A| \le |A| + k = r$.  Remember $A \le_i
C_0$, so $C_0 \subseteq c \ell^r(A,C) \subseteq A^+$, hence $d \in A^+$ a
contradiction.
\medskip

\noindent
\underbar{clause (b)}:  If the conclusion fails, then for some 
$C_1,A \le_i C_1 <_s D' \cap C$.
By monotonicity (= clause (g) of Definition \scite{2.6A}(1)) we have
$\nonforkin{B}{D' \cap C}_{A}^{D}$, so by smoothness (= clause (m) of
Definition \scite{2.6A}(1) we have $D' \cap C \le_s B \cup (D' \cap C)$ but
$\le_s$ is transitive (by \scite{2.6}(10)) so $C_1 <_s B \cup (D' \cap C)$
but $C_1 \subseteq (B \cup C) \cap D' \subseteq B \cup (D' \cap C)$ (as
$C_1 \subseteq (D' \cap C))$ so $C_1 <_s (B \cup C) \cap D' = D'$) (as
$D' \subseteq  B \cup C$) but this contradicts $A \le_i D',A \subseteq C_1$. 
 \hfill$\square_{\scite{2.6B}}$
\enddemo
\bigskip

\definition{\stag{2.6D} Definition}  Let $({\frak K},\nonfork{}{}_{})$ be a
0-1 context (or weakly 0-1 context).  If \nl
$A \le B_\ell \le D$ for $\ell < m$,
we say $\{B_\ell:\ell < m\}$ is free over $A$ inside $D$ (or the $B_\ell$'s
are free over $A$ inside $D$) if for every $\ell,A \le_s B_\ell$ and
$B_\ell \nonfork{}{}_{A} \dsize \bigcup_{k < \ell} B_k$. \newline
Justification is:
\enddefinition
\bigskip

\proclaim{\stag{2.6D1} Claim}  The order (in Definition \scite{2.6D})
is immaterial.
\endproclaim
\bigskip

\demo{Proof}  if $\ell +1 < m$, then (letting $B^*_i = \dsize \bigcup_{j<i}
B_j \cup A$), $B^*_\ell \cup B_{\ell +1} \nonfork{}{}_{B^*_\ell} 
B^*_{\ell +1}$ (by clause (f) of Definition \scite{2.6A}) hence
$B^*_{\ell +1} \nonfork{}{}_{B^*_\ell} B^*_\ell \cup B_{\ell +1}$ (by
clause (l), symmetry, of Definition \scite{2.6A}) that is
$B^*_\ell \cup B_\ell \nonfork{}{}_{B^*_\ell} B^*_\ell \cup B_{\ell +1}$.
We have $\nonfork{B_\ell}{B^*_\ell}_{A}$ so $\nonforkin{B_\ell}{B^*_\ell}
_{A}^{B^*_{\ell +2}}$ (by clause (g) of Definition \scite{2.6A}(1)), and
similarly $\nonforkin{B^*_\ell \cup B_\ell}{B^*_\ell \cup B_{\ell +1}}
_{B^*_\ell}^{B^*_{\ell +2}}$ hence we get 
$\nonforkin{B_\ell}{B^*_\ell \cup B_{\ell +1}}_{A}^{B^*_{\ell +2}}$ (using
clause (j), right transitivity from Definition \scite{2.6A} with
$A,B_\ell,B^*_\ell,B^*_\ell \cup B_\ell,B^*_\ell \cup B_{\ell +1},
B^*_{\ell+2}$ here standing for $A_0,B_0,A_1,B_1,A_2,B_2$ there).

We get that we can permute the order in $\langle B_k:k < m \rangle$ of
$\ell$ and $\ell +1$, but such permutations generates the permutation group
on $\{0,\dotsc,m-1\}$ so we are done. \nl
${{}}$  {\hfill$\square_{\scite{2.6D1}}$
\enddemo
\bn
Now concerning \scite{2.6E} we can say more
\proclaim{\stag{2.6Ea} Claim}  Assume $({\frak K},\nonfork{}{}_{})$ is a 0-1
context which is weakly nice. \nl
1) If $A < B \in {\Cal K}_infty$ then the following are equivalent
\mr
\item "{$(a)$}"  $A <_s B$
\sn
\item "{$(b)$}"  $1 = \text{ Lim}_n\text{ Prob}_{\mu_n}$(for every embedding
$f$ of $A$ into ${\Cal M}_n$ \nl

$\qquad \qquad \qquad$ there are embeddings $g_\ell:B \rightarrow {\Cal M}_n$
extending \nl

$\qquad \qquad \qquad f$ for $\ell < m$ such that $\langle g_\ell:\ell < m
\rangle$ is disjoint over $A$)
\sn
\item "{(c)}"  for every $m < \omega$, for every $k < \omega$ \newline
$1 = \text{ Lim}_n \text{ Prob}_{\mu_n}\biggl($for every embedding $f$ 
of $A$ into ${\Cal M}_n$ there are \newline

$\qquad \qquad \qquad \qquad$ embeddings $g_\ell:B \rightarrow {\Cal M}_n$ 
extending $f$, for $\ell < m$ \newline
 
$\qquad \qquad \qquad \qquad$ such that
$\langle g_\ell:\ell < m \rangle$ is disjoint over $A$ and for \newline

$\qquad \qquad \qquad \qquad \ell < m$ we have $c \ell^k(f(A),{\Cal M}_n)$ 
and the $g_\ell(B)$'s are \newline

$\qquad \qquad \qquad \qquad \nonfork{}{}_{}$-free over 
$g_\ell(A)$ inside ${\Cal M}_n\biggr)$.
\ermn
2)  If $A_0 \le A_1 \in {\Cal K}_\infty$ and $A_0 \le_s B$ \underbar{then}
there are $C \in {\Cal K}_\infty$ and $f$ such that $A_1 \le C,f$ is an 
embedding of $B$ into $C$ over $A$, say $B_1 = f(B)$ and 
$B_1 \cap A_1 = A_0$ and $A_1 \le_s B_1 \cup A_1$.
\endproclaim
\bigskip

\demo{Proof}  1) By \scite{2.6E} clearly $(a) \Leftrightarrow (b)$. \nl
1)  Trivially $(c) \rightarrow (b)$.  Let us prove $(a) \rightarrow (c)$. 
\newline
Let $m,k \in \Bbb N$.  We now choose by induction on $\ell \le m,D_m$ and if
$\ell < m$ also $g_\ell$ such that: $D_0 = A,D_\ell \le D_{\ell +1},
g^0_\ell$ is an embedding of $B$ over $A$ into $D_{\ell +1}$ such that
$\nonforkin{g^0_\ell(B)}{D_\ell}_{A}^{D_{\ell +1}}$, the induction step is
done by clause (i)(existence) of Definition \scite{2.6A}(1).  By clause (h)
of Definition \scite{2.6A}(1) without loss of generality $D_{\ell +1} =
D_\ell \cup g_\ell(B)$.  By clause (m)(smoothness) of Definition
\scite{2.6A}(1) we know that $D_\ell \le_s D_{\ell +1}$ and by clause (e) we
know that $g^0_\ell(B) \cap D_\ell = A$.  So $A \le_s D_m$ so by
``${\frak K}$ weakly nice", if ${\Cal M}_n$ is random enough and $f:A
\rightarrow {\Cal M}_n$ an embedding then we can find an embedding
$g:D_m \rightarrow {\Cal M}_n$ extending $f$.  We let $g_\ell = g \circ
g^0_\ell$ for $\ell < m$.  Now for $\ell_1 < \ell_2 < m,g_{\ell_1}(B) \cap
g_{\ell_2}(B) = g(g^0_{\ell_1}(B) \cap g^0_{\ell_2}(B)) \subseteq
g(g^0_{\ell_1}(B) \cap D_{\ell_2}) = g(g^0_{\ell_1}(B)) = g_{\ell_1}(B)$,
so the disjointness demand holds.  The freeness holds by the construction,
and \scite{2.6D}(1). \nl
\sn
2)  Without loss of generality, (not embedded yet in one model!)
$B \cap A_1 = A_0$.  Let $r^*$ be the number 
of structures $C \in {\Cal K}_\infty$ with set of elements $B \cup A_1$ such
that $B \subseteq C \and A_1 \subseteq C$.  Now let $n_0(\varepsilon)$ be 
such that as $A_1 \in {\Cal K}_\infty$ clearly for some $\varepsilon^* > 0$
\medskip
\roster
\item "{$(*)_2$}"  for arbitrarily large natural number $n$: \newline
$\varepsilon^* \le \text{ Prob}_{\mu_n} \biggl($there is an embedding 
$f$ of $A_1$ into ${\Cal M}_n\biggr)$.
\endroster
\medskip

\noindent
Let (the function) $n_1(\varepsilon)$ be such that:
\medskip
\roster
\item "{$(*)_3$}"  for every $\varepsilon \in \Bbb R^+$,
if $n \ge n_1(\varepsilon)$ then: \newline
$1 - \varepsilon / 4 \le \text{ Prob}_{\mu_n}\biggl($ if $C$ is
embeddable into ${\Cal M}_n$ and $C$ has at most \newline

$\qquad \qquad \qquad \qquad |B| + |A_1|$ elements then 
$C \in {\Cal K}_\infty\biggr)$.
\endroster
\medskip

\noindent
As $A_0 \le_s B$ by part (1) for some function $n_2(\varepsilon,m)$ we have:
\medskip
\roster
\item "{$(*)_4$}"  for every $\varepsilon \in \Bbb R^+$ and $m \in \Bbb N$,
for every $m \ge n_2 (\varepsilon,m)$ we have: \newline
$1 - \varepsilon / 4 \le \text{ Prob}_{\mu_n}\biggl($every embedding $f$ of
$A_0$ into ${\Cal M}_n$ has $m$ extensions $g$ \newline

$\qquad \qquad \qquad \qquad$ to embeddings of $B$ into ${\Cal M}_n$,
pairwise disjoint \newline

$\qquad \qquad \qquad \qquad$ over $A\biggr)$.
\endroster
\medskip

\noindent
Now by \scite{2.6}(11)
\medskip
\roster
\item "{$(*)_5$}"  for $C \in {\Cal K}$ with the set of elements of $C$ being
$B \cup A_1$, such that \newline
$(A_1 \subseteq C) \and \neg(A_1 \le_s C)$ there is $\ell_3 = \ell_3(C)$
such that for each $\varepsilon \in \Bbb R^+$ for some 
$n = n_3(C,\varepsilon)$
we have: for every $n \ge n_3(C,\varepsilon)$, \newline
$1 - \varepsilon /4 \le \text{ Prob}_{\mu_n}\biggl($
there is no sequence $\langle g_\ell:\ell < \ell_3 \rangle$ of embedding of
$C$ \newline

$\qquad \qquad \qquad \qquad$ into ${\Cal M}_n$ pairwise disjoint 
over $A_1\biggr)$.
\endroster
\medskip

\noindent
Let $\varepsilon \in \Bbb R^+$ be given.
\mn
Let $\ell^* = \text{ max}\{\ell_3(C):C \in {\Cal K} \text{ has
set of elements } B \cup A_1 \text{ and extends } B \text{ and } A_1$
\newline

$\qquad \qquad \qquad \quad$ but $\neg(A_1 \le_s C)\}$.
\mn
Let $m^* = r^* \times \ell^* + |A_1 \backslash A_0| + 1$.
\mn
Let $n^*(\varepsilon) = \text{ Max}\{n_1(\varepsilon),
n_2(\varepsilon,m^*),n_3(C,\varepsilon):C \in {\Cal K}$ \newline

$\qquad \qquad \qquad \qquad \qquad \text{ has universe } B \cup A_1
\text{ and extends } B \text{ and } A_1$ \newline

$\qquad \qquad \qquad \qquad \qquad$ but $\neg(A_1 \le_s C)\}$.
\mn
Let $\varepsilon < \varepsilon^*$ and let $n \in Y =: \{n:\text{ the
statement } (*)_2 \text{ holds for } n\}$. \nl
For each $M \in {\Cal K}_n$ choose if possible an embedding $f^M$ of $A_1$ 
into $M$. \newline

Now if $f^M$ is defined, choose if possible a sequence $\langle g^M_m:m <
m^* \rangle$ of embeddings of $B$ into $M$ extending $f^M \restriction A_0$
and pairwise disjoint over $A$.  As $\langle g^M_m(B \backslash A_0):m < m^*
\rangle$ is a sequence of pairwise disjoint subsets of $M$, the
set $u =: \{m < m^*:g^M_m(B \backslash A_0)$ is disjoint to
$f^M(A_1 \backslash A_0)$, equivalently $g^M_m(B) \cap f^M(A_1) = f^M(A_0)\}$
has at least $m^* - |A_1 \backslash A_0|$ members which is
$r^* \times \ell^* +1$.

For $m \in u$ let $C_m$ be the model with universe $B \cup A_1$ such that
$g^M_m \cup f^M$ is an isomorphism from $C_m$ into $M
\restriction (g^M_m(B) \cup f^M(A_1))$ (note: this function is one to one
as $n \in u$ and by the definition of $u$).  By the choice of $r^*$ for some
model $C$ with set of elements $B \cup A_1$ we have $|\{ m \in u:
C = C_m\}| \ge |u|/r^*$.  But $|u| \ge m^* - |A_1 \backslash
A_0| = r^* \times \ell^* + 1$ hence $u' =: \{m \in u_m:C = C_m\}$ has $>
\ell^*$ members.

Now with probability $\ge \varepsilon^* - \varepsilon > 0,{\Cal M}_n$ 
satisfies demands
$(*)_2 - (*)_5$ above hence $f^{{\Cal M}_n}$ is well defined (by $(*)_2$)
and $\langle g^{{\Cal M}_n}_m:m < m^* \rangle$ is well defined
(by $(*)_4$ and the choice of $n(\varepsilon)$).  By $(*)_3$, we know
$C \in {\Cal K}_\infty$, obviously $A_1 \le C,B \le C$ (as each $C_m$
satisfies this).  Lastly by $(*)_5$ the sequence 
$\langle g^{{\Cal M}_n}_m:m < m^* \rangle$
witnesses $A_1 \le_s C$, so we have finished. \hfill$\square_{\scite{2.6E}}$
\enddemo
\bigskip

\fakesubhead {\stag{2.6F} Comment}  \endsubhead
Our real interest is in instances of ${\frak K}$ but $\nonfork{}{}_{}$ is
central in the following way
\medskip
\roster
\item "{$(a)$}"   it is a way to express assumptions on ${\frak K}$ helping
to analyze the limit behaviour (for which having a 0-1 law is a reasonable
criterion)
\smallskip
\noindent
\item "{$(b)$}"  assuming the random enough ${\Cal M}_n$ satisfies $(*)$
below we can define $\nonfork{}{}_{}$ and prove $({\frak K},\nonfork{}{}_{})$
is a 0-1 context where
{\roster
\itemitem{ $(*)$ }  for quantifiers free $\varphi(\bar x,\bar y)$, the numbers
$|\varphi({\Cal M}_n,\bar b)|$ for $\bar b \in {}^{\ell g(\bar y)}
({\Cal M}_n)$ behave regularly enough; where \newline
$$
\varphi({\Cal M}_n,\bar b) = \{\bar a \in {}^{\ell g(\bar x)}({\Cal M}_n):
{\Cal M}_n \models \varphi[\bar a,b]\}.
$$
\endroster}
\endroster
\bigskip

\centerline {$* \qquad * \qquad *$}
\bigskip

\proclaim{\stag{2.18} Claim}  1) If 
${\Cal K}_n = \{([n])\}$ (so $\mu_n$ trivial)
and $\nonforkin{B}{C}_{A}^{D}$ means \newline
$B \cap C = A \and B \cup C \subseteq
D$ \underbar{then} $({\frak K},\nonfork{}{}_{})$ is explicitly nice 0-1 
context.  Also $A \le_s B$ iff $A \subseteq B$ and ${\Cal K}_\infty$ is the
family of finite models. \newline
2) Let ${\Cal K}_n = \{([n],S,c_f,c_\ell)\}$ with $c_f = 1,c_\ell = n$
and $S$ the successor relation 
(so $\mu_n$-trivial).  Then $A <_s B$ means $A \subseteq B \and
(\forall x \in A)(\forall y \in B \backslash A)[\neg x S^By \and \neg y S^B y
\and \neg y S^B x]$ and $A \in {\Cal K}_\infty$ iff $A$ is isomorphic to
$([n],S,c_f,c_\ell) \restriction X$ for some $X \subseteq [n]$. \newline
3) Let in 2), $\nonforkin{B}{C}_{A}^{D}$ means
$A \le_s B \le D \in {\Cal K}_\infty,A \le C \le D,B \cap C = A$ and \newline
$(\forall x \in C)(\forall y \in B \backslash A)
[\neg xS^D y \and \neg y S^D x]$.  \underbar{Then} 
$({\frak K},\nonfork{}{}_{})$ is 0-1 context, ${\frak K}$ is
explicitly almost nice.
\endproclaim
\bigskip

\demo{Proof}  Straight, e.g.  Why \scite{2.6A}(f)? \newline
3)  Let us check clause (f) of Definition \scite{2.6A}(1), (base increasing)
so we have $\nonforkin{B}{C}_{A}^{D}$ and $A \le C_2 \le C$ and we have to
prove $\nonforkin{B \cup C_1}{C}_{C_1}^{D}$.  Now $C_1 \subseteq B \cup C_1$
trivially, and for $C_1 \le_s B \cup C_1$ see the definition of
$\nonfork{}{}_{}$; $B \cup C_1 \subseteq D$ as $B \subseteq D,C_1 \subseteq
C \subseteq D$, and $C_1 \subseteq C$ holds by assumption as $C \subseteq D$
as $\nonforkin{B}{C}_{A}^{D}$.  Lastly, assume $x \in C_1,y \in B \cup
C_1 \backslash C_1$ then $x \in C,y \in B \backslash C$ so by the definition
of $\nonforkin{B}{C}_{A}^{D}$ we have $\neg xS^D y \and \neg y S^Dx$, as
required.  \hfill$\square_{\scite{2.18}}$ 
\enddemo
\bigskip

\centerline {$* \qquad * \qquad *$}
\bigskip

We connect the free amalgamation with \S2, i.e. context ${\frak K}$ which
obeys $\bar h$ which is bound by $h_1$, this gives a natural definition for
free amalgamation. \nl
Continuing Definition \scite{3.1} let
\definition{\stag{3.2a} Definition}  1) For a 0-1 context ${\frak K}$ 
obeying $\bar h$ and function $h$ from
$\dsize \bigcup_n {\Cal K}_n$ to $\Bbb R^+$ we define a four place relation
$\nonfork{}{}_{} = \nonfork{}{}_{h} = \nonfork{}{}_{}[h] =
\nonfork{}{}_{}[\bar h,h]$ on ${\Cal K}_\infty$:
\newline
$\nonfork{}{}_{h}(A,B,C,D) \text{ \underbar{iff}}$
\medskip
\roster
\item "{$(a)$}"   $A \le_s B \le D,A \le C \le D,A = B \cap C$ and,
\sn
\item "{$(b)$}"   letting $D_1 = D \restriction B \cup C$ we have 
$C \le_s D_1$ and
\sn
\item "{$(c)$}"  for every $\varepsilon > 0$ for every random enough 
${\Cal M}_n$ into which $D$ is embeddable we have \footnote{inequality is the
other direction follows from the definitions as proved in clause $(\zeta)$
of \scite{3.4}(2)}
\endroster
\medskip

$$
h^u_{C,D_1}[{\Cal M}_n] \times h[{\Cal M}_n]^\varepsilon \ge h^u_{A,B}
[{\Cal M}_n].
$$
\medskip

\noindent
note that (b) follows by (a), (c)).
\bigskip

We define $\nonfork{}{}_{} = \nonforkin{}{}_{h}^{-} = \nonfork{}{}_{}^{-}
[h] = \nonfork{}{}_{}^{-}[\bar h,\bar h]$ similarly:
$\nonfork{}{}_{}[A,B,C,D]$ iff
\mr
\item "{$(a)$}"  $A \le_s B \le D,A \le C \le D,A = B \cap D$ and
\sn
\item "{$(b)$}"  letting $D_1 = D \restriction (B \cup C)$ we have
$C \le_s D_1$
\sn
\item "{$(c)$}"  for some $m = m(A,B,C,D) \in \Bbb N \backslash \{0\}$,
for every random enough ${\Cal M}_n$ we have
$$
h^u_{C,D_1}[{\Cal M}_n] \times h[{\Cal M}_n]^m \ge h^u_{A,B}[{\Cal M}_n].
$$
\ermn
2)  Let $\nonfork{}{}_{}[\bar h] = \nonfork{}{}_{}[\bar h,h']$ for 
$h'(M) = \|M\|$ (see \scite{3.6A}(3) below).
\enddefinition
\bn
Continuing \scite{3.4} we note
\proclaim{\stag{3.4a} Claim}  In \scite{3.4}(2) we can add
if $\nonforkin{B}{C}_{A}^{D}$ (so $A <_s B \le D,A \le C \le_s D$) and
$D = B \cup C$ \underbar{then} for every $\varepsilon > 0$ for every 
random enough 
${\Cal M}_n$ we have $h[{\Cal M}_n]^\varepsilon \ge h^u_{A,B}[{\Cal M}_n]
/h^u_{C,D}[{\Cal M}_n] \ge 1/(h[{\Cal M}_n])^\varepsilon$ 
(when $C$ is embeddable into ${\Cal M}_n$ of course).
\endproclaim
\bigskip

\proclaim{\stag{3.5} Claim}  Assume ${\frak K}$ obeys $\bar h$ which is 
bounded by $h$ and $\nonfork{}{}_{} = \nonfork{}{}_{\bar h,h}$. \newline
1) $({\frak K},\nonfork{}{}_{})$ is a weak 0-1 context (see Definition
\scite{2.6A}(1) second version). \newline
2) If in addition $(*)$ below holds, \underbar{then} 
$({\frak K},\nonfork{}{}_{})$ is a 0-1 context
\medskip
\roster
\item "{$(*)$}"  if $A <_{pr} B \le D,A < C <_{pr} D,B \cup C = D$ and
$B,C$ are not $\nonfork{}{}_{}$-free by amalgamation over $A$ inside $D$, 
\underbar{then} for some $\varepsilon > 0$ for every random enough 
${\Cal M}_n$, we have
$h^d_{A,B}[{\Cal M}_n]/h^u_{C,D}[{\Cal M}_n] \ge (h[{\Cal M}_n])^\varepsilon$.
\endroster
\medskip

\noindent
3) If condition $(*)$ from part (2) holds then also
\medskip
\roster
\item "{$(*)^+$}"  if $A <_s B \le D,A \le C \le_s D,D = B \cup D$ and
$B,C$ are not $\nonfork{}{}_{}$-free by amalgamation over $A$ inside 
$D$ \underbar{then} for
some $\varepsilon > 0$ for every random enough ${\Cal M}_n$ we have 
$$
h^d_{A,B}[{\Cal M}_n]/h^u_{C,D}[{\Cal M}_n] \ge (h[{\Cal M}_n])^\varepsilon,
$$
\endroster
\endproclaim
\bigskip

\remark{Remark}  1) Note that $(*)$ of \scite{3.5}(2) just say a 
dichotomy; i.e. the condition (c) in Definition \scite{3.1}(5), 
either holds for every random enough ${\Cal M}_n$, or fails for every 
random enough ${\Cal M}_n$. \nl
2) Note that $(*)$ of \scite{3.5} exclude the case of having a successor
relation.
\endremark
\bigskip

\demo{Proof} The order is that first we prove \scite{3.5}(3) (and inside it a
restricted version of transitivity (clause (j) of Definition \scite{2.1}(2))),
and only then we prove \scite{3.5}(1)+(2) by going over all the clauses of
Definition \scite{2.1}(2), so in some clauses there is a difference between
\scite{3.5}(1) and \scite{3.5}(2), and then in the second case we may use
$(*)$ of \scite{3.5}(2) (and so \scite{3.5}(3)). \newline
3) Let $\bar A = \langle A_i:i \le k \rangle$ be a decomposition of
$A <_s B$, so $A_i <_{\text{pr}} A_{i+1}$ for $i < k$.  Let 
$C_i = A_i \cup C$, hence for every $i < k,C_i \le_i C_{i+1}$ or 
$C_i <_{\text{pr}} C_{i+1}$. \newline
[Why?  By \scite{2.6}(4) + \scite{2.6}(8)].
\mn
For each $i < k$ let $\varepsilon(i) \in \Bbb R^{\ge 0}$ be such that:
\mn
\ub{Case 1}:  If $c_i \le_i c_{i+1}$ then $\varepsilon(i) > 0$ and
$h^u_{A_i,A_{i+1}}[{\Cal M}_n] \ge (h[{\Cal M}_n])^{\varepsilon(i)}$ for
every ${\Cal M}_n$ random enough.
\mn
\ub{Case 2}:  If $\neg(c_i \le_i c_{i+1})$ hence $c_i <_{pr} c_{i+1}$ but
$neg\bigl( \nonforkin{A_{i+1}}{C_i}_{A_i}^{C_{i+1}} \bigr)$ then
$\varepsilon(i) > 0$ and
$h^u_{A_i,A_{i+1}}[{\Cal M}_n] / h^d_{c_i,c_{i+1}}[{\Cal M}_n]
\ge h[{\Cal M}_n])^{\varepsilon(i)}$.
\mn
\ub{Case 3}:  If neither Case 1 nor Case 2, then $\varepsilon_i = 0$.
Note that by \scite{3.4B}(2) for ${\Cal M}_n$ random enough
$h^u_{A_i,A_{i+1}}[{\Cal M}_n] \ge h^d_{C_i,C_{i+1}}[{\Cal M}_n]$
(if $D$ embeddable into ${\Cal M}_n$).  Let $w_\ell = \{i < k:
\text{case } \ell \text{ occurs}\}$.  \nl
First assume $\dsize \sum_{i<k} \varepsilon(i) > 0$, and let $\zeta \in
\Bbb R^+$ be $< \dsize \sum_{i<k} \varepsilon(i)/2$. \nl
So let ${\Cal M}_n$ be random enough such that $C$ is embeddable into
${\Cal M}_n$ (hence $D$ is embeddable into ${\Cal M}_n$ hence also the
$C_i$'s are).
\medskip
\roster
\item "{$(a)$}"  $h^u_{A_i,A_{i+1}}[{\Cal M}_n] \ge h^d_{C_i,C_{i+1}}
[{\Cal M}_n]$ when $C_i <_{\text{pr}} C_{i+1}$ \nl
[Why?  See \scite{3.4B}(2).]
\smallskip
\noindent
\item "{$(b)$}"  $1 \le \dsize \prod_{i < k} h^u_{A_i,A_{i+1}}[{\Cal M}_n]
/h^u_{A,B}[M] \le (h{\Cal M}_n])^{\zeta/k}$. \nl
[Why?  The first inequality by the definition in \scite{3.1}(4), as
$h^{-u}_{A,B}(M) = h^u_{A,B}(M)$.  Second inequality by the second inequality
in \scite{3.4}(2)$(\delta)$.]
\endroster
\mn
So

$$
\split
h^d_{C,D}[{\Cal M}] &\le \\
  &\qquad \qquad \qquad \qquad \qquad \qquad \text{[by } \scite{3.4}(2)
(\varepsilon)]
\endsplit
$$

$$
\split
\dsize \prod_{i \in w_2 \cup w_3} h^d_{C_i,C_{i+1}}[{\Cal M}_n] &\le ? \\
  &\qquad \qquad \qquad \qquad \qquad \qquad \text{[trivially]}
\endsplit
$$

$$
\split
(\dsize \prod_{i \in w_1} 1) &\times(\dsize \prod_{i \in w_2} 
h^d_{C_i,C_{i+1}}[{\Cal M}_n]) \\
   &\times \dsize \prod_{i \in w_3} h^d_{C_i,C_{i+1}}[{\Cal M}_n]) \le ? \\
  &\qquad \qquad \qquad \qquad \text{[by the statements in the cases]}
\endsplit
$$

$$
\split
(\dsize \prod_{i \in w_1} h^u_{A_i,A_{i+1}}[{\Cal M}_n]
/(h[{\Cal M}_n])^{\varepsilon(i)} &\times (\dsize \prod_{i \in w_2} 
h^u_{A_i,A_{i+1}}[{\Cal M}_n]/h[{\Cal M}_n]^{\varepsilon(i)}) \\
  &\times (\dsize \prod_{i \in w_i} h^u_{A_i,A_{i+1}}[{\Cal M}_n]) \\
  &= (\dsize \prod_{i \le k} h^u_{A_i,A_{i+1}}[{\Cal M}_n]) \times
     (h[{\Cal M}_n])^{- \Sigma \varepsilon(i)} \\
  &\le (\dsize \prod_{i \le k} h^u_{A_i,A_{i+1}}[{\Cal M}_n]) \times
   (h[{\Cal M}_n])^\zeta \times (h[{\Cal M}_n])^{- \Sigma \varepsilon(i)} \\
  &\qquad \qquad \qquad \qquad \text{[trivially]}
\endsplit
$$
\mn
So $\dsize \sum_{i < k} \varepsilon(i) - \zeta > \Bbb R^+$ can serve as the
desired $\varepsilon$.
\sn
Second, assume $\dsize \sum_{i<k}\varepsilon(i) = 0$ then, 
as $C_k = D$, we have
\medskip
\roster
\item "{$(c)$}"  $1 \le \dsize \prod_{i < k} h^u_{C_i,C_{i+1}}[M]/h^u_{C,D}[M]
\le (h[M])^\varepsilon$ \newline
and by Definition \scite{3.1}(5)
\sn
\item "{$(d)$}"  for every $\varepsilon \in \Bbb R^+$ for every random
enough ${\Cal M}_n$ into which $C_{i+1}$ is embeddable
$$
(h[{\Cal M}_n])^{-\varepsilon}1 \le h^u_{A_i,A_{i+1}}[{\Cal M}_n]
/h^u_{C,D}[{\Cal M}_n] \le (h[{\Cal M}])^\varepsilon.
$$
\medskip

\noindent
Using (b)+(c)+(d) we can get ``$B,C$ are $\nonfork{}{}_{}$-free over $A$
inside $D$" so an assumption of $(*)^+$ fail.
\endroster
\enddemo
\bigskip

\fakesubhead{Proof of (1)+(2)} \endsubhead  
The proof is split according to the clause in Definition \scite{2.6A}(1); i.e.
to $(d),\dotsc,(m)$.
\medskip

\noindent
\underbar{Clause $(d)$}: \newline

Trivial.
\mn
\underbar{Clause $(e)$}: \newline

Trivial.
\mn
\underbar{Clause $(f)$}: \ub{Monotonicity in $C$} \newline

The point to notice is that ``$D$" is embeddable into ${\Cal M}_n$" does not
imply ``$D$ is embeddable into ${\Cal M}_n$" but read \scite{3.1}(3),
\scite{3.1}(3A).
\bn
\centerline {$* \qquad * \qquad *$}
\bn
\underbar{Clause $(h)$}:  (Reading Definition \scite{3.1}(5)), 
but here there is a difference
between \scite{3.5}(1) and \scite{3.5}(2).  We have to restrict ourselves to
models into which $D$ can be embeddable.  But there may be many into which
$D'$ is embeddable but $D$ is not.

As $D \in {\Cal K}_\infty$, for some $\zeta \in \Bbb R^+$ for arbitrarily
large $n$, the probability that ${\Cal M}_n$ satisfies the desired 
inequalities is $\ge \zeta$, so for \scite{3.5}(2), the assumption $(*)^+$
there guarantees they hold for every random enough ${\Cal M}_n$.

For the weakly version in Definition \scite{2.6}(1) clause $(h)$ we have 
restricted ourselves to the case $\boxtimes_{D',D}$ (see Definition
\scite{2.6}(1)). 
\bn
\ub{Symmetry}:  So assume $\nonforkin{B}{C}_{A}^{D},A <_s C$, and we should
prove $\nonforkin{C}{C}_{A}^{D}$, \nl
without loss of generality $A \ne B,A \ne C$.  So we have
$A <_s B,B \cap C = A,C <_s B$.  We should prove $B <_s B \cup C$, and the
inequality.  Let $m^* \in \Bbb N$.  let $\varepsilon \in \Bbb R^+$ be small
enough.

So let ${\Cal M}_n$ be random enough and $f_0$ an embedding of $A$ into
${\Cal M}_n$.  So

$$
F_1 = \{g_1:f_1 \text{ is an embedding of } C \text{ into } {\Cal M}_n
\text{ extending } f_0\}.
$$
\mn
Clearly $h \ge h^d_{A,C}[{\Cal M}_n]$ and for each $f_1 \in F_1$,

$$
\align
F^2_{f_1} = \{f_2:&f_2 \text{ is an embedding of } C \cup B \text{ into} \\
  &{\Cal M}_n \text{ extending } f_1\}
\endalign
$$
\mn
has $\ge h^d_{C,C \cup B}[{\Cal M}_n]$ members. \nl
So $F_2 = \cup\{F^1_{f_1}:f_1 \in F_1\}$ has $\ge h^d_{A,C}[{\Cal M}_n] \times
h^d_{C,B \cup C}[{\Cal M}_n]$ member.  Consider the mapping $G_2:F_2
\rightarrow F^1 =:\{f^1:f^1$ is an embedding of $B$ into ${\Cal M}_n$
extending $f_0\}$. \nl
So $G_2$ is a mapping from $F_2$ into $F^1$.  Also $F_1$ has at most
$h^u_{A,B}[{\Cal M}_n]$ which is $\ge h^u_{C,B \cup C}[{\Cal M}_n]/
(h[{\Cal M}_n])^\varepsilon$ as $\nonfork{}{}_{}(A,B,C,D)$.  So the number
$|\{f_2 \in F^2:f_2 \restriction B = f^1\}|$ averaging on all $f^1 \in F^1$
is $\ge h^d_{A,C}[{\Cal M}_n]/(h[{\Cal M}_n])^\varepsilon$ which is
$> m^*$ (as ${\Cal M}_n$ is random enough).  So for some $f^1 \in F^1$ the
actual number is $> m^*$, hence $B <_s B \cup C$.  Similarly we get the
inequality.
\bn
\underbar{Clause $(g)$}: \newline

\noindent
I.e. we assume $A \subseteq B' \le_s B,A \subseteq C' \subseteq C$ and
$\nonforkin{B}{C}_{A}^{D}$ (and we should prove $\nonforkin{B'}{C'}_{A}^{D}$).
Now $\nonforkin{B}{C}_{A}^{D}$ means that $A \le_s B \le D,A \le C \le D,
B \cap C = A$ (this is clause (a) of \scite{3.11}(5)) and letting 
$D_1 = B \cup C$, also $C \le_s D_1$ (this is clause (b) of \scite{3.1}(5))
and similarly $C' \le_s B \cup C'$ and for every $\varepsilon_1 > 0$ for 
every random enough ${\Cal M}_n$ we have $1 \le h^u_{A,B}[{\Cal M}_n]/
h^u_{C,D_1}[{\Cal M}_n] < (h[{\Cal M}_n])^{\varepsilon_1}$ (this is clause (c)
of Definition \scite{3.7}(5)).
Looking at the desired conclusion (in particular the definition of
$\nonfork{}{}_{}$) we can restrict ourselves to
${\Cal M}_n$ into which $D$ is embeddable; hence $C',C$ are embeddable.
\newline
Without loss of generality $B' = B$ or $C' = C$.
\bigskip

In the first case, we note that for random enough ${\Cal M}_n$ into which $D$
is embeddable, we have that by clause \scite{3.4B}(1) above:

$$
h^u_{A,B}[{\Cal M}_n] \ge h^u_{C',B \cup C'}[{\Cal M}_n]/h({\Cal M}_n)
^\varepsilon \ge h^u_{C,D_1}[{\Cal M}_n]/h[{\Cal M}_n]^{2 \varepsilon}
$$

\noindent
and by the definition of $\nonfork{}{}_{}$ also $h^u_{A,B}({\Cal M}_n) \le
h^u_{C,D_1}({\Cal M}_n)(h({\Cal M}_n))^\varepsilon$, together we have the
desired inequality.
\bigskip

In the second case $(C' = C)$ we have
$A <_s B' <_s B$, and similar inequalities using clause $(\varepsilon)$ of
\scite{3.4}(2) give the result.
\bn
\underbar{Clause $(i)$}: \newline

Note that the number of members of ${\Cal K}_\infty$ with
\smallskip

$\qquad \le |B| + |C| - |A|$ elements and, for simplicity, set of elements

$\qquad \subseteq \{1,\dotsc,|B| + |C| - |A|\}$, is finite.
\mn
So let $D_j,f_j$ for $j < j^*$ be such that
\medskip
\roster
\item "{$(a)$}"  $C \le D_j \in {\Cal K}_\infty$
\sn
\item "{$(b)$}"  $f_j$ is an embedding of $B$ into $D_j$ over $C$
\sn
\item "{$(c)$}"  $D_j = C \cup f_j(B)$
\sn
\item "{$(d)$}"  for $j_1 \ne j_2, (f_{j_2} \circ f^{-1}_{j_1} \cup
\text{ id}_C):D_{j_1} \rightarrow D_{j_2}$ is not an isomorphism
\sn
\item "{$(e)$}"  under (a)-(d), $j^*$ is maximal.
\endroster
\medskip

\noindent
So as said above, $j^* \in \Bbb N$.  Now let $M \in {\Cal K}$ be any model
to which $C$ is embeddable, say by $g_M:C \rightarrow M$ such that there are 
distinct embeddings $f_i$ of $B$ into $M$ for $i < h^d_{A,B}(M)$ extending 
$g_M \restriction A$ (remembering $A \le C$); i.e. $M$ is random enough.
So for each $i < h^d_{A,B}(M)$ for some unique $j_i < j^*$ we have
$(f_i \circ f^{-1}_{j_i}) \cup g_M$ is an isomorphism from $D_{j_i}$ onto
$M \restriction (\text{Rang}(f_i) \cup \text{ Rang}(g_M)$.
So for some $j = j_M < j^*$ we have

$$
w =: \{ i < h^d_{A,B}(M):(f_i \circ f^{-1}_j) \cup g_M \text{ embed }
D_j \text{ into } M \text{ equivalently } j_i = j_M\}
$$
\medskip

\noindent
has $\ge h^d_{A,B}[M]/j^*$ members.  So for some $j_1$ not for every random
enough ${\Cal M}_n$ do we have $j_{{\Cal M}_n} \ne j$.  Hence by 
\scite{3.5}(3) we are done.
\mn
\underbar{Clause $(j)$}: \newline

Restricting ourselves to models into which $B_2$ can be embeddable, by (g)
+ (h) we can deal with $B'_\ell \,(\ell \le 2),B'_\ell = A_\ell \cup B_0$, and
for this look at the proof of \scite{3.5}(3) above.

What about the models into which $B_2$ is not embeddable?  Look at the proof
of clause (h) above.
\mn
\underbar{Clause $(k)$}: \newline

Like the proof of $(j)$.
\mn
\underbar{Clause $(\ell)$}: \newline

Straight by computing $h_{A,D_1}({\Cal M}_n)$ (approximately) in two ways

where $D_1 = B \cup C$.
\mn
\underbar{Clause $(m)$}:  smoothness

Follows from (d) and (f). \hfill$\square_{\scite{3.4}}$
\medskip

We deal in \scite{3.6A} - \scite{3.6D} with the polynomial case;
in the general case we deal in \scite{3.6}.
\bigskip
\noindent
Concerning \scite{4.8} we add
\definition{\stag{4.9a} Definition}  We define a four place 
relation \nl
$\nonfork{}{}_{} =
\nonfork{}{}_{}[p]$ on ${\Cal K}^{+p}_\infty:\nonfork{}{}_{}(A^+,B^+,C^+,D^+)$
\underbar{iff} \nl
$A^+ \le^+ B^+ \le^+ D^+,A^+ \le^+ C^+ \le D^+,B^+ 
\cap C^+ = A^+$ and \nl
$\beta(C^+,D^+) = \beta(A^+,B^+)$.  We normally write
$\nonfork{}{}_{}$ when no confusion arises and may write this
$\nonforkin{B^+}{C^+}_{A^+}^{D^+}$.
\enddefinition
\bigskip
\noindent
Concerning \scite{4.10} we add
\demo{\stag{4.10A} Fact}  Using 
$<^{+p}_x$ for $x = b,j,t,qr$ instead $<^{{\frak K}^+}_x$ for
$x = a,i,s,pr$ respectively $({\frak K}^+,\nonfork{}{}_{}[p])$ satisfies
Definition \scite{2.6A} and satisfies $(*)$ of \scite{3.5}(2).  Moreover, it
satisfies $(*)$ of \scite{3.6A}(1) if $\bar h$ is polynomial and
$\otimes_3 + \otimes_4$ of \scite{3.6B}, in general.
\enddemo
\bigskip

\proclaim{\stag{4.12} Claim}  Assume $(*)_{\text{ap}}$ of \scite{4.8} and 
$\otimes_1 + \otimes_2$ of \scite{4.11}. \newline
1) A sufficient condition for ``$({\frak K}^+,\nonfork{}{}_{})$ is nice" is
\medskip
\roster
\item "{$\otimes_3$}"  $({\frak K},\nonfork{}{}_{})$ is nice.
\endroster
\medskip

\noindent
2) A sufficient condition for ``$({\frak K}^+,\nonfork{}{}_{})$ obeys $\bar h$
very nicely" (see Definition \scite{3.15})
\medskip
\roster
\item "{$\otimes'_3$}"  $({\frak K},\nonfork{}{}_{})$ 
obeys some $\bar h$ very nicely
\endroster
\medskip

\noindent
3) A sufficient condition for ``${\frak K}^+$ is polynomially nice"
is
\medskip
\roster
\item "{$\otimes^-_3$}"  $({\frak K},\nonfork{}{}_{})$ is polynomially nice.
\endroster
\endproclaim
\bigskip

\demo{Proof}  1) We use \scite{3.6B}(4).  Now $(*)_4$ there holds by 
\scite{4.10}(6) or use \sciteu{3.14}(2).  \newline
2) Straight. \newline
3) Define the conditions for non-polynomial cases.

$$
\align
&\le \biggl( \Pi\{p^d_{C^+}({\Cal M}_n):C^+ \subseteq B^+,C^+ \nsubseteq
A^+\} \biggr)^m \\
  &\le (h_2({\Cal M}_n))^{t_2}\|{\Cal M}_n\|^{\sum\{\beta(C^+):
C^+ \subseteq B^+,C^+ \nsubseteq A^+\}}.
\endalign
$$

\noindent
So the expected value of the number of such $\bar g$'s is

$$
\align
(h_1({\Cal M}_n),h_2({\Cal M}_n))^{t_2+t_2} &\cdot (\|M_n\|^
{\alpha(A,B)})^m \cdot
(\|{\Cal M}_n\|^{\sum\{\beta(C^+):C^+ \subseteq B^+,C^+ \subseteq A^+\}})^n \\
  &= (h_1({\Cal M}_n))^{t_1+t_2},\|{\Cal M}_n\|^{m\beta(A,B)}.
\endalign
$$

\noindent
So the expected number of such $\bar g$ for some $f$ is

$$
\align
\le (\{f:f &\text{ embeds } A \text{ into } {\Cal M}_n\}) \times (h_1
({\Cal M}_n) \times h_2({\Cal M}_n))^{t_1+t_2} \times 
\|{\Cal M}_n\|^{m,\beta(A,B)} \\
  &< \|(h_1({\Cal M}_n) \times h_2({\Cal M}_n))^{t_1+t_2} 
\times \|{\Cal M}_n\|^{m \beta(A,B) + |A|}.
\endalign
$$

\noindent
This converges to zero if $m \times \beta(A,B) > |A|$ which holds for $m$
large enough.
\enddemo
\newpage

\head {\S7 Variants of Nice} \endhead  \resetall
\bigskip

Below we define some properties of $({\frak K},\nonfork{}{}_{})$, variants
of nice, we will use semi-nice to prove elimination of quantifiers (hence 
0-1 laws), it is essentially the weakest among those discussed below.  But 
it will be natural in various contexts to verify stronger ones.
\medskip

The reader may e.g. ignore Definition \scite{2.7}(1),(2),(4),(5),(7)-(10),
\scite{2.8}(7),(8) and the version of \scite{2.8}(5),(6) with ``explicit";
in \cite{Sh:637} we present only one variant.
\bigskip

\definition{\stag{2.7} Definition}  1) We say $({\frak K},\nonfork{}{}_{})$ 
is explicitly \footnote{the ``nice" appears in \scite{2.7}(8) below} 
\underbar{nice} if we have: 
\medskip
\roster
\item "{$(*)_2$}"  for every $A <_{\text{pr}} B$ (in ${\Cal K}_\infty$) and 
$k \in \Bbb N$ for some $r = r(A,B,k) \in \Bbb N$ for every $C,D$ such that
$\nonforkin{B}{C}_{A}^{D}$ and $D = B \cup C$ we have: \newline
\smallskip
\noindent
$1 = \text{ Lim}_n \text{ Prob}_{\mu_n}\biggl($for every embedding $f$ of 
$C$ into ${\Cal M}_n$ satisfying \newline

$\qquad \qquad \qquad \qquad c \ell^r(f(A),{\Cal M}_n) \subseteq f(C)$ 
there is $g:D \rightarrow {\Cal M}_n$ \newline

$\qquad \qquad \qquad \qquad$ extending $f$ such that 
$c \ell^k(g(B),{\Cal M}_n) \subseteq g(D) \biggr)$.
\endroster
\medskip

\noindent
If $r(A,B,k) = k + |A|$ or $k$ we say $({\frak K},\nonfork{}{}_{})$ is
explicitly$^+$ nice or explicitly$^{++}$ nice respectively. \newline
[Note: to deal with e.g. successor functions we need slightly more.] \newline
\smallskip
\noindent
2)  We say $(A,A_0,B,D)$ is an almost $(k,r)$-good (quadruple) if:
\medskip

$(*)^{k,r}_{A,A_0,B,D}$  $\quad A_0 \le A \le D \in {\Cal K}_\infty$ and 
$B \le D$ and for every random enough ${\Cal M}_n$  \newline

$\qquad \qquad \qquad \,\,$ we have:
\medskip
\roster
\item "{$(**)$}"  every embedding $f:A \rightarrow {\Cal M}_n$ satisfying
\newline

\smallskip
\noindent
$\qquad c \ell^r(f(A_0),{\Cal M}_n) \subseteq f(A)$ has an extension
$g:D \rightarrow {\Cal M}_n$ satisfying \newline

$c \ell^k(g(B),{\Cal M}_n) = g(c \ell^k(B,D))$.
\endroster
\medskip

\noindent
If $r=k$ we may write $k$ instead of $(k,r)$. \newline
3) We say $(A^+,A,B,D)$ is a semi $(k,r)$-good quadruple if:
\medskip

$(*)^{k,r}_{A^+,A,B,D}$  $\quad A \le A^+ \in {\Cal K}_\infty$ and 
$A \le D,B \le D \in {\Cal K}_\infty$ \newline 

$\qquad \qquad \qquad$ and for every random enough ${\Cal M}_n$ we have:
\medskip
\roster
\item "{$(**)$}"  for every embedding $f:A^+ \rightarrow {\Cal M}_n$ 
satisfying \newline
$c \ell^r(f(A),{\Cal M}_n) \subseteq f(A^+)$ there is an extension 
$g$ of $f \restriction A$, embedding \newline
$D$ into ${\Cal M}_n$ such that \newline
$c \ell^k(g(B),{\Cal M}_n) = g(c \ell^k(B,D))$.
\endroster
\medskip

\noindent
If $r=k$ we may write $k$ instead of $(k,r)$. \newline
4) We say $(A,B,D)$ is semi $(k,r)$-good if: $A \le A^+ \in {\Cal K}_\infty
\Rightarrow (A^+,A,B,D)$ is semi-$(k,r)$-good. \newline
5)  We say ${\frak K}$ is almost nice \underbar{if} it is weakly nice and, 
for every $A \in {\Cal K}_\infty$ and $k$ for some $\ell,m,r$ we have:
\medskip
\roster
\item "{$(*)$}"  for every random enough ${\Cal M}_n$, for every 
$f:A \rightarrow {\Cal M}_n$ we have:
\smallskip
\noindent
\item "{$(**)$}"  for every $b \in {\Cal M}_n \backslash c \ell^m(f(A),
{\Cal M}_n)$ we can find $A_0,A^+,B^+$ such that:
{\roster
\itemitem{ $(\alpha)$ }  $f(A) \subseteq A_0 \subseteq A^+ \subseteq c 
\ell^m(f(A),{\Cal M}_n)$,
\itemitem{ $(\beta)$ }  $A^+ \cup \{b\} \subseteq B^+ \subseteq {\Cal M}_n$
\itemitem{ $(\gamma)$ }  $|B^+| \le \ell$
\itemitem{ $(\delta)$ }  $(A^+,A_0,A_0 \cup \{b\},B^+)$ is almost 
$(k,r)$-good
\itemitem{ $(\varepsilon)$ }  $c \ell^r(A_0,{\Cal M}_n) \subseteq A^+$
\itemitem{ $(\zeta)$ }  $c \ell^k(f(A) \cup \{b\},{\Cal M}_n) \subseteq
B^+$.
\endroster}
\endroster
\medskip

\noindent
6) We say that ${\frak K}$ is semi-nice \underbar{if} it is weakly nice and
for every $A \in {\Cal K}_\infty$ and $k$ for some $\ell,m,r$ we have:
\medskip
\roster
\item "{$(*)$}"  for every random enough ${\Cal M}_n$, and embedding
$f:A \rightarrow {\Cal M}_n$ and $b \in {\Cal M}_n$ we can find $A_0,A^+,B,D$ 
such that:
{\roster
\itemitem{ $(\alpha)$ }  $f(A) \le A_0 \le A^+ \le c \ell^m(f(A),
{\Cal M}_n)$ \newline
[note that we can have finitely many possibilities for $(\ell,m,r)$; does not
matter]
\itemitem{ $(\beta)$ }  $f(A) \cup \{b\} \subseteq B \subseteq D \subseteq 
{\Cal M}_n$
\itemitem{ $(\gamma)$ }  $|D| \le \ell$
\itemitem{ $(\delta)$ }  $(A^+,A_0,B,D)$ is semi $(k,r)$-good
\itemitem{ $(\varepsilon)$ }  $c \ell^r(A_0,{\Cal M}_n) \subseteq A^+$
\itemitem{ $(\zeta)$ }    $c \ell^k(B,{\Cal M}_n) \subseteq D$. 
\endroster}
\endroster
\medskip

\noindent
7) We say the pair $(A,B)$ is $(k,r)$-good when $A \le_s B$ and: 
if $B \le D,A \le_s D$ then $(A,B,D)$ is semi-$(k,r)$-good.  
We say $(A,B)$ is $k$-good if it is
$(k,k)$-good; good if it is $(k,k)$-good for every $k$; and $*$-good, if
for every $k$ for some $r$ it is $(k,r)$-good. \newline
8) We say ${\frak K}$ is nice if $A <_s B$ implies $(A,B)$ is $*$-good.
\newline
9) We say ${\frak K}$ is explicitly almost nice when it is weakly nice and
for every $k,\ell$ for some $r,m$, for some random enough ${\Cal M}_n$, 
(i.e. $0 < \text{ lim }
\sup \text{ Prob}_{\mu_n}$), if $A_0 \le A <_s D,A_0 \le B \le D \subseteq
{\Cal M}_n$ and $|D| \le \ell$ are such that $c \ell^m(B,{\Cal M}_n) = 
c \ell^m(B,D)$ and $c \ell^r(A_0,{\Cal M}_n) \subseteq A$ \underbar{then} 
$(A,A_0,B,D)$ is almost $(k,r)$-good. \newline
10) We say ${\frak K}$ is explicitly semi-nice when it is weakly nice and
for every $\ell$ and $k$ for some $r$ for some random enough ${\Cal M}_n$ 
we have:
\medskip
\roster
\item "{${(*)}$}"  if $A <_s D,B \le D \subseteq {\Cal M}_n,|D| \le \ell,
A^+ = c \ell^r(A,{\Cal M}_n)$ and \newline
$c \ell^k(B,{\Cal M}_n) \subseteq D$ \underbar{then} $(A^+,A,B,D)$ is 
semi $(k,r)$-good.
\endroster
\enddefinition
\bigskip

\remark{\stag{2.7A} Remark}  1) We may consider other candidates to
\scite{2.7}(7), \scite{2.7}(8). \nl
2) We may consider in the Definition of semi-good (or explicitly) 
semi-good/nice to split $k$ to two: in assumption and in conclusion. \nl
3) Also there to demand $b \notin c \ell^m(A,{\Cal M}_n)$.
\endremark
\bigskip

\fakesubhead{\stag{2.8} Fact} \endsubhead  1) In Definition \scite{2.7}(1) 
(of explicitly nice) we can replace \newline
$A <_{\text{pr}} B$ by $A <_s B$ and/or replace
$c \ell^k(g(B),{\Cal M}_n) \subseteq g(D)$ by \newline
$c \ell^k(g(B),{\Cal M}_n) = g(c \ell^k(B,D))$. 
If $({\frak K},\nonfork{}{}_{})$ is explicitly$^{++}$ nice then it is
explicitly$^+$ nice which implies it is explicitly nice. \newline
2) If $(A,A_0,B,D)$ is almost $(k,r)$-good, \underbar{then} $(A,A_0,B,D)$ is
semi $(k,r)$-good. \newline
3) $({\frak K},\nonfork{}{}_{})$ being explicitly nice implies 
${\frak K}$ is weakly nice. \newline
4) If $(A,A_0,B,D)$ is semi $(k,r)$-good and $A \cup c \ell^k(B,D) 
\subseteq D' \le D$ \underbar{then} $(A,A_0,B,D')$ is semi $(k,r)$-good.
\newline
5) If the definition of semi-nice we can demand $B = f(A) \cup \{b\}$.
\bigskip

\demo{Proof}  1) For the first phrase, clearly the new version implies the
old as \newline
$A <_{\text{pr}} B \Rightarrow A <_s B$.  So assume the old version
and let $A <_s B$, so by \scite{2.6} we can find $n$ and $A_0 <_{\text{pr}}
A_1 <_{\text{pr}} \ldots <_{\text{pr}} A_n,A_0 = A,A_n = B$.  By
\scite{2.6A}, \newline
$\nonforkin{A_{\ell +1}}{C \cup A_\ell}_{A_\ell}^{C \cup
A_{\ell+1}}$.  Define $k(\ell)$ for $\ell \subseteq n$ by downward induction
on $\ell$.  For $\ell = 0$ let $k(0) = k$, for $\ell$ let it be the
$r(r_{\ell +1},A_\ell,A_{\ell +1})$ guaranteed by \scite{2.7}(1).  Now for
random enough ${\Cal M}_n$ and embedding $f:C \rightarrow {\Cal M}_n$ such
that $c \ell^k(f(A),{\Cal M}_n) \subseteq f(C)$, we choose by induction on
$\ell \le n$ an embedding $f_\ell:C \cup A_\ell \rightarrow {\Cal M}_n$
increase with $\ell$ such that $c \ell(f_\ell(C \cup A_\ell) \subseteq
C \cup A_{\ell +1}$.  For $\ell =0$ this is given, for $\ell +1$ use the
choice of $r_\ell$.  \newline

For the second phrase, clearly $c \ell^k(g(B),{\Cal M}_n) = g(c \ell^k(B,D))$
implies \newline
$c \ell^k(g(B),{\Cal M}_n) \subseteq g(D)$, and for the other
direction remember \newline
$c \ell^k(A,N_2) \le N_1 \le N_2 \Rightarrow
c \ell^k(A,N_1) = c \ell^k(A,N_2)$. \newline
In the second sentence, first implication holds as $c \ell^n(A',B')$
increase with $m$; the second implication holds by the definition.  
\newline
2) Read the definitions. \newline
3) $({\frak K},\nonfork{}{}_{})$ explicitly nice 
$\Rightarrow {\frak K}$ is weakly nice. \newline
Let $A <_{pr} B$ and $m \in \Bbb N$ and $\varepsilon > 0$ be given.
Let $r$ be as guaranteed by \scite{2.7}(1) and let $m^*$ be such that
$A \le A' \Rightarrow (c \ell^r(A,A')) \le m^*$ (exists by \scite{2.6}(13)).
Let $\{(C_i,D_i):i < i^*\}$ list with no repetitions (up to isomorphism over
$B$) of the pairs $(C,D)$ such that the quadruple $(A,B,C,D)$ is as in 
$(*)_2$ of Definition \scite{2.7}(1) with $|D| \le |B| \times m + m^*$, 
and let $n^* \in \Bbb N$ be large enough such that
for every $n \ge n^*$ the probability of the event ${\Cal E}^i_n =$ ``for
every embedding $f:C_i \rightarrow {\Cal M}_n$ such that
$c \ell^r((A),{\Cal M}_n) \subseteq f(C_i)$ there is an embedding
$g:D_i \rightarrow {\Cal M}_n$ extending $f$" is $\ge 1 - \varepsilon/i^*$.
So for $n \ge n^*$ the probability that all the events 
${\Cal E}^0_n,\dotsc,{\Cal E}^{i^*-1}_n$ occur is $\ge 1 - \varepsilon$, 
and it suffices to prove that for such
${\Cal M}_n$ for every embedding $f:A \rightarrow {\Cal M}_n$, there are
$m$ disjoint extensions to $g:B \rightarrow {\Cal M}_n$.  Choose by
induction on $\ell$ an embedding $g_\ell:B \rightarrow {\Cal M}_n$ extending
$f$ with Rang$(g_\ell) \backslash \text{ Rang}(f)$ disjoint to
$\dsize \bigcup_{i < \ell} \text{ Rang}(g_i)$.  If we succeed to get
$g_0,\dotsc,g_{m-1}$, we are done, so assume we are stuck for some 
$\ell < m$.  By Definition \scite{2.6A}(2)(i) (existence for 
$\nonfork{}{}_{}$) we can find
$D \in {\Cal K}_\infty$ such that $C = {\Cal M}_n \restriction
\left( \dsize \bigcup_{j < \ell} \text{ Rang}(g_j) \cup c \ell^r(f(A),
{\Cal M}_n) \right) \le D$ and there is an embedding $f^+:B \rightarrow D$ 
extending $f$ such that $f^+(B) \backslash f(A)$ is disjoint to $C$ and 
$D = C \cup f^+(B)$ and $\nonforkin{f^+(B)}{C}_{f(A)}^{D}$.
Now there is $i < i^*$ such that $(D,C) \cong (D_i,C_i)$ more exactly there
is an isomorphism $h$ from $D$ onto $D_i$ such that $h(C) = C_i$ and
$f^+ = h \restriction B$ and apply ``${\Cal E}^i_n$ occurs to ${\Cal M}_n$" 
to get contradiction. \newline
4) Read definitions.  \hfill$\square_{\scite{2.8}}$
\enddemo
\bigskip

\proclaim{\stag{2.8A} Claim}  1) Assume $({\frak K},\nonfork{}{}_{})$ is an
explicitly nice 0-1 context and $r(A,B,k)$ is as in $(*)_2$ of
Definition \scite{2.7}(1).
\medskip
\roster
\item "{$(a)$}"  if $\nonforkin{(A \cup B)}{C}_{A}^{D},r = r(A,A \cup B,k)$
(of Definition \scite{2.7}(1)$(*)_2$), \underbar{then} $(C,A,B,D)$ is almost
$(k,r)$-good
\smallskip
\noindent
\item "{$(b)$}"  if $A \le A^+ \in {\Cal K}_\infty$ and $A \le_s D$ 
and $B \le D \in {\Cal K}_\infty$ \newline
and $k \le r,r(A,D,k) \le r$ where $r(A,D,k)$ is as guaranteed by
\scite{2.8}(1) then $(A^+,A,B,D)$ is semi $(k,r)$-good
\smallskip
\noindent
\item "{$(c)$}"  if $A \le_s D,B \le D \in {\Cal K}_\infty$ \underbar{then} 
$(A,B,D)$ is semi $(k,r)$-good.
\endroster
\medskip

\noindent
2) In Definition \scite{2.7}(5), in $(**)$ we can allow any $b \in 
{\Cal M}_n$. \newline
3) In Definition \scite{2.7}(6), in $(*)$ we can restrict ourselves to
$b \in {\Cal M}_n \backslash c \ell^m(f(A),{\Cal M}_n)$ and/or replace
in clause $(\gamma)$ the demand
``$|D| \le \ell$" by ``$|A^+ \cup D| \le \ell|$".
\endproclaim
\bigskip

\demo{Proof}  1)
\mr
\item "{$(a)$}"  Reread Definition \scite{2.7}(1) particularly $(*)_2$ and
Definition \scite{2.7}(2).
\sn
\item "{$(b)$}"  By clause (i) of Definition \scite{2.6A}(1) 
($=$ existence) without loss of generality 
$\nonforkin{D}{A^+}_{A}^{D^+}$ for some $D^+ \in {\Cal K}_\infty$.
By clause (a) we know $(A^+,A,D,D^+)$ is almost $(k,r)$-good.  
By \scite{2.8}(2) we get the desired conclusion. 
\sn
\item "{$(c)$}"  Left to the reader.
\ermn
2) Follows by part (2) (and the Definition \scite{2.7}(4)). \newline
3) Left to the reader.  \hfill$\square_{\scite{2.8A}}$
\enddemo
\bigskip

\proclaim{\stag{2.8B} Claim}  1) If $({\frak K},\nonfork{}{}_{})$ is
explicitly nice, \underbar{then} ${\frak K}$ is explicitly semi-nice, and
${\frak K}$ is semi-nice. \newline
2) If $({\frak K},\nonfork{}{}_{})$ is explicitly nice and \underbar{then}
${\frak K}$ is explicitly almost nice and, if in addition, it has the strong
finite basis property it is almost nice. \newline
3) If ${\frak K}$ is explicitly semi-nice, \underbar{then} ${\frak K}$ is
semi-nice. \newline
4) If ${\frak K}$ is explicitly almost nice and $({\frak K},
\nonfork{}{}_{})$ has the strong finite basis property \underbar{then}
${\frak K}$ is almost nice. \newline
5) If ${\frak K}$ is almost nice \underbar{then} ${\frak K}$ is semi-nice.
\endproclaim
\bigskip

\demo{Proof}  1) So assume that $({\frak K},\nonfork{}{}_{})$ is 
explicitly nice, so by \scite{2.8}(1), even if just \newline
$A <_s B,k \in \Bbb N$
then for some $r = r(A,B,k) \ge k$ we have $(*)_2$ of \scite{2.7}(1).
Let us prove that ${\frak K}$ is explicitly semi-nice; i.e.
Definition \scite{2.7}(10), so let $\ell$ and $k$ be given, and we should
provide $r$ as required there.  Let
$r = \text{ Max}\{r(A,B,k):A <_s D \in {\Cal K}_\infty,|D| \le \ell\}$.

We should verify $(*)$ of Definition \scite{2.7}(10), so let ${\Cal M}_n$
be random enough, \newline
$A <_s D,B \le D \subseteq {\Cal M}_n,|D| \le \ell,
A^+ = c \ell^r(A,{\Cal M}_n)$ and assume that $c \ell^k(B,{\Cal M}_n)
\subseteq D$.  We should prove that $(A^+,A,B,D)$ is semi-$(k,r)$-good.
For this it suffices to verify the assumptions of \scite{2.8A}(1)(b) but
they are obvious.  To finish the proof note that by \scite{2.8B}(3) below 
${\frak K}$ is explicitly semi-nice (see Definition \scite{2.7}(10)) implies 
that ${\frak K}$ is semi-nice. \newline
2) Similar to part (1) using this time \scite{2.8A}(1)(a) above and
\scite{2.8B}(4) below. \newline
3) So let $A \in {\Cal K}_\infty$ and $k \in \Bbb N$ be given and we have
to find $\ell,m,r$ as required in Definition \scite{2.7}(6).

Let $\bold f:\Bbb N \times \Bbb N \rightarrow \Bbb N$ be such that: 
for any $i,j$ if ${\Cal M}_n$ is random enough and $A' \subseteq {\Cal M}_n,
|A'| \le j$ then $c \ell^i(A',{\Cal M}_n)$ has at most $\bold f(i,j)$ 
elements.

Let $\ell^* = \bold f(k,|A|+1)$, now let $r^*$ be the $r$ guaranteed to exist
in Definition \scite{2.7}(10) for $k$ and $\ell$.  Now define by induction
on $i \le \ell^* + 1$ a number $m_i$ as follows: $m_0 = |A|,m_{i+1} =
\bold f(r,\ell^*) \times \dot (m_i)^{\ell^*}$ and lastly let 
$m =: m_{\ell^* +1}$.

So we have chosen $\ell,m,r$ and have to show that they are as required in
Definition \scite{2.7}(6). \newline
So let ${\Cal M}_n$ be random enough and $f:A \rightarrow {\Cal M}_n$ and
$b \in {\Cal M}_n$.  We define by induction on $i \le \ell^* + 1$ a set $A_i$
increasing with $i$ as follows:

$$
A_0 = f(A),A_{i+1} = A_i \cup \bigcup \{ c \ell^r(A',{\Cal M}_n):A' \subseteq
A_i \text{ and } |A'| \le \ell\}.
$$
\medskip

Clearly we can prove by induction on $i$ that $A_0 \le_i A_i$ and
$|A_i| \le m_i$ hence $A_i \subseteq c \ell^{m_i}(f(A),{\Cal M}_n)$.

As $c \ell^k(A_0 \cup \{b\},{\Cal M}_n) \backslash A_0$ has $\le \ell^*$ 
members necessarily for some $i < \ell + 1$ we have: $c \ell^k(A_0 \cup \{b\},
{\Cal M}_n)$ is disjoint to $A_{i+1} \backslash A_i$, and choose the minimal
such $i$.  So $c \ell^k(A_0 \cup \{b\},{\Cal M}_n) \backslash A_i$ has at
most $\ell - |A|-i$ members, so by the definition of $A_{i+1}$ there is no
$A'$, such that $A_i \cap c \ell^k(A_0 \cup \{b\},{\Cal M}_n)
<_i A' \le c \ell^k(A \cup \{b\},{\Cal M}_n)$ hence letting $A^* =
c \ell^k(A_0 \cup \{b\},{\Cal M}_n) \cap A_i$ we have $A^* \le_s c \ell^k
(A_0 \cup \{b\},{\Cal M}_n)$.  Let $D =: c \ell^k(A_0 \cup \{b\},
{\Cal M}_n)$, and
$B = A_0 \cup \{b\}$ and $A^+ = c \ell^r(A^*,{\Cal M}_n)$ hence
$A^+ \subseteq A_{i+1}$.  Now we use
Definition \scite{2.6}(10) (i.e. the choice of $r$) with $(A^+,A^*,B,D),k,
\ell,r$ here standing for $A^+,A,B,D,k,\ell,r$ there, clearly the
assumption of \scite{2.6}(10)$(*)$ holds (i.e. $A^* \le A^+ \le {\Cal M}_n,
A^* \le_s D \le {\Cal M}_n,B \le D,|D| \le \ell,f(A) = A_0 \subseteq
A^* \subseteq A^+ = c \ell^r(A^*,{\Cal M}_n) \subseteq A_{m_{\ell+1}} =
c \ell^{m_{\ell +1}}(f(A),{\Cal M}_n)$ and $r,m^*$ as required for $k,\ell$).
Hence we get $(A^+,A^*,B,D)$ is semi $(k,r)$-nice.

Now let us check requirements $(\alpha)-(\zeta)$ in $(*)$ of
Definition \scite{2.7}(6).  Now clauses $(\alpha),(\beta),(\gamma),
(\varepsilon)$ holds by the suitable choices, and clause $(\delta)$ holds 
by a previous sentence and $(\zeta)$ as $D = c \ell^k(B,{\Cal M}_n)$.
\newline
4) Similar to the proof of part (3). \newline
So assume ${\frak K}$ is explicitly almost nice and $({\frak K},
\nonfork{}{}_{})$ has the strong finite basis property and we should prove 
that ${\frak K}$ is almost nice.  So we have to check Definition 
\scite{2.7}(5),
clearly ${\frak K}$ is weakly nice, so we are given $A \in {\Cal K}_\infty$
and $k \in \Bbb N$ and we should find $\ell,m,r$ satisfying $(*)$ from
Definition \scite{2.7}(5).

Choose $r^*,m^*$ as in \scite{2.7}(9) for with $(k,|A|+1)$ here standing
for $(k,\ell)$ there.  Let $i(*)$ be such that:
\medskip
\roster
\item "{$(*)$}"  if $D \in {\Cal K}_\infty,A_i \subseteq D$ for $i \le i(*),
A_i \subseteq A_{i+1},B \subseteq D,|B| \le \bold f(m^*,|A|+1)$ then for some
$i < i(*)$ there is $B^*$ such that $A_{i+1} \cap B \subseteq B^* \subseteq
A_i$ and $\nonforkin{B^* \cup B}{c \ell^r(B^*,A_{i+1})}_{B^*}^{D}$ \newline
[why $i(*)$ exists?  as $({\frak K},\nonfork{}{}_{})$ has the finite basis
property.] \newline
Define inductively $m_i$ for $i \le i(*)$ by $m_0 = |A|$ and
$m_{i+1} = \bold f(r,m_i)$.  Lastly, let $m = m_{i(*)},r = r^*$ and $\ell =
m + f(m,|A|+1)$.
\endroster
\medskip

So let ${\Cal M}_n$ be random enough and let $f:A \rightarrow {\Cal M}_n$ and
$b \in {\Cal M}_n \backslash c \ell^m(f(A),{\Cal M}_n)$ and we should find
$A_0,A^+,B$ satisfying $(\alpha)-(\varepsilon)$ of $(**)$ of \scite{2.7}(5).
Let $B_1 = c \ell^m(A \cup \{b\},{\Cal M}_n)$.

We define by induction on $i$, a set $A^*_i$ as follows:

$$
A^*_0 = f(A),A^*_{i+1} = c \ell^r(A^*_i,{\Cal M}_n).
$$
\medskip

\noindent
Clearly $|A_i| \le m_i$.  As $({\frak K},\nonfork{}{}_{})$ has the finite
basis property (and the choice of above) we can find $i$ such that
$\nonforkin{A^*_i \cup B_1}{A^*_{i+1}}_{A^*_i}^{A^*_{i+1} \cup B_1}$.

Now choose $A_0 = A^*_i,A^+ = A^*_{i+1},B^+ = B_1 \cup A^*_{i+1}$ and let
us check clause $(\alpha) - (\varepsilon)$ of $(**)$ of \scite{2.7}(5),
now $(\alpha),(\beta)$ hold by the choices $A_0,A^+,B^+$ and clause
$(\varepsilon)$ holds by the choice of $A^*_{i+1} = A^+$, also clause
$(\gamma)$ holds as \newline
$|B^+| \le | c \ell^m(A \cup \{b\},{\Cal M}_n)| +
|A^*_{i+1}| \le \bold f(m,|A|+1) + m^*_{i+1} \le \bold f(m,|A|) + m$.
\newline
5) So assume that ${\frak K}$ is almost nice (i.e. Definition \scite{2.7}(5))
and let us check Definition \scite{2.7}(6).  Obviously, ${\frak K}$ is weakly
nice.  So let $A \in {\Cal K}_\infty$ and $k \in \Bbb N$ be given and we
should find $\ell,m,r$ such that $(*)$ of Definition \scite{2.7}(6) holds.
We just choose them as in Definition \scite{2.7}(5).  Now let ${\Cal M}_n$
be random enough and $f:A \rightarrow {\Cal M}_n$ and we should find
$A_0,A^+,B,D$ as in $(*)$ of Definition \scite{2.7}(6).  Let $A_0,A^+,B^+$ 
be as required in $(**)$ of Definition \scite{2.7}(5).  Let us choose
$B = A_0 \cup \{b\},D = B^+$, and we have to check clause 
$(\alpha)-(\delta)$ of
Definition \scite{2.7}(6).  Now they hold by the respective clauses in
Definition \scite{2.7}(5)$(**)$, but for $(\delta)$ we have to use
\scite{2.8}(2). \newline
${}$ \hfill$\square_{\scite{2.8B}}$
\enddemo
\bigskip

\remark{\stag{2.8C} Remark} Let ``${\frak K}$ is explicitly$^-$ semi-nice"
means
\medskip
\roster
\item "{$\bigotimes$}"  ${\frak K}$ is weakly nice and for every $\ell,k$
there is $m$ such that for every $\ell_1$ there is $r$ such that:
{\roster
\itemitem{ $(*)$ }  if ${\Cal M}_n$ is random enough, $A <_s D,B \le D \le
{\Cal M}_n,|A \cup B| \le \ell_0,|D| \le \ell_1$ and $A^+ = c \ell^r(A,
{\Cal M}_n),c \ell^m(B,{\Cal M}_n) \subseteq D$ \underbar{then} $(A^+,A,B,D)$
is semi $(k,r)$-good.
\endroster}
\endroster
\medskip

\noindent
Now explicitly semi-nice $\Rightarrow$ explicitly$^-$ semi-nice
$\Rightarrow$ semi-nice.
\endremark
\bigskip

\proclaim{\stag{3.6} Claim}  Assume that ${\frak K}$ is a 0-1 context obeying
$\bar h = (h^d,h^u)$ which is bounded by $h^*$, and $\nonfork{}{}_{} =
\nonfork{}{}_{\bar h,h^*}$ and assume $(*)$ of \scite{3.5}(2). \newline
1) ${\frak K}$ is weakly nice. \newline
2)  A sufficient condition for ``$({\frak K},\nonfork{}{}_{})$ is
explicitly$^+$ nice" is
\medskip
\roster
\item "{$\bigotimes_2$}"  if $A <_s B <_i D$ and $A \le_s D$, \ub{then} for
some $\varepsilon \in \Bbb R^+$ for every random enough ${\Cal M}_n$ we have
$(h[{\Cal M}_n])^{-\varepsilon} > h^u_{A,D}[{\Cal M}_n]/
h^d_{A,B}[{\Cal M}_n]$.
\endroster
\medskip

\noindent
3) A sufficient condition for $(A^+,A,B,D)$ being semi $(k,r)$-good is
\medskip
\roster
\item "{$\bigotimes_3$}"  $A \le A^+,B \le D \in {\Cal K}_\infty$ and
$A \le_s D$ and $r \ge k + |A|$ and:
{\roster
\itemitem{ $\bigoplus^{k,r}_{A,B,D}$ }  if $D \le D' = D \cup C,|C| \le k,
D' \ne D \cup c \ell^\infty(A,D')$ and \nl
$C \cap B <_i C,A \le_i A' \le_s D'$ \ub{then} for every 
$\varepsilon \in \Bbb R^+$ for every random enough
${\Cal M}_n$ we have $\varepsilon > h^u_{A,D}[{\Cal M}_n]/h^d_{A',D'}
[{\Cal M}_n]$.
\endroster}
\endroster
\medskip

\noindent
4) A sufficient condition for ${\frak K}$ being explicitly semi nice 
(see Definition \scite{2.7}(10) hence semi nice, by \scite{2.8B}(3)) is
[Saharon copied]
\medskip
\roster
\item "{$\bigotimes_4$}"  if $B \le D \in {\Cal K}_\infty,A \le_s D,k \in
\Bbb N$ and $r \ge |A| + k$ then
{\roster
\itemitem{ $\bigoplus^{k,r}_{A,B,D}$ }
if $D \le D' = D \cup C,|C| \le k,D' \ne D \cup c \ell^\infty(A,D')$ and \nl
$C \cap B <_i C,A \le_i A' \le_s D'$
\ub{then} for every $\varepsilon \in \Bbb R^+$ for every random enough
${\Cal M}_n$ we have $\varepsilon > h^u_{A,D}[{\Cal M}_n]/h^d_{A',D'}
[{\Cal M}_n]$.
\endroster}
\endroster
\medskip

\noindent
5) A sufficient condition for ${\frak K}$ being explicitly semi-nice (hence
semi-nice by \scite{2.8B}(3)) is
\medskip
\roster
\item "{$\bigotimes_5$}"  if $k,\ell \in \Bbb N,{\Cal M}_n$ is random enough,
$B \le D \le {\Cal M}_n$ and $|D| \le \ell$ and $c \ell^k(B,{\Cal M}_n) 
\subseteq D$ and $A <_s D$ and $r = k + |A|$ then
{\roster
\itemitem{ $\bigoplus^{k,r}_{A,B,D}$ }
if $D \le D' = D \cup C,|C| \le k,D' \ne D \cup c \ell^\infty(A,D')$ and \nl
$C \cap B <_i C,A \le_i A' \le_s D'$
\ub{then} for every $\varepsilon \in \Bbb R^+$ for every random enough
${\Cal M}_n$ we have $\varepsilon > h^u_{A,D}[{\Cal M}_n]/h^d_{A',D'}
[{\Cal M}_n]$.
\endroster}
\endroster
\medskip

\noindent
6) A sufficient condition for ${\frak K}$ being explicitly$^-$ semi-nice
(see \scite{2.8C}) is
\medskip
\roster
\item "{$\bigotimes_6$}"  for any $\ell,k \in \Bbb N$ there are $m,r \in
\Bbb N$ such that: if ${\Cal M}_n$ is random enough, $A \le B \le {\Cal M}_n,
|B| \le \ell,D = c \ell^m(B,{\Cal M}_n)$ and $D \le_i D'$, \nl
$D \cup C,|C| \le k,B \cap C <_i C,A \le_i A' \le_s D'$ then 
$\bigoplus^{k,r}_{A',B,D}$ holds.
\endroster
\endproclaim
\bigskip

\remark{Remark}  1) When $c \ell^k(A,M)$ has no bound by $k,|A|$, we may
consider cases \nl
$c\ell^{k,\ell+1}(A,M)$, is repeated closure under
$c\ell^k$ where $\ell$ is the length of the iteration.  What about the 
bound on the size of $c\ell^{k,\ell}(A,M)$?  Consider $log[\|{\Cal M}_n\|]$.
\nl
2) We may consider giving to $A <_i B$ possibly a negative $\alpha(A,B)$ so
measure how few such cases arise, this may help make additive formulas
meaningful.
\endremark
\bigskip

\demo{Proof}  1) By \scite{3.4A}(2). \newline
2) In Definition \scite{2.7}(1), the ``$({\frak K},\nonfork{}{}_{})$ is a
0-1 context" holds by \scite{3.5}(2) as we are assuming $(*)_2$ of 
\scite{3.5}(2).
Let $r(A,B,k) = |A|+k$.  So assume $\nonforkin{B}{C}_{A}^{D}$ and $r = |A|
+k$ and $k$ be given and \wilog $D = B \cup C$.  
Let ${\Cal M}_n$ be random enough and $f:C \rightarrow
{\Cal M}_n$.  Now we know the ``order of magnitude" of

$$
G = \{g:g \text{ is an embedding of } D \text{ into } {\Cal M}_n
\text{ extending } f\}.
$$
\medskip

\noindent
Also by the definition of $\nonfork{}{}_{}$ (in \scite{3.1}(5))

$$
G_1 = \bigl\{ g \in G:\neg \bigl( \nonfork{g(D)}{c \ell^r(f(C),{\Cal M}_n)}
_{f(C)} \bigr\}
$$
\medskip

\noindent
has smaller order of magnitude and also, by the assumption $\bigotimes_2$

$$
G_2 = \{g \in G:c \ell^r(g(D),{\Cal M}_n) \subseteq g(D) \cup 
c \ell^r(f(C),{\Cal M}_n)\}.
$$
\medskip

\noindent
Now every $g \in G \backslash G_1 \backslash G_2$ is as required, by
\scite{2.18}(1) Shmuel?. \newline
3) Similar proof. \newline
4) So let $\ell,k$ be given, we have to choose $r$ as required in
Definition \scite{2.7}(10).  Let $r = \ell + k$ and let ${\Cal M}_n$ be
random enough, so we have just to check $(*)$ from Definition \scite{2.7}(10).
So assume ${\Cal M}_n$ is random enough, $A <_s D,B \le D \le {\Cal M}_n,
|D| \le \ell,c \ell^k(B,{\Cal M}_n) \subseteq D$ and 
$A^+ = c \ell^r(A,{\Cal M}_n)$, and it suffices to prove
that the quadruple $(A^+,A,B,D)$ is semi $(k,r)$-good.  For this we use the
criterion from part (3) which holds (by assumptions above and) the
assumption $\bigotimes_4$. \nl
5), 6)  Shmuel - details.
\enddemo
\bigskip

\proclaim{\stag{3.6A} Claim}  Assume that 
${\frak K}$ obeys $\bar h$ and $\bar h$ is polynomial over $h$. \newline
1) If $A <_s B$ \underbar{then} for some $\alpha(A,B) \in \Bbb R^+$ and 
$k = k(A,B)$ (see Definition \scite{3.3}) every random enough ${\Cal M}_n$ 
satisfies:
\medskip
\roster
\item "{$(*)$}"  for every embedding $f$ of $A$ into ${\Cal M}_n$,
$$
\|{\Cal M}_n\|^{\alpha(A,B)}/h[{\Cal M}_n]^k \le \text{ nu}(f,A,B,
{\Cal M}_n) \le \|{\Cal M}_n\|^{\alpha(A,B)}h[{\Cal M}_n]^k.
$$
\endroster
\medskip

\noindent
2) In fact if $A = A_0 <_{\text{pr}} A_1 <_{\text{pr}} \cdots <_{\text{pr}} 
A_k = B$ \underbar{then} we can let $k(A,B) = \dsize \sum_{i<k} k(A_i,
A_{i+1})$ and $\alpha(A,B) = \dsize \sum_{\ell < k} \alpha(A_\ell,
A_{\ell + 1})$. \newline
So the sum $\alpha(A,B)$ does not depend on the choice of the decomposition,
i.e. of $\langle A_\ell:\ell \le k \rangle$ but only on $(A,B)$ (and
${\frak K}$). \newline
3) Some $h^*$ bounds $\bar h$ (see Definition \scite{3.1}(3)), we can
choose $h^*(n) = n$, but even can demand that it satisfies: for every 
$\varepsilon > 0$ for
every random enough ${\Cal M}_n$ we have $h^*({\Cal M}_n) < 
\|{\Cal M}_n\|^\varepsilon$. \newline
4) If $A <_s B <_s C$, \underbar{then} (see \scite{3.6A}(2) above):

$$
\alpha(A,C) = \alpha(A,B) + \alpha(B,C).
$$

\noindent
5) If $A <_s B \le_i D,A <_s D$, then $\alpha(A,B) \ge \alpha(A,D)$.
\newline
6) If $A <_s B \le D = B \cup C$ and $A \le C <_s B$ then

$$
\alpha(A,B) \ge \alpha(C,D).
$$
\endproclaim
\bigskip

\demo{Proof}  Straightforward.
\enddemo 
\bigskip
\noindent
Saharon: about almost nice?
\proclaim{\stag{3.6B} Claim}  Assume that ${\frak K}$ is a 0-1 context
obeying $\bar h$ and $\bar h$ is polynomial over $h$, so $\bar h$ is bound
by some $h^*$ (see \scite{3.4A}) and let $\nonfork{}{}_{} = 
\nonfork{}{}_{\bar h,h^*}$. \newline
1) ${\frak K}$ is weakly nice. \newline
2) $\nonforkin{B}{C}_{A}^{D}$ iff $A \le_s B \le D,A \le C \le_s D,B \cap C
=A$ and $\alpha(A,B) = \alpha(C,C \cup B)$. \newline
3) If the condition $\otimes_3$ below holds then condition $(*)_2$ of
\scite{3.5}(2) holds hence $({\frak K},\nonfork{}{}_{})$ is a 0-1 context,
where
\medskip
\roster
\item "{$\bigotimes_3$}"  if $A \le B \le D,A \le C \le D,A <_s B$ and
$\neg(C <_s D)$ or $\neg(\nonforkin{B}{C}_{A}^{D})$ and $D = B \cup C$ then
$\alpha(A,B) > \alpha(C,D)$.
\endroster
\medskip

\noindent
4) A sufficient condition for ``$({\frak K},\nonfork{}{}_{})$ is
explicitly$^+$ nice" is
\medskip
\roster
\item "{$\bigotimes_4$}"  if $A <_s B <_i D$ and $A \le_s D$, \underbar{then}
$\alpha(A,B) > \alpha(A,D)$.
\endroster
\medskip

\noindent
5) A sufficient \footnote{still it does not cover all} condition for
``${\frak K}$ is semi-nice" is
\medskip
\roster
\item "{$\bigotimes_5$}"  for every $\ell \in \Bbb N$ for some $ r \in
\Bbb N$ we have: \newline
for every random enough ${\Cal M}_n$, and every $A \le B \le {\Cal M}_n,
|B| \le \ell$, letting $D = c \ell^r(B,{\Cal M}_n),A \le_i A_0 \le_s D$ we
have
\endroster
\medskip

$\bigotimes^{k,\ell}_{A_0,B,D} \quad$  
if $D \le D_1 = D \cup C,|C| \le k,C \cap B <_i C$, \newline

$\qquad \qquad \quad A_0 \le_i A_1 <_s D,A_1$ embeddable into 
${\Cal M}_n$ over $A$ \newline

$\qquad \qquad \quad$ \underbar{then} $\alpha(A_0,D) > \alpha(A_1,D_1)$.
\medskip

\noindent
6) A sufficient condition for $\bigotimes_3$ of \scite{3.6B}(3) above, is
\medskip
\roster
\item "{$\bigotimes_6$}"  if $A <_{\text{pr}} B \le D,A \le C <_{\text{pr}} D
= B \cup C$ and \newline
$\neg(\nonforkin{B}{C}_{A}^{D})$, \underbar{then}
$\alpha(A,B) > \alpha(C,D)$.
\endroster
\endproclaim
\bigskip

\demo{Proof}
1) Should be clear. \newline
2) Should be clear. \newline
3) Read $(*)_2$ of \scite{3.5}(2). \newline
4) Read the definition. \newline
5) Easy. \newline
6) Easy. \newline
\enddemo
\bigskip

\definition{\stag{3.6C} Definition}  ${\frak K}$ is polynomially nice if
it obeys $\bar h$ which is polynomial over $h$ and satisfies
$\bigotimes_3 + \bigotimes_4$ of \scite{3.6B}.
\enddefinition
\bigskip

\definition{\stag{3.6D} Definition}  We say 
${\frak K}$ is polynomially semi-nice if it obeys $\bar h$ which is 
polynimial over $h$ and satisfies $\bigotimes_3,\bigotimes_5$ of \scite{3.6B}.
\enddefinition
\bigskip

\definition{\stag{3.15} Definition}  1) We say ${\frak K}$, a 0-1 law 
context obeys $\bar h$ very nicely it obeys $\bar h$ and 
satisfies $(*)_2$ of \scite{3.5}(2) and $(*)_2$ of \scite{3.6}(2) 
(hence it is explicitly nice), with the $\nonfork{}{}_{}$ of Definition
\scite{3.1}(5), of course. \newline
2)  We say ${\frak K}$, a 0-1 law context obeys a polynomial $\bar h$ very
nicely if it satisfies $\bigotimes_3$ of \scite{3.6B}(3) and $\bigotimes_4$
of \scite{3.6B}(4). \newline
3) By \scite{3.4A}(2), clearly ${\frak K}$ is weakly nice.  As for $(*)$ of
\scite{3.5}(2), reread the definitions.
\enddefinition
\bigskip

\definition{\stag{3.16} Definition}  We cay $({\frak K},c \ell,
\nonfork{}{}_{})$ is polynomially almost nice if:
\mr
\item "{$(*)_4$}"  for every $k,\ell$ for some $m,t,s$ for every
$A,B \le D,A \le_s D$ we have: \nl
for every random enough ${\Cal M}_n$, if \nl
$f:D \rightarrow {\Cal M}_n,A \le_s D,B \le D,c \ell^s(f(B),{\Cal M}_n)
\subseteq f(D)$, \nl
$D \cap e \ell^m(f(A),{\Cal M}_n) \equiv f(A)$ then $(A,B,D)$ is
$k$-good.
\endroster
\enddefinition
\shlhetal

------------------------------------------------------------

\newpage
    
REFERENCES.  
\bibliographystyle{lit-plain}
\bibliography{lista,listb,listx,listf,liste}

\enddocument

\bye

%% file: mathdefs.tex
\expandafter\ifx\csname mathdefs.tex\endcsname\relax
  \expandafter\gdef\csname mathdefs.tex\endcsname{}
\else \message{Hey!  Apparently you were trying to
  \string\input{mathdefs.tex} twice.   This does not make sense.} 
\errmessage{Please edit your file (probably \jobname.tex) and remove
any duplicate ``\string\input'' lines} \fi




\catcode`\X=12\catcode`\@=11

\def\n@wcount{\alloc@0\count\countdef\insc@unt}
\def\n@wwrite{\alloc@7\write\chardef\sixt@@n}
\def\n@wread{\alloc@6\read\chardef\sixt@@n}
\def\r@s@t{\relax}\def\v@idline{\par}\def\@mputate#1/{#1}
\def\l@c@l#1X{\firstpart.#1}\def\gl@b@l#1X{#1}\def\t@d@l#1X{{}}

\def\crossrefs#1{\ifx\all#1\let\tr@ce=\all\else\def\tr@ce{#1,}\fi
   \n@wwrite\cit@tionsout\openout\cit@tionsout=\jobname.cit 
   \write\cit@tionsout{\tr@ce}\expandafter\setfl@gs\tr@ce,}
\def\setfl@gs#1,{\def\@{#1}\ifx\@\empty\let\next=\relax
   \else\let\next=\setfl@gs\expandafter\xdef
   \csname#1tr@cetrue\endcsname{}\fi\next}
\def\m@ketag#1#2{\expandafter\n@wcount\csname#2tagno\endcsname
     \csname#2tagno\endcsname=0\let\tail=\all\xdef\all{\tail#2,}
   \ifx#1\l@c@l\let\tail=\r@s@t\xdef\r@s@t{\csname#2tagno\endcsname=0\tail}\fi
   \expandafter\gdef\csname#2cite\endcsname##1{\expandafter
     \ifx\csname#2tag##1\endcsname\relax?\else\csname#2tag##1\endcsname\fi
     \expandafter\ifx\csname#2tr@cetrue\endcsname\relax\else
     \write\cit@tionsout{#2tag ##1 cited on page \folio.}\fi}
   \expandafter\gdef\csname#2page\endcsname##1{\expandafter
     \ifx\csname#2page##1\endcsname\relax?\else\csname#2page##1\endcsname\fi
     \expandafter\ifx\csname#2tr@cetrue\endcsname\relax\else
     \write\cit@tionsout{#2tag ##1 cited on page \folio.}\fi}
   \expandafter\gdef\csname#2tag\endcsname##1{\expandafter
      \ifx\csname#2check##1\endcsname\relax
      \expandafter\xdef\csname#2check##1\endcsname{}%
      \else\immediate\write16{Warning: #2tag ##1 used more than once.}\fi
      \multit@g{#1}{#2}##1/X%
      \write\t@gsout{#2tag ##1 assigned number \csname#2tag##1\endcsname\space
      on page \number\count0.}%
   \csname#2tag##1\endcsname}}
\def\multit@g#1#2#3/#4X{\def\t@mp{#4}\ifx\t@mp\empty%
      \global\advance\csname#2tagno\endcsname by 1 
      \expandafter\xdef\csname#2tag#3\endcsname
      {#1\number\csname#2tagno\endcsnameX}%
   \else\expandafter\ifx\csname#2last#3\endcsname\relax
      \expandafter\n@wcount\csname#2last#3\endcsname
      \global\advance\csname#2tagno\endcsname by 1 
      \expandafter\xdef\csname#2tag#3\endcsname
      {#1\number\csname#2tagno\endcsnameX}
      \write\t@gsout{#2tag #3 assigned number \csname#2tag#3\endcsname\space
      on page \number\count0.}\fi
   \global\advance\csname#2last#3\endcsname by 1
   \def\t@mp{\expandafter\xdef\csname#2tag#3/}%
   \expandafter\t@mp\@mputate#4\endcsname
   {\csname#2tag#3\endcsname\lastpart{\csname#2last#3\endcsname}}\fi}
\def\t@gs#1{\def\all{}\m@ketag#1e\m@ketag#1s\m@ketag\t@d@l p
   \m@ketag\gl@b@l r \n@wread\t@gsin
   \openin\t@gsin=\jobname.tgs \re@der \closein\t@gsin
   \n@wwrite\t@gsout\openout\t@gsout=\jobname.tgs }
\outer\def\localtags{\t@gs\l@c@l}
\outer\def\globaltags{\t@gs\gl@b@l}
\outer\def\newlocaltag#1{\m@ketag\l@c@l{#1}}
\outer\def\newglobaltag#1{\m@ketag\gl@b@l{#1}}

\newif\ifpr@ 
\def\m@kecs #1tag #2 assigned number #3 on page #4.%
   {\expandafter\gdef\csname#1tag#2\endcsname{#3}
   \expandafter\gdef\csname#1page#2\endcsname{#4}
   \ifpr@\expandafter\xdef\csname#1check#2\endcsname{}\fi}
\def\re@der{\ifeof\t@gsin\let\next=\relax\else
   \read\t@gsin to\t@gline\ifx\t@gline\v@idline\else
   \expandafter\m@kecs \t@gline\fi\let \next=\re@der\fi\next}
\def\pretags#1{\pr@true\pret@gs#1,,}
\def\pret@gs#1,{\def\@{#1}\ifx\@\empty\let\n@xtfile=\relax
   \else\let\n@xtfile=\pret@gs \openin\t@gsin=#1.tgs \message{#1} \re@der 
   \closein\t@gsin\fi \n@xtfile}

\newcount\sectno\sectno=0\newcount\subsectno\subsectno=0
\newif\ifultr@local \def\ultralocal{\ultr@localtrue}
\def\firstpart{\number\sectno}
\def\lastpart#1{\ifcase#1 \or a\or b\or c\or d\or e\or f\or g\or h\or 
   i\or k\or l\or m\or n\or o\or p\or q\or r\or s\or t\or u\or v\or w\or 
   x\or y\or z \fi}

\def\resetall{\global\advance\sectno by 1\subsectno=0
   \gdef\firstpart{\number\sectno}\r@s@t}
\def\resetsub{\global\advance\subsectno by 1
   \gdef\firstpart{\number\sectno.\number\subsectno}\r@s@t}
\def\newsection#1\par{\resetall\vskip0pt plus.3\vsize\penalty-250
   \vskip0pt plus-.3\vsize\bigskip\bigskip
   \message{#1}\leftline{\bf#1}\nobreak\bigskip}
\def\subsection#1\par{\ifultr@local\resetsub\fi
   \vskip0pt plus.2\vsize\penalty-250\vskip0pt plus-.2\vsize
   \bigskip\smallskip\message{#1}\leftline{\bf#1}\nobreak\medskip}

\def\t@gsoff#1,{\def\@{#1}\ifx\@\empty\let\next=\relax\else\let\next=\t@gsoff
   \def\@@{p}\ifx\@\@@\else
   \expandafter\gdef\csname#1cite\endcsname##1{\zeigen{##1}}
   \expandafter\gdef\csname#1page\endcsname##1{?}
   \expandafter\gdef\csname#1tag\endcsname##1{\zeigen{##1}}\fi\fi\next}
\def\verbatimtags{\ifx\all\relax\else\expandafter\t@gsoff\all,\fi}
\def\zeigen#1{\hbox{$\langle$}#1\hbox{$\rangle$}}

\def\(#1){\edef\dot@g{\ifmmode\ifinner(\hbox{\noexpand\etag{#1}})
   \else\noexpand\eqno(\hbox{\noexpand\etag{#1}})\fi
   \else(\noexpand\ecite{#1})\fi}\dot@g}

\newif\ifbr@ck
\def\eat#1{}
\def\[#1]{\br@cktrue[\br@cket#1'X]}
\def\br@cket#1'#2X{\def\temp{#2}\ifx\temp\empty\let\next\eat
   \else\let\next\br@cket\fi
   \ifbr@ck\br@ckfalse\br@ck@t#1,X\else\br@cktrue#1\fi\next#2X}
\def\br@ck@t#1,#2X{\def\temp{#2}\ifx\temp\empty\let\neext\eat
   \else\let\neext\br@ck@t\def\temp{,}\fi
   \def\teemp{#1}\ifx\teemp\empty\else\rcite{#1}\fi\temp\neext#2X}
\def\resetbr@cket{\gdef\[##1]{[\rtag{##1}]}}
\def\references{\resetbr@cket\newsection References\par}

\newtoks\symb@ls\newtoks\s@mb@ls\newtoks\p@gelist\n@wcount\ftn@mber
    \ftn@mber=1\newif\ifftn@mbers\ftn@mbersfalse\newif\ifbyp@ge\byp@gefalse
\def\defm@rk{\ifftn@mbers\n@mberm@rk\else\symb@lm@rk\fi}
\def\n@mberm@rk{\xdef\m@rk{{\the\ftn@mber}}%
    \global\advance\ftn@mber by 1 }
\def\rot@te#1{\let\temp=#1\global#1=\expandafter\r@t@te\the\temp,X}
\def\r@t@te#1,#2X{{#2#1}\xdef\m@rk{{#1}}}
\def\b@@st#1{{$^{#1}$}}\def\str@p#1{#1}
\def\symb@lm@rk{\ifbyp@ge\rot@te\p@gelist\ifnum\expandafter\str@p\m@rk=1 
    \s@mb@ls=\symb@ls\fi\write\f@nsout{\number\count0}\fi \rot@te\s@mb@ls}
\def\byp@ge{\byp@getrue\n@wwrite\f@nsin\openin\f@nsin=\jobname.fns 
    \n@wcount\currentp@ge\currentp@ge=0\p@gelist={0}
    \re@dfns\closein\f@nsin\rot@te\p@gelist
    \n@wread\f@nsout\openout\f@nsout=\jobname.fns }
\def\m@kelist#1X#2{{#1,#2}}
\def\re@dfns{\ifeof\f@nsin\let\next=\relax\else\read\f@nsin to \f@nline
    \ifx\f@nline\v@idline\else\let\t@mplist=\p@gelist
    \ifnum\currentp@ge=\f@nline
    \global\p@gelist=\expandafter\m@kelist\the\t@mplistX0
    \else\currentp@ge=\f@nline
    \global\p@gelist=\expandafter\m@kelist\the\t@mplistX1\fi\fi
    \let\next=\re@dfns\fi\next}
\def\symbols#1{\symb@ls={#1}\s@mb@ls=\symb@ls} 
\def\bigsymbol{\textstyle}
\symbols{\bigsymbol\ast,\dagger,\ddagger,\sharp,\flat,\natural,\star}
\def\ftnumbers{\ftn@mberstrue} \def\ftsymbols{\ftn@mbersfalse}
\def\paginal{\byp@ge} \def\resetftnumbers{\ftn@mber=1}
\def\ftnote#1{\defm@rk\expandafter\expandafter\expandafter\footnote
    \expandafter\b@@st\m@rk{#1}}

\long\def\jump#1\endjump{}
\def\ssum{\mathop{\lower .1em\hbox{$\textstyle\Sigma$}}\nolimits}

\def\qed{\nobreak\kern 1em \vrule height .5em width .5em depth 0em}
\def\newneq{\hbox{\rlap{\hbox to 1\wd9{\hss$=$\hss}}\raise .1em 
   \hbox to 1\wd9{\hss$\scriptscriptstyle/$\hss}}}
\def\subsetne{\setbox9 = \hbox{$\subset$}\mathrel{\hbox{\rlap
   {\lower .4em \newneq}\raise .13em \hbox{$\subset$}}}}
\def\supsetne{\setbox9 = \hbox{$\subset$}\mathrel{\hbox{\rlap
   {\lower .4em \newneq}\raise .13em \hbox{$\supset$}}}}

\def\vbar{\mathchoice{\vrule height6.3ptdepth-.5ptwidth.8pt\kern-.8pt}
   {\vrule height6.3ptdepth-.5ptwidth.8pt\kern-.8pt}
   {\vrule height4.1ptdepth-.35ptwidth.6pt\kern-.6pt}
   {\vrule height3.1ptdepth-.25ptwidth.5pt\kern-.5pt}}
\def\f@dge{\mathchoice{}{}{\mkern.5mu}{\mkern.8mu}}
\def\b@c#1#2{{\rm \mkern#2mu\vbar\mkern-#2mu#1}}
\def\b@b#1{{\rm I\mkern-3.5mu #1}}
\def\b@a#1#2{{\rm #1\mkern-#2mu\f@dge #1}}
\def\bb#1{{\count4=`#1 \advance\count4by-64 \ifcase\count4\or\b@a A{11.5}\or
   \b@b B\or\b@c C{5}\or\b@b D\or\b@b E\or\b@b F \or\b@c G{5}\or\b@b H\or
   \b@b I\or\b@c J{3}\or\b@b K\or\b@b L \or\b@b M\or\b@b N\or\b@c O{5} \or
   \b@b P\or\b@c Q{5}\or\b@b R\or\b@a S{8}\or\b@a T{10.5}\or\b@c U{5}\or
   \b@a V{12}\or\b@a W{16.5}\or\b@a X{11}\or\b@a Y{11.7}\or\b@a Z{7.5}\fi}}

\catcode`\X=11 \catcode`\@=12

%% file: citeadd.tex

\expandafter\ifx\csname citeadd.tex\endcsname\relax
\expandafter\gdef\csname citeadd.tex\endcsname{}
\else \message{Hey!  Apparently you were trying to
\string\input{citeadd.tex} twice.   This does not make sense.} 
\errmessage{Please edit your file (probably \jobname.tex) and remove
any duplicate ``\string\input'' lines} \fi

\def\sciteu{\sciteerror{undefined}}

\def\sciteerror#1#2{{\mathortextbf{\scite{#2}}}\complainaboutcitation{#1}{#2}}
\def\mathortextbf#1{\hbox{\bf #1}}
\def\complainaboutcitation#1#2{%
\vadjust{\line{\llap{---$\!\!>$ }\qquad scite$\{$#2$\}$ #1\hfil}}}

%% file: alice2jlem.tex

\expandafter\ifx\csname alice2jlem.tex\endcsname\relax
  \expandafter\gdef\csname alice2jlem.tex\endcsname{}
\else \message{Hey!  Apparently you were trying to
\string\input{alice2jlem.tex}  twice.   This does not make sense.}
\errmessage{Please edit your file (probably \jobname.tex) and remove
any duplicate ``\string\input'' lines} \fi

\input bib4plain

\def\widestnumber#1#2{}

\def\rm{\fam0 \tenrm}

\def\fakesubhead#1\endsubhead{\bigskip\noindent{\bf#1}\par}

\input rsfs

%% file: bib4plain.tex
\expandafter\ifx\csname bib4plain.tex\endcsname\relax
  \expandafter\gdef\csname bib4plain.tex\endcsname{}
\else \message{Hey!  Apparently you were trying to \string\input
  bib4plain.tex twice.   This does not make sense.}
\errmessage{Please edit your file (probably \jobname.tex) and remove
any duplicate ``\string\input'' lines} \fi

\def\renewcommand{\newcommand}	       
\edef\cite{\the\catcode`@}%
\catcode`@ = 11
\let\@oldatcatcode = \cite
\chardef\@letter = 11
\chardef\@other = 12
%
%
%
%
\def\@innerdef#1#2{\edef#1{\expandafter\noexpand\csname #2\endcsname}}%
%
%
\@innerdef\@innernewcount{newcount}%
\@innerdef\@innernewdimen{newdimen}%
\@innerdef\@innernewif{newif}%
\@innerdef\@innernewwrite{newwrite}%
%
%
%
\def\@gobble#1{}%
%
%
%
\ifx\inputlineno\@undefined
   \let\@linenumber = \empty 
\else
   \def\@linenumber{\the\inputlineno:\space}%
\fi
%
%
%
\def\@futurenonspacelet#1{\def\cs{#1}%
   \afterassignment\@stepone\let\@nexttoken=
}%
\begingroup 
\def\\{\global\let\@stoken= }%
\\ 
\endgroup
\def\@stepone{\expandafter\futurelet\cs\@steptwo}%
\def\@steptwo{\expandafter\ifx\cs\@stoken\let\@@next=\@stepthree
   \else\let\@@next=\@nexttoken\fi \@@next}%
\def\@stepthree{\afterassignment\@stepone\let\@@next= }%
%
%
%
\def\@getoptionalarg#1{%
   \let\@optionaltemp = #1%
   \let\@optionalnext = \relax
   \@futurenonspacelet\@optionalnext\@bracketcheck
}%
%
%
\def\@bracketcheck{%
   \ifx [\@optionalnext
      \expandafter\@@getoptionalarg
   \else
      \let\@optionalarg = \empty
      \expandafter\@optionaltemp
   \fi
}%
\def\@@getoptionalarg[#1]{%
   \def\@optionalarg{#1}%
   \@optionaltemp
}%
%
%
%
\def\@nnil{\@nil}%
\def\@fornoop#1\@@#2#3{}%
\def\@for#1:=#2\do#3{%
   \edef\@fortmp{#2}%
   \ifx\@fortmp\empty \else
      \expandafter\@forloop#2,\@nil,\@nil\@@#1{#3}%
   \fi
}%
\def\@forloop#1,#2,#3\@@#4#5{\def#4{#1}\ifx #4\@nnil \else
       #5\def#4{#2}\ifx #4\@nnil \else#5\@iforloop #3\@@#4{#5}\fi\fi
}%
\def\@iforloop#1,#2\@@#3#4{\def#3{#1}\ifx #3\@nnil
       \let\@nextwhile=\@fornoop \else
      #4\relax\let\@nextwhile=\@iforloop\fi\@nextwhile#2\@@#3{#4}%
}%
%
%
%
\@innernewif\if@fileexists
\def\@testfileexistence{\@getoptionalarg\@finishtestfileexistence}%
\def\@finishtestfileexistence#1{%
   \begingroup
      \def\extension{#1}%
      \immediate\openin0 =
         \ifx\@optionalarg\empty\jobname\else\@optionalarg\fi
         \ifx\extension\empty \else .#1\fi
         \space
      \ifeof 0
         \global\@fileexistsfalse
      \else
         \global\@fileexiststrue
      \fi
      \immediate\closein0
   \endgroup
}%
%
%
%
%
\def\bibliographystyle#1{%
   \@readauxfile
   \@writeaux{\string\bibstyle{#1}}%
}%
\let\bibstyle = \@gobble
%
%
\let\bblfilebasename = \jobname
\def\bibliography#1{%
   \@readauxfile
   \@writeaux{\string\bibdata{#1}}%
   \@testfileexistence[\bblfilebasename]{bbl}%
   \if@fileexists
      \nobreak
      \@readbblfile
   \fi
}%
\let\bibdata = \@gobble
%
%
\def\nocite#1{%
   \@readauxfile
   \@writeaux{\string\citation{#1}}%
}%
\@innernewif\if@notfirstcitation
%
%
\def\cite{\@getoptionalarg\@cite}%
%
%
\def\@cite#1{%
   \let\@citenotetext = \@optionalarg
   \printcitestart
   \nocite{#1}%
   \@notfirstcitationfalse
   \@for \@citation :=#1\do
   {%
      \expandafter\@onecitation\@citation\@@
   }%
   \ifx\empty\@citenotetext\else
      \printcitenote{\@citenotetext}%
   \fi
   \printcitefinish
}%
\def\@onecitation#1\@@{%
   \if@notfirstcitation
      \printbetweencitations
   \fi
   \expandafter \ifx \csname\@citelabel{#1}\endcsname \relax
      \if@citewarning
         \message{\@linenumber Undefined citation `#1'.}%
      \fi
      \expandafter\gdef\csname\@citelabel{#1}\endcsname{%
\strut
\vadjust{\vskip-\dp\strutbox
\vbox to 0pt{\vss\parindent0cm \leftskip=\hsize 
\advance\leftskip3mm
\advance\hsize 4cm\strut\openup-4pt 
\rightskip 0cm plus 1cm minus 0.5cm ?  #1 ?\strut}}
         {\tt
            \escapechar = -1
            \nobreak\hskip0pt
            \expandafter\string\csname#1\endcsname
            \nobreak\hskip0pt
         }%
      }%
   \fi
   \csname\@citelabel{#1}\endcsname
   \@notfirstcitationtrue
}%
%
%
\def\@citelabel#1{b@#1}%
%
%
\def\@citedef#1#2{\expandafter\gdef\csname\@citelabel{#1}\endcsname{#2}}%
%
%
%
\def\@readbblfile{%
   \ifx\@itemnum\@undefined
      \@innernewcount\@itemnum
   \fi
   \begingroup
      \def\begin##1##2{%
         \setbox0 = \hbox{\biblabelcontents{##2}}%
         \biblabelwidth = \wd0
      }%
      \def\end##1{}
      %
      %
      \@itemnum = 0
      \def\bibitem{\@getoptionalarg\@bibitem}%
      \def\@bibitem{%
         \ifx\@optionalarg\empty
            \expandafter\@numberedbibitem
         \else
            \expandafter\@alphabibitem
         \fi
      }%
      \def\@alphabibitem##1{%
         \expandafter \xdef\csname\@citelabel{##1}\endcsname {\@optionalarg}%
         \ifx\biblabelprecontents\@undefined
            \let\biblabelprecontents = \relax
         \fi
         \ifx\biblabelpostcontents\@undefined
            \let\biblabelpostcontents = \hss
         \fi
         \@finishbibitem{##1}%
      }%
      \def\@numberedbibitem##1{%
         \advance\@itemnum by 1
         \expandafter \xdef\csname\@citelabel{##1}\endcsname{\number\@itemnum}%
         \ifx\biblabelprecontents\@undefined
            \let\biblabelprecontents = \hss
         \fi
         \ifx\biblabelpostcontents\@undefined
            \let\biblabelpostcontents = \relax
         \fi
         \@finishbibitem{##1}%
      }%
      \def\@finishbibitem##1{%
         \biblabelprint{\csname\@citelabel{##1}\endcsname}%
         \@writeaux{\string\@citedef{##1}{\csname\@citelabel{##1}\endcsname}}%
         \ignorespaces
      }%
      %
      %
      \let\em = \bblem
      \let\newblock = \bblnewblock
      \let\sc = \bblsc
      \frenchspacing
      \clubpenalty = 4000 \widowpenalty = 4000
      \tolerance = 10000 \hfuzz = .5pt
      \everypar = {\hangindent = \biblabelwidth
                      \advance\hangindent by \biblabelextraspace}%
      \bblrm
      \parskip = 1.5ex plus .5ex minus .5ex
      \biblabelextraspace = .5em
      \bblhook
      \input \bblfilebasename.bbl
   \endgroup
}%
%
%
\@innernewdimen\biblabelwidth
\@innernewdimen\biblabelextraspace
%
%
%
\def\biblabelprint#1{%
   \noindent
   \hbox to \biblabelwidth{%
      \biblabelprecontents
      \biblabelcontents{#1}%
      \biblabelpostcontents
   }%
   \kern\biblabelextraspace
}%
%
%
%
\def\biblabelcontents#1{{\bblrm [#1]}}%
%
%
\def\bblrm{\rm}%
%
%
\def\bblem{\it}%
%
%
\def\bblsc{\ifx\@scfont\@undefined
              \font\@scfont = cmcsc10
           \fi
           \@scfont
}%
%
%
\def\bblnewblock{\hskip .11em plus .33em minus .07em }%
%
%
\let\bblhook = \empty
%
%
%
\def\printcitestart{[}
\def\printcitefinish{]}
\def\printbetweencitations{, }
\def\printcitenote#1{, #1}
%
%
%
\let\citation = \@gobble
%
%
%
\@innernewcount\@numparams
%
%
\def\newcommand#1{%
   \def\@commandname{#1}%
   \@getoptionalarg\@continuenewcommand
}%
%
%
\def\@continuenewcommand{%
   \@numparams = \ifx\@optionalarg\empty 0\else\@optionalarg \fi \relax
   \@newcommand
}%
%
%
\def\@newcommand#1{%
   \def\@startdef{\expandafter\edef\@commandname}%
   \ifnum\@numparams=0
      \let\@paramdef = \empty
   \else
      \ifnum\@numparams>9
         \errmessage{\the\@numparams\space is too many parameters}%
      \else
         \ifnum\@numparams<0
            \errmessage{\the\@numparams\space is too few parameters}%
         \else
            \edef\@paramdef{%
               \ifcase\@numparams
                  \empty  No arguments.
               \or ####1%
               \or ####1####2%
               \or ####1####2####3%
               \or ####1####2####3####4%
               \or ####1####2####3####4####5%
               \or ####1####2####3####4####5####6%
               \or ####1####2####3####4####5####6####7%
               \or ####1####2####3####4####5####6####7####8%
               \or ####1####2####3####4####5####6####7####8####9%
               \fi
            }%
         \fi
      \fi
   \fi
   \expandafter\@startdef\@paramdef{#1}%
}%
%
%
%
%
\def\@readauxfile{%
   \if@auxfiledone \else 
      \global\@auxfiledonetrue
      \@testfileexistence{aux}%
      \if@fileexists
         \begingroup
            \endlinechar = -1
            \catcode`@ = 11
            \input \jobname.aux
         \endgroup
      \else
         \message{\@undefinedmessage}%
         \global\@citewarningfalse
      \fi
      \immediate\openout\@auxfile = \jobname.aux
   \fi
}%
%
%
\newif\if@auxfiledone
\ifx\noauxfile\@undefined \else \@auxfiledonetrue\fi
%
%
%
%
\@innernewwrite\@auxfile
\def\@writeaux#1{\ifx\noauxfile\@undefined \write\@auxfile{#1}\fi}%
%
%
%
\ifx\@undefinedmessage\@undefined
   \def\@undefinedmessage{No .aux file; I won't give you warnings about
                          undefined citations.}%
\fi
%
%
\@innernewif\if@citewarning
\ifx\noauxfile\@undefined \@citewarningtrue\fi
%
%
%
\catcode`@ = \@oldatcatcode

%% file: rsfs.tex

%
%
%

%

\font\textrsfs=rsfs10
\font\scriptrsfs=rsfs7
\font\scriptscriptrsfs=rsfs5

\newfam\rsfsfam
\textfont\rsfsfam=\textrsfs
\scriptfont\rsfsfam=\scriptrsfs
\scriptscriptfont\rsfsfam=\scriptscriptrsfs

\edef\oldcatcodeofat{\the\catcode`\@}
\catcode`\@11

\def\Cal@@#1{\noaccents@ \fam \rsfsfam #1}

\catcode`\@\oldcatcodeofat

%% file: 550.bbl
\begin{thebibliography}{ShSp 304}
\bibitem[Sh 551]{Sh:551}Saharon Shelah.
\newblock {In the random graph $G(n,p),p=n^{-a}$: if $\psi$ has probability
  $0(n^{-\varepsilon})$ for every $\varepsilon > 0$ then it has probability
  $0(e^{-n^\varepsilon})$ for some $\varepsilon > 0$}.
\newblock {\em {Annals of Pure and Applied Logic}}, {\bf 82}:97--102, 1996.
\end{thebibliography}
